\documentclass[11pt,a4paper]{article}
\usepackage[T1]{fontenc}
\usepackage[utf8]{inputenc}
\usepackage[english]{babel}
\usepackage{lmodern}
\usepackage{xcolor}
\usepackage{pict2e}
\usepackage{amsmath,amssymb,amsthm,mathtools}
\usepackage[a4paper,margin=2.5cm]{geometry}
\usepackage{enumitem,booktabs,tabularx}
\usepackage[expansion=false]{microtype}
\usepackage[hidelinks]{hyperref}
\numberwithin{equation}{section}
\newtheorem{theorem}{Theorem}[section]
\newtheorem{proposition}[theorem]{Proposition}
\newtheorem{lemma}[theorem]{Lemma}
\newtheorem{corollary}[theorem]{Corollary}
\theoremstyle{definition}
\newtheorem{definition}[theorem]{Definition}
\theoremstyle{remark}
\newtheorem{remark}[theorem]{Remark}
\newcommand{\R}{\mathbb R}
\newcommand{\N}{\mathbb N}
\newcommand{\Prob}{\mathbb P}
\newcommand{\E}{\mathbb E}
\newcommand{\dd}{\,\mathrm d}
\newcommand{\Ind}{\mathbf 1}
\newcommand{\phin}{\varphi_{\mathrm{in}}}
\newcommand{\phout}{\varphi_{\mathrm{out}}}
\newcommand{\pbar}{\bar\varphi}
\newcommand{\ov}{\operatorname{ov}}
\newcommand{\sgn}{\operatorname{sgn}}
\newcommand{\F}{\mathcal F}

\newcommand{\TV}{d_{\mathrm{TV}}}
\newcommand{\GBG}{\texttt{GBG}}
\DeclareMathOperator{\dist}{dist}
\setlist{itemsep=2pt,topsep=4pt}

 \newcommand{\Q}{{\bf Q}}  
 
 \newcommand{\M}{   \mathsf M}  
\title{\textbf{Distinguishability and community recovery\\
in the weight-dependent random connection model}}
\author{Rémi Bernard, Raphaël Lachièze-Rey}
\date{}
\begin{document}
\maketitle

 \begin{abstract}
Weight-dependent random connection models (WDRCMs) are a family of geometric random graphs that exhibit power-law degrees, cliques behaviours, and several other features of real social networks. We study community detection in the WDRCM with two communities {\it preferential attachment} kernel
 with parameter $\gamma\in(0,1)$. In the assortative case, we show that, without knowing profile connection functions and preferential attachment marks, (i) the model is distinguishable in every dimension: a test  almost surely decides correctly whether communities are present once the finite observation window is large enough; (ii) for nonincreasing, integrable profiles in dimension $d\geq2$, strong recovery is possible, meaning that an arbitrarily large fraction of nodes can be classified correctly as the intensity $ \lambda $ of nodes tends to infinity. The second result uses a non parametric reconstruction procedure that does not know the ages either. These results apply in the more realistic asymptotic  {\it sparse} regime, i.e.  when the connection rule does not change as the window grows; result (i) requires neither monotonicity nor integrability of the profiles. The overlap deficit is at most $\exp[-\lambda/(\log\lambda)^{2d+3}]$ for this agnostic procedure. If on the other hand   ages and parameters are known, a similar method gives a deficit that decreases exponentially  in $e^{ -c\lambda }$.    These rates concern the high-intensity limit, taken after the volume limit at fixed intensity. We discuss the results in relation to Euclidean geometric models, in particular \cite{abbe}, adapting and extending some of their proof techniques to heavy-tailed degrees and explaining the role of percolation in recovery.
\end{abstract}

 {\bf Keywords:} random geometric graphs, community detection, random connection models, preferential attachment, heavy tail 

\section{Introduction}

 {Community detection in a large graph aims to recover groups of nodes that tend to be more connected between them than with other nodes. This question arises in both social and biological networks, but a high local density of connections does not by itself reveal a community: a set of vertices may have many internal connections because its members occupy the same region of space even in the absence of communities. The challenge is to identify a community  phenomenon within a cluster beyond these effects of proximity and varying degrees of heterogeneity. These features, central to network science \cite{newman}, motivate the joint study of community structure, geometry and node heterogeneity.}

The stochastic block model (SBM) is the canonical model to perform community detection, and has been intensively studied, with optimal theoretical results identifying information theoretic detection thresholds \cite{abbe-sbm}. The light tail (Poisson) typical degree distribution and non-geometric nature makes it unadapted to large social networks, where power-law degree, cliques behaviour, and many other characteristic features have been consistently observed  \cite{barabasi,newman}. In the WDRCMs, nodes are randomly placed in a Euclidean framework, and connection probabilities depend on both distance and node weights. For the preferential attachment kernel considered here, these weights are represented by marks, called  {\it ages}, with older nodes being more attractive, reflecting the fact that nodes degrees are very heterogeneous in real networks. This mechanism produces heavy-tailed degrees while retaining geometric effects on triangles and percolation \cite{gracar2019,gracar2021}, and it is mathematically more tractable than the classical preferential attachment mechanism, where the attractivity of a node is random, built along the graph evolution. It therefore provides a setting in which detection can be studied
beyond a description based only on groups.

In the Euclidean setting, Abbe, Baccelli and Sankararaman \cite{abbe} developed  a recovery procedure based on counting triangles, i.e. common neighbours of given pairs of points,   in the so-called Good-Bad-Grid algorithm, denoted by \GBG,  that we explain below.

 We tackle here the two following problems: \begin{itemize}
\item The question of  {\it distinguishability } is that of determining if there are indeed communities or not, or said differently if the probability to connect to a member of the other communities can be different from that of connecting to a member of the same communities. In statistical terms, we therefore give an algorithm able to detect in finite time if the  {\it null hypothesis} of the absence of communities is true or not.
\item At the other end of the spectrum, {\it strong recovery} aims to show that, when there are indeed communities with different behaviours, then the proportion of nodes misclassified by our algorithm goes to $ 0$  as the intensity of nodes in a fixed Euclidean substrate goes to $ \infty $.
\end{itemize}  
Let us briefly present the proof strategies.
For distinguishability, a study of the truncated degrees  allows to find the unique candidate profile without communities that is compatible with the observations,
and an excess of triangles {
rules it out under the assortativity and non-degeneracy assumption (H).
}   The first ingredient is {
the injectivity of the map $\varphi\mapsto q_\varphi$, where
$q_\varphi(R)=\int_{(0,1)^2}\varphi(Rg(u,v))\dd u\dd v$, it represents the connection probability
averaged over all possible ages
}.  For a typical vertex, the set of statistics of mean number of neighbours within the ball of radius \(R\) for all \(R\) completely determine this functional and so the profile \(\varphi\) almost everywhere. The only profile without communities who fit this set of statistics is the mean profile  \(\bar\varphi=(\varphi_{\mathrm{in}}+\varphi_{\mathrm{out}})/2\), however, the two-community model has a strictly larger mean triangle count for at least one radius, which {
rules out this remaining candidate under (H).
} When the profiles are unknown and only measurable,  the procedure is said to be  {\it agnostic}. The alternative must be separated from the whole null family. The joint density of the connection arguments makes the moments continuous in the weak-star topology, giving the required compactness. A uniform fluctuation bound turns this separation into a finite-window test. Appendix~\ref{app:finite-references} gives an implementation based on comparison with a finite list of references, without an efficiency guarantee at moderate sample sizes.\\

For the recovery problem, the profiles are  supposed to be integrable and nonincreasing, which means that nodes connect preferentially to closer nodes.   The \GBG~algorithm first compares pairs of nodes through their common neighbours, then combines these decisions on a spatial grid: at high intensity, cells with reliable classifications form a large component on which labels can be recovered. Here, this first, partial reconstruction provides the starting point for a voting procedure that progressively assigns the remaining labels. Heterogeneous attractiveness, however, requires us to control ages in the local tests and then to work without this information when ages are not observed, but only estimated.
We first slice the graph into  overlapping age slices, and check that \GBG\ applies on each slice. Then the classifications obtained  are aligned thanks to the  only possibility consistent over the overlaps.  A voting procedure then proceeds towards smaller ages; comparison with an oracle controls the effect of errors on each new slice, without uncontrolled accumulation. The universal procedure replaces ages by ranks of calibration degrees and known thresholds by quartiles of counts, choosing its scales from the estimated intensity in the same way for every admissible model.

The global dependence of the empirical ranks and thresholds is the main difficulty in this adaptation. We construct a local auxiliary field that guarantees all the decisions whenever the global calibrations are accurate enough, then quantify the passage from reliable cells to recovery. The total amplification of errors in voting remains bounded, while the initialisation introduces a logarithmic loss in the exponent; scales chosen for the model parameters recover an exponential bound, even without ages. Section~\ref{sec:distinguability} treats distinguishability, Section~\ref{sec:recovery} constructs the recovery procedures, and Section~\ref{sec:rate} establishes their rates. Proofs specific to ranks and quartiles are given in Appendix~\ref{app:universal}.

\section{Model and main results}\label{sec:model}
\subsection{Marked process and connection kernel}
Let $d\geq1$ and let $\bar \eta$ be a Poisson point process on $\R^d\times(0,1]$ with intensity $\lambda\dd x\dd u$, where $\lambda>0$.
{We refer to \cite[Chapter~5]{lastpenrose} for background on marked Poisson processes.} We write its points as $X=(x,u)$, with position $x$ and age $u$.   The  {\it ground process} is denoted by $ \eta= \{x: (x,u)\in  \bar \eta \}$.  Denote $ \theta _{ z}(x,u) = (x + z,u)$ spatial shift operators, for $ z\in \mathbb{R}^{ d}.$ The processes $ \eta $ and $  \bar \eta $ are space invariant, or stationary, i.e. $ \theta _{ z}\eta \stackrel{(d)}{=}  \eta , \theta _{ z}  \bar \eta \stackrel{(d)}{=}   \bar \eta .  $

Conditionally on $\bar \eta$, nodes receive  i.i.d.  uniform labels $\sigma_X\in\{-1,+1\},   {X = (x,u)\in  \bar \eta }$. Set
\begin{equation}\label{eq:kernel}
g(u,v)=(u\wedge v)^\gamma(u\vee v)^{1-\gamma},\qquad
\kappa_{XY}=\frac{g(u,v)|x-y|^d}{\beta},
\quad \gamma\in(0,1),\quad\beta>0.
\end{equation}
Conditionally on $  \bar \eta$ and $\sigma$, edges are independent and satisfy
\begin{equation}\label{eq:connection}
\Prob(X\sim Y\mid\eta,\sigma)=
\begin{cases}
\phin(\kappa_{XY}),&\sigma_X=\sigma_Y,\\
\phout(\kappa_{XY}),&\sigma_X\neq\sigma_Y.
\end{cases}
\end{equation}

 Since the connection rule does not depend on the spatial coordinate, the law of the resulting geometric graph is also invariant under spatial translations.

The goal of community recovery is to recover each node's label, up to a global exchange of the two communities.
The profiles take values in $[0,1]$. We denote the edge configuration by $E$ and the marked graph by $G=(\bar\eta,E)$; its community labels are not observed. It will be useful to distinguish
\begin{equation}\label{eq:signal}
\psi=\phin-\phout,\qquad \pbar=\frac{\phin+\phout}{2},
\end{equation}
 the  {\it signal} $\psi$ from the average profile $\pbar$. We always reserve $g$ for the attachment kernel.

\begin{remark}[Normalising the density parameter]\label{rem:beta}
The dilation $x\mapsto c^{1/d}x$ transforms $(\lambda,\beta)$ into $(\lambda/c,c\beta)$ without changing the connections. Taking $c=\beta^{-1}$ therefore allows us to set $\beta=1$ in calculations with known parameters, with effective intensity $\lambda\beta$. The universal procedure uses the observed coordinates and does not apply this dilation. To analyse it with an unknown fixed $\beta$, we instead absorb $\beta$ into the profiles $t\mapsto\varphi(t/\beta)$; the constants and the intensity threshold may then depend on $\beta$.

\end{remark}

In the recovery part, the profiles are nonincreasing, as in the usual definition of the model. We use the following assumptions:
\begin{description}[leftmargin=1.5cm,style=nextline]
\item[(A0)] The profiles do not depend on the observation volume and
\begin{equation}\label{eq:sparse}
\int_{\R^d}\int_0^1\int_0^1
(\phin+\phout)\bigl(g(u,v)|x|^d\bigr)
\dd u\dd v\dd x<\infty.
\end{equation}
\item[(H)] The signal is assortative and non-degenerate: $0\leq\phout\leq\phin\leq1$ and \[\operatorname{Leb}\{t>0:\psi(t)>0\}>0.\]
\end{description}
Distinguishability only requires measurable profiles and (H), with no monotonicity or assumption (A0), in every dimension $d\geq1$. Recovery results instead assume nonincreasing profiles, (A0), (H) and $d\geq2$.

\subsection{Observations and inference problems}
We consider the growing windows  
\[
B_n=[-n^{1/d}/2,n^{1/d}/2)^d,\qquad
\bar \eta _{ n}=\bar \eta\cap(B_n\times(0,1]),\qquad N_n=\#\bar \eta _{ n}.
\]
Thus $\operatorname{Leb}(B_n)=|B_n|=n$ and $N_n\sim\operatorname{Pois}(\lambda n)$. The data consist of the induced graph $(\bar \eta _{ n},E_n)$ when ages are observed, or $(\eta _{ n},E_n)$ when ages are hidden; all positions are known, including those of isolated nodes. We first study recovery with known ages and parameters, then without ages, profiles or structural parameters, except for the dimension. The universal procedure estimates the intensity itself; {
Section~\ref{sec:estimation} gives the intensity estimator and its
role in the ideal test, while Section~\ref{sec:rate} discusses the
separate question of retaining the observed-age exponential rate
when intensity is unknown.
} To make algorithmic choices unambiguous, we order grid cells by increasing $\ell_\infty$ norm of their indices, using lexicographic order to break ties; for every finite selection of nodes, we use lexicographic order of their positions.

\subsection{Weak and strong recovery}

\begin{definition}[Overlap]\label{def:overlap}
For an estimator $\hat\sigma^{(n)}$ based on the available data, define
\[
\ov_n(\sigma,\hat\sigma^{(n)})=
\frac1{N_n}\left|\sum_{X\in\bar \eta _{ n}}\sigma_X\hat\sigma_X^{(n)}\right|,
\]
with the convention $1$ if $N_n=0$.  {A  {\it trivial}} estimator $ \hat \sigma  ^{ (n)}$, i.e.  independent of labels and connections, gives an average asymptotic recovery $  \ov_{ n} \xrightarrow[ n\to \infty ]{a.s.}0$ by the strong law of large numbers (the $ \sigma _{ X}$ are  i.i.d.) . 
The absolute value accounts for the undetermined global sign. If $p$ is the fraction of nodes with correct labels in a given orientation, the overlap is $|2p-1|$. We keep $\gamma$ for the attachment exponent and use $q$ or $\Gamma$ for an overlap level.

\end{definition}
\begin{definition}[Weak recovery]
    We say that weak recovery is possible if there exist deterministic  or randomised estimators and a constant $q>0$ such that $\Prob(\ov_n\geq q)\to1$ as $n\to\infty$. 
\end{definition}

\begin{definition}[Strong recovery in the high-intensity limit]\label{def:strong}
We say that strong recovery is possible in the high-intensity limit if, for every $\varepsilon>0$, there exists $\lambda(\varepsilon)<\infty$ such that, for each $\lambda\geq\lambda(\varepsilon)$, deterministic or randomised estimators can be chosen with $\Prob(\ov_n\geq1-\varepsilon)\to1$. Here the limit $n\to\infty$ is taken at fixed $\lambda$.
\end{definition}

 Strong recovery requires neither that the overlap tends to one at a fixed finite intensity, nor exact recovery, which would require $\Prob(\ov_n=1)\to1$. 

Another framework seen in the literature is when, at fixed intensity $ \lambda $, the degree of the typical node goes to infinity with $ n$, typically at a logarithmic speed. We rather consider that the model is fixed, and evaluate the effectiveness of the method as $ \lambda $ increases.  Note that exact recovery is typically impossible at fixed $ \lambda $ as there will always be isolated points, impossible to classify with certainty.

\subsection{Distinguishability}

Distinguishability asks whether the graph has community structure or whether $\phin=\phout$ almost everywhere  i.e. if (H) is violated. We seek decisions that are almost surely correct once the window is large enough, without knowing the profiles or ages; the index from which they become correct is random and is not supplied by the test. We denote by $\mathbf{P}_\varphi $ the null model, without communities, with profile $\phin=\phout=\varphi$, and by $\Q_{\phin,\phout}$  the two-community model.  Note that distinguishability is strictly weaker than weak recovery.

Define
$\F=\{f\in L^\infty(0,\infty):0\leq f\leq1\text{ a.e.}\}$.
Under the null model $\mathbf{P}_\varphi $, edges are conditionally independent given $\bar \eta$, with probabilities $\varphi(\kappa_{XY})$; under $   \Q =\Q_{\phin,\phout}$, they follow \eqref{eq:connection}. The tests do not observe communities or ages, but know $(\lambda,d,\gamma)$;  {
we give an  estimation method for $ \lambda $ and $ \gamma $  in  Propositions~\ref{prop:intensity-estimation}  and Remark~\ref{rem:gamma-estimate},
 under the additional assumption (A0)
}. Lemma~\ref{lem:density} ensures that changing a profile on a null set does not change these laws.

\begin{definition}[Distinguishability]
The two-community model is distinguishable from $\mathbf{P}_\varphi $ if there exists a sequence of tests based on the observations that, under equal prior probabilities, identifies the generating model with success probability at least $1/2+q-o(1)$ for some $q>0$. A consistent test, whose two error probabilities tend to zero, satisfies this definition with success probability tending to one.
\end{definition}  {
This definition only requires a positive asymptotic advantage
over random guessing and allows a nonvanishing error probability.
Our result is stronger: under each fixed admissible model,
the tests are almost surely correct for all sufficiently
large observation windows. The index from which all decisions
are correct is random and is not supplied by the procedure. In an ergodic setting, it is expected that the success probability is $ 0$ or $ 1$; the previous general definition also applies to some particular non-ergodic settings; e.g. in the planted partition model.
}

\subsection{Main results}
The following test does not know the ages or profiles; it only assumes known structural parameters and requires neither integrability nor monotonicity of the profiles.

\begin{theorem}[A universal distinguishability test]\label{th:disc}
Fix $\lambda>0$, $d\geq1$ and $\gamma\in(0,1)$.
Let $\mathcal O_n=\sigma(\eta _{ n},E_n)$ be the
sigma-field generated by positions and induced edges in $B_n$.
There exists a sequence of $\mathcal O_n$-measurable tests
$\tau_n\in\{0,1\}$,
 whose definition does not depend on $ \lambda $, the connection
profiles $\phin$, $\phout$ or $\varphi$,
such that, for every $\varphi\in\F$ and every $\phin,\phout\in\F$
satisfying (H),
\begin{gather*}
\mathbf{P}_\varphi
\bigl(\exists n_0\in\N,\ \forall n\geq n_0,\ \tau_n=0\bigr)=1
,\\
\Q_{\varphi_{\mathrm{in}},\varphi_{\mathrm{out}}}
\bigl(\exists n_0\in\N,\ \forall n\geq n_0,\ \tau_n=1\bigr)=1
  .
\end{gather*}
There is also a measurable event $\mathcal A$ of the infinite
spatial graph such that $\mathbf{P}_\varphi(\mathcal A)=1$ for every $\varphi$
and $\Q_{\varphi_{\mathrm{in}},\varphi_{\mathrm{out}}}(\mathcal A)=0$ for every admissible alternative. 
The laws are therefore mutually singular:
$\mathbf{P}_\varphi\perp \Q_{\varphi_{\mathrm{in}},\varphi_{\mathrm{out}}}$.  Furthermore, this separating event does not depend on profile functions.
\end{theorem}  
Section~\ref{sec:test} defines this test through the distance to the null family. Appendix~\ref{app:finite-references} gives an explicit procedure comparing observed moments with a finite family of reference profiles. Their moments are approximated by simulation with controlled error, while preserving consistency for each fixed alternative. The required volume depends on the signal; no uniform guarantee is claimed as the signal tends to zero.

{
For recovery, we first consider a procedure that uses the ages
and model parameters to choose its settings, and then an
agnostic procedure that uses only the observed positions,
edges and dimension, together with independent auxiliary
randomness, without requiring knowledge of the ages,
connection profiles or model parameters. The same algorithm
applies to all models satisfying the assumptions stated below,
although the intensity above which its recovery guarantee
holds may vary from one model to another, and this threshold
is not an input to the algorithm.

The auxiliary randomness is independent of the graph,
its ages and its community labels, and is used in our
construction to assign independent auxiliary colours
to the nodes, thereby splitting them into groups used
at different stages of the procedure.
}
\begin{theorem}[Quantitative recovery and agnostic recovery]
\label{th:main-rates}
Suppose $d\geq2$, the profiles $\phin,\phout$ are nonincreasing,
and (A0), (H) hold, with fixed $\gamma\in(0,1)$ and $\beta>0$.
Then the following statements hold.
\begin{enumerate}
\item \textit{ [Known ages and model parameters.]}
There exist constants $C,c,\lambda_0>0$ and a sequence of
deterministic label estimators $\hat\sigma^{(n)}$  {{measurable  with respect to}} $(\bar \eta _{ n},E_n)$
  such that, for each
$\lambda\geq\lambda_0$,
\[
\Prob\!\left(
\ov_n(\sigma,\hat\sigma^{(n)})\geq1-Ce^{-c\lambda}
\right)\longrightarrow1
\qquad(n\to\infty).
\]
The constants may depend on the fixed profiles and on
$d,\gamma,\beta$, but not on $\lambda$ or $n$.

\item \textit{[Agnostic recovery.]}
There exists a sequence of randomised label estimators
$\hat\sigma^{(n)}$
{defined from $ d$, the observed positions, induced edges 
and independent auxiliary randomness, whose definition does not depend
on the ages, connection profiles or $\gamma,\beta,\lambda$.
For every fixed pair of nonincreasing profiles $\phin,\phout$
satisfying (A0) and (H), and every fixed $\gamma\in(0,1)$ and
$\beta>0$, there exists $\lambda_0>1$ such}
that, for each $\lambda\geq\lambda_0$,
\[
\Prob\!\left(
\ov_n(\sigma,\hat\sigma^{(n)})
\geq1-\exp\!\left[
-\frac{\lambda}{(\log\lambda)^{2d+3}}
\right]
\right)\longrightarrow1
\qquad(n\to\infty).
\]    
The probability includes the auxiliary randomisation.
The threshold $\lambda_0$ may depend on the fixed profiles
and on $d,\gamma,\beta$; no uniformity over these choices
is claimed.
\end{enumerate}
In both statements, the limit $n\to\infty$ is taken at fixed
intensity $\lambda$.
\end{theorem}
The agnostic procedure is described in
Definition~\ref{def:hidden-algorithm} and uses an independent auxiliary coloration, the role of which will be clarified during its construction. For each fixed model, an average argument then allows us to fix this randomness while preserving the asymptotic guarantees, {
which also establishes deterministic existence for each fixed model.
}
 When parameters are known but ages are hidden, Corollary~\ref{cor:tuned-hidden} recovers the exponential bound by choosing fixed scales in the same rank-based construction. Limits are always taken at fixed intensity, and we make no claim that the rates are optimal.

  {
\subsection{Estimation of \texorpdfstring{$\lambda$ and $\gamma$}{lambda and gamma}}\label{sec:estimation}
The positions of all nodes, including isolated ones, are observed,
so the intensity can be estimated directly from their number. Intensity estimation is a super classical topic of spatial statistics, but we remind it here for completeness.
Estimating the attachment exponent also uses the graph, through its
mean degree and its proportion of isolated nodes. Those results are solely for the interested reader, they are not used for our main results.

\begin{proposition}[Estimation of the intensity]\label{prop:intensity-estimation}
Under each fixed model with intensity $\lambda>0$, the estimator
\begin{equation}\label{eq:lambda-estimate}
\widehat\lambda_n=\frac{N_n}{n}
\end{equation}
satisfies
\[
\E\widehat\lambda_n=\lambda,\qquad
\operatorname{Var}(\widehat\lambda_n)=\frac{\lambda}{n},\qquad
\widehat\lambda_n\longrightarrow\lambda\quad\text{almost surely}.
\]
For the windows $W_L=[-L,L)^d$, let
$\widehat\lambda_L^{\rm win}=N(W_L)/|W_L|$; then
\[
\Prob\!\left(
\left|\widehat\lambda_L^{\rm win}/\lambda-1\right|>L^{-1/8}
\right)\leq C_{\lambda,d}L^{1/2-2d},\qquad L\geq1.
\]
These estimators use neither ages nor community labels.
\end{proposition}
\begin{proof}
Since $N_n\sim\operatorname{Pois}(\lambda n)$, its first two
moments give the expectation and variance, while the strong law for
the Poisson counts in expanding cubes gives almost-sure convergence.
For $V=|W_L|=(2L)^d$, the fourth centred Poisson moment gives
\[
\E\left|\frac{N(W_L)}{\lambda V}-1\right|^4
=(\lambda V)^{-3}+3(\lambda V)^{-2}.
\]
Markov's inequality at threshold $L^{-1/8}$ proves the last bound,
which is summable in $L$ for every $d\geq1$.
\end{proof}

\begin{corollary}[Unknown intensity in the ideal distinguishability test]
\label{cor:unknown-intensity-test}
The almost-sure eventual correctness of the ideal test in
Theorem~\ref{th:disc} remains valid when $\lambda$ is unknown:
in the rule~\eqref{eq:test}, replace the intensity in both the
normalisations and the null reference moments by
$\widehat\lambda_L^{\rm win}$, and use threshold $L^{-1/16}$.
The parameters $d$ and $\gamma$ remain known, and the test returns
zero if $\widehat\lambda_L^{\rm win}=0$.
\end{corollary}
Its proof is given after the construction of the test in
Section~\ref{sec:test}; this statement concerns the ideal distance
to the null family, whereas Appendix~\ref{app:finite-references}
establishes its finite-comparison implementation for known
structural parameters.

For the degree calculation below, we use $\beta=1$, absorbing a fixed
$\beta$ into the profiles as in Remark~\ref{rem:beta}. The Palm law
places a typical node at the origin, leaving the other marked Poisson
points unchanged; its full definition is recalled in
Section~\ref{sec:local}.

\begin{lemma}[Mean degree]\label{lem:degree}
Let $\omega_d$ be the volume of the Euclidean unit ball and set $I_{\pbar}=\omega_d\int_0^\infty\pbar(t)\dd t$. Under (A0), the degree of a vertex of age $u$, under the Palm law of the infinite graph, is Poisson with mean
\begin{equation}\label{eq:mean-degree}
\lambda I_{\pbar}\Theta(u),\qquad
\Theta(u)=\frac1{1-\gamma}+\frac{u^{-\gamma}-1}{\gamma}.
\end{equation}
This law does not depend on its community. In particular, $\E^0\deg(0)=2\lambda I_{\pbar}/(1-\gamma)<\infty$.
\end{lemma}
\begin{proof}
Under Palm, independent thinning of potential neighbours, after averaging over their signs, gives the Poisson law with profile $\pbar$. The change of variables $y=g(u,v)^{1/d}x$ gives the mean $\lambda I_{\pbar}\int_0^1g(u,v)^{-1}\dd v$. The integrals over $(0,u)$ and $(u,1)$ equal $1/(1-\gamma)$ and $(u^{-\gamma}-1)/\gamma$, respectively, yielding \eqref{eq:mean-degree}. Finally, $\int_0^1\Theta(u)\dd u=2/(1-\gamma)$. This calculation also shows that (A0) is equivalent to integrability of $\phin+\phout$ on $(0,\infty)$, for fixed $\gamma$.
\end{proof}

\begin{remark}[Estimating the attachment exponent $ \gamma $]\label{rem:gamma-estimate}
Under (A0), with a nonzero mean profile, the mean degree $m>0$ and the proportion $p_0$ of isolated vertices also contain information about $\gamma$, even when ages and profiles are unknown. The mixed Poisson degree law, established in Lemma~\ref{lem:degree}, gives
\[p_0=F_m(\gamma),\qquad F_m(\theta)=\int_0^1e^{-m r_\theta(u)}\dd u,\qquad r_\theta(u)=\frac12+\frac{1-\theta}{2\theta}(u^{-\theta}-1).\]
{
For fixed $m>0$, $F_m$ extends continuously to $[0,1]$
by dominated convergence. To verify strict monotonicity,
write, for $0<\theta<1$,
\[
R_\theta(a):=\int_0^a r_\theta(u)\dd u
=a+\frac{a^{1-\theta}-a}{2\theta}
=a+\frac a2\int_0^{-\log a}e^{\theta s}\dd s.
\]
Hence $\partial_\theta R_\theta(a)>0$ for $0<a<1$,
whereas $r_\theta'(u)<0$. Differentiating under the integral
and integrating by parts therefore gives
\[
F_m'(\theta)
=
m\int_0^1
\frac{\dd}{\dd u}\bigl(e^{-m r_\theta(u)}\bigr)
\,\partial_\theta R_\theta(u)\dd u
>0.
\]
Differentiation is justified locally in $\theta\in(0,1)$
by an integrable bound of the form
$C(1+u^{-b})(1+|\log u|)$ with $b<1$.
The boundary terms vanish since
$\partial_\theta R_\theta(a)\to0$ as $a\downarrow0$
and $\partial_\theta R_\theta(1)=0$.
Thus $F_m$ is strictly increasing on $[0,1]$, with limiting values 
\[
F_m(0)=\frac{e^{-m/2}}{1+m/2},
\qquad
F_m(1)=e^{-m/2}.
\]
}  Thus $\widehat m_n=2|E_n|/N_n$ and $\widehat p_{0,n}=N_n^{-1}\#\{X:\deg_{G_n}(X)=0\}$ yield an estimator of $\gamma$ by inversion of $F_{\widehat m_n}$, after projecting $\widehat p_{0,n}$ onto its range; if $N_n=0$ or $\widehat m_n=0$, set the estimator to zero. Under (A0), ergodicity and control of edges lost at the boundary give convergence of these two statistics, and then of the estimator. We mention this possibility without developing its proof here. It suggests extending the distinguishability test with finite-radius isolation proportions to handle an unknown $\gamma$; such an extension would require checking the corresponding separation and uniform approximations. We claim no rate for estimating $\gamma$. Consistency alone does not justify substituting this estimator into a test with a decreasing threshold, and Theorem~\ref{th:disc} remains stated for known $\gamma$. Universal reconstruction does not need this estimate.
\end{remark}
}

\section{Distinguishability}\label{sec:distinguability}

 Distinguishing between the two models is done by observing the triangles count statistics. Note that observing only the degrees through the marginal law would not allow to distinguish both models.

In the calculations below, we take $\beta=1$ as in Remark~\ref{rem:beta}. With unknown profiles, this convention replaces them by $t\mapsto\varphi(t/\beta)$ and requires no transformation of the observations.
\subsection{Local counts and ergodic limits}\label{sec:local}
Let $\eta=\sum_{(x,u)\in\bar\eta}\delta_x$
be the spatial projection that forgets ages, and let
$\eta_n=\eta|_{B_n}$ be its restriction to the observation
window $B_n$. For distinct nodes $x,y\in\eta$, write
$x\sim y$ if they are connected, and set the degree
and triangle counts 
\begin{align}
d_R(x,G)&=\sum_{y\in{\eta}\setminus\{x\}}\Ind_{\{|y-x|\leq R\}}\Ind_{\{x\sim y\}},\label{eq:degree-score}\\
t_R(x,G)&=\sum_{\substack{y,z\in{\eta}\setminus\{x\}\\y\neq z}}\Ind_{\{|y-x|\leq R,\ |z-x|\leq R\}}\Ind_{\{x\sim y,\ x\sim z,\ y\sim z\}}.\label{eq:triangle-score}
\end{align}
Triangles are counted through ordered pairs, hence twice each.
 Introduce the Palm versions $ \eta ^{ \circ } = \eta \cup \{0\},  \bar \eta ^{ \circ } =  \bar \eta \cup \{(0,U)\}$ where $U $ is an independent  variable uniform on $ (0,1)$. Let $ \sigma _{ 0}\in \{-1,1\}$ an independent Rademacher label and $   \mathbf P^{ \circ }_{ \varphi},   \Q ^{ \circ }_{ \varphi _{ in},\varphi _{ out}} $ the laws of the graphs where connections between $ 0$ and points of $ \eta  $ (resp. $  \bar \eta  $) are drawn independently according to the probability  \eqref{eq:connection}.
We use the generic notation $   \mathsf M$ is either $  {\bf P}_{ \varphi},  {\bf Q}$ and $   \mathsf M^{ \circ }$ is defined accordingly, called  {\it Palm measure}. Let $ G^{ \circ }$ the corresponding random graph.

By Slivnyak's theorem for the marked Poisson process, this
amounts to adding a node with an independent uniform age at
the origin, while the remaining nodes retain their original
Poisson law. 
Let $\omega_d=\operatorname{Leb}_d(B(0,1))$ be the volume of
the Euclidean unit ball. 
The number $ K_{ R}$ of points in  in $B(0,R)$ is Poisson with parameter $ \lambda \omega _{ d}R^{ d}$, therefore the truncated degree and number of triangles are trivially bounded 
$$d_R(0,G)\leq K_R\qquad t_R(0,G)\leq K_R(K_R-1),$$
where triangles are counted through ordered pairs of distinct
neighbours. The Poisson factorial moments give
\begin{equation}\label{eq:palm-bounds}
\begin{aligned}
0\leq D_R(\mathsf M)
&:=\E_{\mathsf M}^0[d_R(0,G)]
\leq\lambda\omega_dR^d,\\
0\leq\Delta_R(\mathsf M)
&:=\E_{\mathsf M}^0[t_R(0,G)]
\leq(\lambda\omega_dR^d)^2.
\end{aligned}
\end{equation}
The expectations are taken under the Palm law of the model
being considered, whereas the upper bounds hold for both
models. These local expectations remain finite even when
the expected total degree is infinite.

For a bounded region $C$, write
$N(C)=\#(\eta _{ n}\cap C)$ and define
\begin{equation}\label{eq:estimators}
\widehat D_{R,C}
=\frac{\sum_{x\in\eta\cap C}d_R(x,G)}{N(C)},
\qquad
\widehat\Delta_{R,C}
=\frac{\sum_{x\in\eta\cap C}t_R(x,G)}{N(C)},
\end{equation}
with value zero when $N(C)=0$. These averages include all
nodes, including isolated ones. They are observable in a
window $B$ whenever $C+B(0,R)\subset B$, since every node
needed to compute a score centred in $C$ then lies in the
observed window.

To state the ergodic convergence of these empirical averages,
fix $R>0$ and, for integers $L>\lceil R\rceil$, use the
observation windows
\[
W_L=[-L,L)^d=B_{(2L)^d}
\]
and the interior cubes
\[
C_{L,R}
=[-L+\lceil R\rceil,L-\lceil R\rceil)^d.
\]
We will let $L\to\infty$ while keeping $R$ fixed.
Removing this boundary strip ensures that the radius-$R$
neighbourhood of every node in $C_{L,R}$ is contained in
$W_L$. Nodes in the boundary strip may still serve as
neighbours; only the nodes over which we compute the
averages are restricted to the interior cube.

Recall that in both models, the law of the geometric graph is invariant under spatial translations. Call  {\it translation-invariant event} any $ \{0,1\}$-valued $ \M$-measurable functional $ \mathcal{B} $ such that $ \M$-a.s., $  \mathcal{B}(\theta _{ z}\M) =\mathcal{B}(\M),z\in \mathbb{R}^{ d} . $ Denote by $ \theta _{ z}\mathcal{B} = \mathcal{B}(\theta _{ z}\,.\,).$

\begin{proposition}[Ergodic limits of the counts]\label{prop:ergodic}

Let $\mathsf M$ denote either a fixed null-model law
$\mathbf P_\varphi$ or a fixed two-community law $\Q $.  Then $   \mathsf M$
 is ergodic, meaning that every translation invariant event has probability zero or one.
Consequently, for every fixed $R>0$,
\[
\widehat D_{R,C_{L,R}}\longrightarrow D_R(\mathsf M),
\qquad
\widehat\Delta_{R,C_{L,R}}\longrightarrow\Delta_R(\mathsf M)
\quad \mathsf M\text{-a.s.}
\]
as $L\to\infty$ through integers $L>\lceil R\rceil$.
For any fixed countable family of radii, these convergences
hold simultaneously on an event of $\mathsf M$-probability one.
\end{proposition}

\begin{proof} Call induced graph on some region $ A\subset \mathbb{R}^{ d}$ the set of nodes $ x\in A$ and edges connecting only those nodes $ x\in A$.
Induced
graphs in disjoint bounded regions are independent: they
use independent restrictions of the marked Poisson process
and independent variables for their internal edges.
{ 
Let $\mathcal F_m$ be the sigma-field generated by the graph
induced in $[-m,m)^d$. These sigma-fields are increasing
and generate the sigma-field of the infinite spatial graph,
since every node and every edge is eventually contained
in one of these boxes.
For any measurable event $\mathcal A$, the martingale
convergence theorem therefore gives
\[
p_m:=\E_{\mathsf M}
[\mathbf 1_{\mathcal A}\mid\mathcal F_m]
\longrightarrow\mathbf 1_{\mathcal A}
\quad\text{in }L^1.
\]
Thus the events $\mathcal A_m=\{p_m>1/2\}\in\mathcal F_m$
satisfy
\[
\mathsf M(\mathcal A\triangle\mathcal A_m)
\leq 2\E_{\mathsf M}
\bigl[|p_m-\mathbf 1_{\mathcal A}|\bigr]
\longrightarrow0,
\]
which proves the required approximation by events
depending on bounded regions.
}
After a sufficiently large translation, two such local
events are independent. Approximation therefore gives
\[
\mathsf M(\mathcal A\cap\theta_z^{-1}\mathcal B)
\longrightarrow
\mathsf M(\mathcal A)\mathsf M(\mathcal B)
\qquad (|z|\to\infty).
\]
For a translation-invariant event $\mathcal A$, taking
$\mathcal B=\mathcal A$ yields
$\mathsf M(\mathcal A)=\mathsf M(\mathcal A)^2$,
which proves ergodicity.

Fix $R>0$ and let $\xi(x,G)$ denote either $d_R(x,G)$
or $t_R(x,G)$, viewed as a score assigned to the node $x$.
Set
\[
M_\xi(C)=
\sum_{x\in\eta\cap C}\xi(x,G).
\]
Both scores are measurable, have finite Palm expectations
by \eqref{eq:palm-bounds}, and satisfy
$\xi(x+z,\theta_zG)=\xi(x,G)$: translating the graph
does not change the score of a node.
Thus $M_\xi$ is a stationary random measure with
\[
\E_{\mathsf M}[M_\xi(C)]
=\lambda|C|\,\E_{\mathsf M}^0[\xi(0,G)]
\]
by {the refined Campbell formula}.
{
Here $\E_{\mathsf M}^0$ denotes expectation under the Palm
law introduced above, so that $\E_{\mathsf M}^0[\xi(0,G)]$
is the mean score of a typical node, placed at the origin.
The formula therefore expresses the expected total score
in $C$ as the mean number of nodes in $C$ multiplied
by the mean score of a typical node.
}  {
We apply the almost-sure spatial ergodic theorem for
stationary random measures \cite{nguyen1979ergodic}.
For a stationary ergodic random measure with finite intensity,
this theorem states that its mass per unit volume in expanding
cubes converges almost surely to its intensity.
Here, $M_\xi$ and $N$ are stationary and ergodic, as
translation-covariant functions of the underlying ergodic graph,
and their intensities are respectively
$\lambda\E_{\mathsf M}^0[\xi(0,G)]$ and $\lambda$.
Since $C_{L,R}=[-t,t)^d$ with
$t=L-\lceil R\rceil\to\infty$, it follows that
}
\[
\frac{M_\xi(C_{L,R})}{|C_{L,R}|}
\longrightarrow\lambda\E_{\mathsf M}^0[\xi(0,G)],
\qquad
\frac{N(C_{L,R})}{|C_{L,R}|}
\longrightarrow\lambda>0
\quad\mathsf M\text{-a.s.}
\]
Taking the ratio proves the asserted convergence for
each of the two scores.
Finally, let $\mathcal E_R$ be the event on which both
empirical averages converge to their respective Palm
expectations at radius $R$. We have
$\mathsf M(\mathcal E_R)=1$, so for any countable family
$(R_q)_{q\geq1}$,
\[
\mathsf M\!\left(\bigcap_{q\geq1}\mathcal E_{R_q}\right)=1.
\]
This gives simultaneous convergence for that family.

\end{proof}
For $\varphi\in\F$, let $\mathbf P_\varphi$ denote the law
of the model without communities with connection profile $\varphi$,
that is, $\phin=\phout=\varphi$. Conditional on positions and ages,
edges are independent and a pair of distinct nodes $X,Y$ is connected
with probability $\varphi(\kappa_{XY})$.
{Community labels play no role under this law, since
connection probabilities do not depend on whether the labels agree.}

For $\varphi\in\F$, write    
\[q_\varphi(t)=\int_{(0,1)^2}\varphi\bigl(tg(u,v)\bigr)\dd u\dd v,\quad t>0.\]
The positions $0,x,y$ and ages $u_0,u_x,u_y$ give the effective distances $\kappa_{0x}=|x|^dg(u_0,u_x)$, $\kappa_{0y}=|y|^dg(u_0,u_y)$ and $\kappa_{xy}=|x-y|^dg(u_x,u_y)$. The Campbell--Mecke formula gives  
\begin{align}
D_R(\mathbf{ P}_\varphi)&=\lambda\int_{B(0,R)}q_\varphi(|x|^d)\dd x,\label{eq:edge-moment}\\
\Delta_R(\mathbf{ P}_\varphi)&=\lambda^2\int_{B(0,R)^2}\left(
\int_{(0,1)^3}\varphi(\kappa_{  {\bf 0}x})\varphi(\kappa_{  {\bf 0}Y})\varphi(\kappa_{XY})\dd\mathbf u
\right)\dd x\dd y.\label{eq:triangle-moment}
\end{align}
 where $ X = (x,u_{ x}),Y = (y,u_{ y}),  {\bf 0} = (0,u_{ 0}),  {\bf u} = (u_{ 0},u_{ x},u_{ y})$.
The factorial measure $\lambda^2\dd x\dd y$ corresponds
to ordered pairs; Tonelli's theorem and
\eqref{eq:palm-bounds} justify the integrations.
Since the factors in \eqref{eq:triangle-moment} share ages,
their joint integral does not generally reduce to a product
of $q_\varphi$ terms. Triangles may already be present under the
null model, so their presence alone does not distinguish
the models: we must identify an excess relative to the
null model compatible with the truncated-degree statistics.
\subsection{Known profiles: identification and distinguishability}\label{sec:identification}
The set statistic of truncated degrees for every range $R$ identify the only null candidate $\varphi=\pbar$; we then show that its triangles cannot match those of the alternative with communities. This step provides both a simple test for known profiles and the separation needed for a common test.
\begin{lemma}[Injectivity without monotonicity]\label{lem:injectivity}\label{cor:degree-identification}
For $\varphi\in\F$, the function $q_\varphi$ is continuous on $(0,\infty)$.  For $ \varphi_{ 1},\varphi_{ 2}\in \F$, if $q_{\varphi_1}=q_{\varphi_2}$ almost everywhere, then $\varphi_1=\varphi_2$ almost everywhere. This conclusion also holds if $D_R(   \mathbf P_{ \varphi_{ 1}})=D_R(\mathbf{ P}_{\varphi_2})$ for every rational radius $R>0$.
\end{lemma}
\begin{proof}
 By symmetry, the contribution of the domains $u<v$ and $ v<u$ each  account for half the value of $ q_{ \varphi}(t)$.  The change of variable $u = \alpha v, z=\alpha^\gamma v$ gives 
\begin{align}\notag   \displaystyle q_{ \varphi}(t) =     & 2 \int_{(0,1)^{ 2}}\varphi(tu^{ \gamma }v^{ 1-\gamma }){\bf 1}\{ u<v \}dudv=   \int_0^1\varphi(tz)\ell_\gamma(z)\dd z\label{eq:age-density}\\
\ell_\gamma(z)=&2z\int_{z^{1/\gamma}}^1\alpha^{-2\gamma}\dd\alpha,\quad 0<z<1
\end{align}
For $F(x)=\varphi(e^x)$ and $E_b(t)=e^{-bt}\Ind_{\{t>0\}}$, we obtain
\begin{equation}\label{eq:convolution}
q_\varphi(e^x)=(F*k_\gamma)(x),\qquad k_\gamma(t)=\ell_\gamma(e^{-t})e^{-t}\Ind_{\{t>0\}}=\frac2\gamma E_2*E_{1/\gamma}(t).
\end{equation}
Since $k_\gamma\in L^1$, continuity of its translations in $L^1$ and the inequality $|(F*k_\gamma)(x+h)-(F*k_\gamma)(x)|\leq\|F\|_\infty\|k_\gamma(\cdot+h)-k_\gamma\|_1$ give continuity of $q_\varphi$.

Apply this representation to $\varphi=\varphi_1-\varphi_2$. For every bounded function $F$, $F*E_b=e^{-bx}\int_{-\infty}^xF(t)e^{bt}\dd t$ is locally absolutely continuous and satisfies $(F*E_b)'+b(F*E_b)=F$. Setting $V=F*E_{1/\gamma}$ and $Z=V*E_2$, the equality $q_\varphi=0$ successively implies $Z=0$, $V=Z'+2Z=0$ and $F=V'+\gamma^{-1}V=0$ almost everywhere. Thus $\varphi_1=\varphi_2$ almost everywhere.

Finally, $D_R(\mathbf{ P}_\varphi)=\lambda\omega_d\int_0^{R^d}q_\varphi(t)\dd t$. Equality of these primitives at rational radii extends by continuity to all radii; differentiation and injectivity then conclude the proof.
\end{proof}
Truncated degrees therefore identify the profile. It remains to rule out the candidate $\pbar$.

 Let us lay out the strategy heuristically.
Conditionally on positions and ages, write the edge probability as $\pbar(\kappa_{ij})+\tfrac12\psi(\kappa_{ij})\sigma_i\sigma_j$. In the product of the three probabilities, terms with one or two signal factors have expectation zero because a sign appears to an odd power; the product of all three edge signs, however, equals one. The triangle probability, averaged over communities, is therefore   
\begin{equation}\label{eq:triangle-identity}
\prod_{ij\in\{12,13,23\}}\pbar(\kappa_{ij})+\frac18\prod_{ij\in\{12,13,23\}}\psi(\kappa_{ij}).
\end{equation}
To make this excess strictly positive without regularity of the profiles, we use the additional freedom provided by ages.
\begin{lemma}[Regularisation by ages]\label{lem:density}  Let $ x_{ 1},x_{ 2},x_{ 3}\in \mathbb{R}^{ d}$ distinct, with i.i.d. marks $u_1,u_2,u_3$ {uniform on $(0,1)$} and the corresponding effective distances $ \kappa _{ ij}.$ The triple $(\kappa_{12},\kappa_{13},\kappa_{23})$ has a law absolutely continuous with respect to Lebesgue measure on $(0,\infty)^3$.
\end{lemma}
These variables are connection arguments, not ages: at fixed positions, $0<\kappa_{ij}\leq|x_i-x_j|^d$. Having a density does not mean that it is positive everywhere in the orthant.
\begin{proof}
The density of a single argument follows from \eqref{eq:age-density} by rescaling. For the triple, separate the six age orderings. On $u_1<u_2<u_3$, the logarithms of the arguments are, up to spatial constants, the image of the logarithms of the ages under
\[M_\gamma=\begin{pmatrix}\gamma&1-\gamma&0\\\gamma&0&1-\gamma\\0&\gamma&1-\gamma\end{pmatrix},\qquad \det M_\gamma=-\gamma(1-\gamma)\neq0.\]
This transformation, followed by exponentiation, is a diffeomorphism on the region considered. The other five orderings follow by permutation, and equal ages have probability zero, giving the joint density.
\end{proof}
{
The joint density in this lemma is obtained by averaging
over the ages alone, while keeping the three distinct
positions fixed. This remains true in dimension one,
 although the ordinary Euclidean distances between three
points on a line satisfy an additive relation: the largest
distance is the sum of the other two. The age-dependent
factors remove this constraint for the effective distances,
whose joint law is absolutely continuous in three dimensions.
The specific feature used here is therefore regularisation
by ages at fixed positions; in other models, a joint density
may instead arise from averaging over random positions.
}
\begin{lemma}[Strict positivity in every dimension]\label{lem:positive}
If $\psi:(0,\infty)\to[0,1]$ is measurable and $|\{\psi>0\}|>0$, there exists an integer $R\geq1$ such that
\begin{equation}\label{eq:positive-integral}
\int_{B(0,R)^2}\int_{(0,1)^3}\psi(\kappa_{0x})\psi(\kappa_{0y})\psi(\kappa_{xy})\dd\mathbf u\dd x\dd y>0.
\end{equation}
\end{lemma}
\begin{proof}
Choose $b>0$ and a bounded measurable set $A\subset(0,\infty)$ of positive measure on which $\psi\geq b$. The Lebesgue density theorem allows us to choose $c>0$ such that
\[\lim_{\delta\downarrow0}\frac{|A\cap(c-\delta,c+\delta)|}{2\delta}=1.\]
We construct an interior configuration whose three arguments equal this density point $c$. Set
\[\rho=2^{\gamma/d},\quad v=\frac12,\quad t=v(1+\rho^{-1})^{-d/(1-\gamma)},\quad s=t/2,\quad \ell=\left(\frac{c}{t^\gamma v^{1-\gamma}}\right)^{1/d}.\]
Take ages $(u_0,u_x,u_y)=(s,v,t)$ and positions $x=\rho\ell e_1$, $y=(\rho+1)\ell e_1$. Then $0<s<t<v<1$ and
\[s^\gamma v^{1-\gamma}(\rho\ell)^d=t^\gamma v^{1-\gamma}\ell^d=s^\gamma t^{1-\gamma}((\rho+1)\ell)^d=c.\]
The configuration is colinear, so it also works in dimension one.

Fix an integer $R>(\rho+1)\ell$. Near this configuration, ages remain strictly ordered and positions distinct, so the map $\Phi(x,y,\mathbf u)=(x,y,\boldsymbol\kappa)$ is $C^1$ with nonzero Jacobian by Lemma~\ref{lem:density}. The inverse function theorem provides a neighbourhood $V$ on which $\Phi$ is a diffeomorphism onto an open set containing $(x_0,y_0,c,c,c)$. Choose $O\times J^3$ relatively compact in this open set, with $O\subset B(0,R)^2$ and $J$ around $c$; by continuity, $|\det D\Phi^{-1}|\geq b_0>0$ on this product. The integral is therefore at least $b_0|O|b^3|A\cap J|^3>0$. Only the positive measure of $A\cap J$ is used: $\psi$ need not be positive on a whole interval.
\end{proof}
\begin{proposition}[Reference profile and excess triangles]\label{prop:excess}
Under the two-community model,
\begin{align}
D_R(\Q)&=D_R(  {\bf P}_{\pbar}),\quad R>0,\label{eq:edge-alternative}\\
\Delta_R(\Q)-\Delta_R({\bf P}_{\pbar})&=\frac{\lambda^2}{8}\int_{B(0,R)^2}\int_{(0,1)^3}\psi(\kappa_{0x})\psi(\kappa_{0y})\psi(\kappa_{xy})\dd\mathbf u\dd x\dd y.\label{eq:triangle-excess}
\end{align}
The excess is nonnegative for every $R$ and, under (H), strictly positive for some integer $R$.
\end{proposition}
\begin{proof}
The marginal probability of an edge is $\pbar(\kappa_{XY})$, which gives equality of the truncated degrees. For triangles, integrating \eqref{eq:triangle-identity} gives \eqref{eq:triangle-excess}; the coefficient $\lambda^2/8$ accounts for ordered pairs and averaging over signs, while Lemma~\ref{lem:positive} ensures strict positivity at a suitable integer radius.
\end{proof}
\begin{remark}[Test with known profiles]\label{rem:known-test}
When $\varphi,\phin,\phout$ are known,
if {$\varphi$ and $\pbar$ differ on a set of positive
Lebesgue measure},
a truncated degree at a rational radius separates their means
by Lemma~\ref{cor:degree-identification};
otherwise a triangle score at an integer radius works
by Proposition~\ref{prop:excess}.  For this score $\xi$, write $m_0=\E_{\mathbf{P}_\varphi }^0\xi$, $m_1=\E_Q^0\xi$, and let $\widehat m_L$ be its empirical average. The test
\begin{equation}\label{eq:known-test}
T_L^{\rm known}=\Ind_{\{(m_1-m_0)(\widehat m_L-(m_0+m_1)/2)>0\}}
\end{equation}
chooses the closer theoretical mean; ergodic convergence ensures almost surely correct decisions from some index onwards. If, instead, $\psi=0$ almost everywhere, then $Q=P_{\pbar}$ and no distinction is possible. To remove knowledge of the profiles, we still need to separate the alternative from the entire null family at once.
\end{remark}
\subsection{Unknown profiles: from individual separation to a common test}\label{sec:unknown}

{
Without known profiles, we cannot choose in advance the score that
separates the models, so we collect the moments in a vector
$\mu(\mathsf M)$, where $\mathsf M$ denotes either a null law
$\mathbf P_\varphi$ or the fixed alternative
$\Q=\Q_{\phin,\phout}$. Our aim is to show that the alternative's
moment vector $\mu(\Q)$ has positive distance from the whole null family
$\mathcal C_0=\{\mu(\mathbf P_\varphi):\varphi\in\F\}$.
Identification by degrees and the triangle excess rule out each null
profile individually; compactness will turn this exclusion into separation.
}

The description through countably many observables recalls subgraph densities \cite{lovaszszegedy,lovasz} and nonparametric network inference \cite{bickelchen}, but we use Banach--Alaoglu weak-star compactness, rather than graphon compactness in the cut metric. The joint density of the connection arguments is essential: if all three arguments were the same uniform random variable, the moment would be $\int \varphi^3$, which is not weak-star continuous, as oscillating indicators of mean $1/2$ show. Marginal densities alone would therefore not suffice. Truncation of observables, standard in nonparametric testing \cite{ingstersuslina}, will be justified by a fluctuation bound independent of the profiles.

Keep degree scores at rational radii, which are needed for identification by continuity, and triangle scores at integer radii, since a positive triangle excess persists as the radius increases. Enumerate these scores as $(\xi_j,R_j)$ with $R_j\leq j$: at step $k$, add $t_k$ and the $d_{p/q}$ with $\max(p,q)=k$. Set
\[c_j=\begin{cases}1+\lambda\omega_dR_j^d,&\xi_j=d_{R_j},\\1+(\lambda\omega_dR_j^d)^2,&\xi_j=t_{R_j},\end{cases}\quad \mu_j(\mathsf M)=\E_{\mathsf M}^0(\xi_j/c_j)\in[0,1].\]
Write $\mu(\mathsf M)=(\mu_j(\mathsf M))_{j\geq1}$ and $\mu(\varphi)=\mu(\mathbf{ P}_\varphi)$. Normalisation bounds the expectations, rather than the individual scores, and controls the tail of the moment sequence. On $[0,1]^{\N^*}$, we use
\begin{equation}\label{eq:profile-metric}
d_*(a,b)=\sum_{j\geq1}2^{-j}|a_j-b_j|,\qquad d_{*,J}(a,b)=\sum_{j=1}^J2^{-j}|a_j-b_j|.
\end{equation}
Weak-star convergence $\varphi_k\to \varphi$ in $L^\infty(0,\infty)$ means here that
\[\int \varphi_k(t)a(t)\dd t\longrightarrow\int \varphi(t)a(t)\dd t\qquad\text{for every }a\in L^1(0,\infty).\]
This does not require pointwise convergence. Banach--Alaoglu will give compactness of $\F$, while the next lemma handles the nonlinearity of triangle moments.
\begin{lemma}[Continuity of integrated products]\label{lem:tensor}
If $\varphi_k,\varphi\in\F$ and $\varphi_k\to \varphi$ in the weak-star topology $\sigma(L^\infty,L^1)$, then, for every $p\geq1$ and every $W\in L^1((0,\infty)^p)$,
\[\int W(\mathbf t)\prod_{i=1}^p\varphi_k(t_i)\dd\mathbf t\longrightarrow\int W(\mathbf t)\prod_{i=1}^p\varphi(t_i)\dd\mathbf t.\]
\end{lemma}
\begin{proof}
For a product $W(\mathbf t)=\prod_iw_i(t_i)$, with $w_i\in L^1$, the integral factorises and weak-star convergence applies to each factor. The same holds for finite sums of products, which are dense in $L^1$ by approximation with simple functions on bounded rectangles. Replacing $W$ by such an approximation $W_0$ changes the difference of the integrals by at most $2\|W-W_0\|_1$, uniformly in $k$, since the products are bounded by one.
\end{proof}
\begin{proposition}[Compactness and strictly positive separation]\label{prop:compactness}\label{prop:separation}
The set $\mathcal C_0$ is compact for $d_*$. For each fixed admissible alternative,
\begin{equation}\label{eq:separation-distance}
\Delta_{\Q}^*:=\dist_{d_*}(\mu(\Q),\mathcal C_0)>0.
\end{equation}
\end{proposition}
\begin{proof}
The set $\F$ is weak-star closed in the unit ball of $L^\infty$: the inequalities $\varphi\geq0$ and $1-\varphi\geq0$ are preserved by integration against nonnegative $L^1$ functions. The Banach--Alaoglu theorem \cite{alaoglu} therefore gives compactness. Separability of $L^1(0,\infty)$ is not needed for this compactness theorem; it makes the weak-star topology metrisable on this bounded family, so that continuity can be checked using sequences.

By Lemma~\ref{lem:density}, each normalised edge coordinate has the form $\mu_j(\varphi)=\int W_j(t)\varphi(t)\dd t$, with $W_j\geq0$ in $L^1$. Each triangle coordinate similarly has the form
\[\mu_j(\varphi)=\int_{(0,\infty)^3}W_j(\mathbf t)\varphi(t_1)\varphi(t_2)\varphi(t_3)\dd\mathbf t,\qquad W_j\geq0,\quad W_j\in L^1.\]
To see this, integrate the conditional densities over the positions: coincident positions form a negligible set, and Tonelli's theorem applies. In both cases, $\|W_j\|_1=\mu_j(1)\leq1$. Lemma~\ref{lem:tensor} then gives continuity of each coordinate. The tail of the series \eqref{eq:profile-metric} is uniformly bounded by $2^{-J}$, so $\varphi\mapsto\mu(\varphi)$ is continuous for $d_*$, and its image is compact.

If a null profile $\varphi$ were to satisfy $\mu(\Q)=\mu(\varphi)$, all truncated-degree coordinates would agree. Identity~\eqref{eq:edge-alternative} and Lemma~\ref{cor:degree-identification} would force $\varphi=\pbar$ almost everywhere. This sole candidate is excluded by Proposition~\ref{prop:excess}, which provides an integer radius at which the triangle moments differ, a contradiction. The distance from $\mu(\Q)$ to the compact set $\mathcal C_0$ is therefore strictly positive.
\end{proof}
This separation can be detected with finitely many coordinates, since
\[\inf_{\varphi\in\F}d_{*,J}(\mu(\Q),\mu(\varphi))\geq\Delta_{\Q}^*-2^{-J}.\]
As soon as $2^{-J}<\Delta_{\Q}^*/2$, the first $J$ moments separate every null profile from the alternative. The test will increase $J$ without knowing this distance.
\paragraph{Separation on the infinite graph.}
Let $\pi(t)=\min(1,\max(0,t))$ and set    
\begin{equation}\label{eq:normalized-estimator}
\widehat\mu_{j,C}=\pi\left(\frac{\sum_{x\in\eta\cap C}\xi_j(x,G)}{c_jN(C)}\right),
\end{equation}
with value zero if $N(C)=0$. {
The scores are nonnegative, so the lower truncation at zero is
inactive here; we retain it to write $\pi$ as the projection of the
real line onto $[0,1]$. The empirical average divided by $c_j$ may
exceed one, whereas its limiting Palm mean $\mu_j(\mathsf M)$ lies
in $[0,1]$, so the upper truncation keeps all coordinates in the
space used for the comparison.
} For any target $m\in[0,1]$, the inequality $|\pi(t)-m|\leq|t-m|$ shows that clipping cannot increase the error. It changes neither the individual scores nor the counts. Proposition~\ref{prop:ergodic} thus gives $\widehat\mu_{j,C_{L,R_j}}\to\mu_j(\mathsf M)$ almost surely, simultaneously for all $j$, under each fixed model. Define the event
\begin{equation}\label{eq:separator}
{
\mathcal A=\left\{
\begin{array}{l}
\text{for every }j\geq1,\quad
 a_j:=\displaystyle\lim_{L\to\infty}\widehat\mu_{j,C_{L,R_j}}
 \text{ exists},\\[2pt]
\text{and }(a_j)_{j\geq1}\in\mathcal C_0
\end{array}\right\}.
}
\end{equation}
{
This definition uses only the observed limits, not the unknown law;
under any fixed law $\mathsf M$, their vector equals $\mu(\mathsf M)$
almost surely. Countability of the coordinates and closedness of
$\mathcal C_0$ make
} $\mathcal A$ measurable; ergodicity and separation give $\mathbf{P}_\varphi (\mathcal A)=1$ and $\Q(\mathcal A)=0$, proving the singularity stated in Theorem~\ref{th:disc}.
\subsection{Unknown profiles: the finite-window test}\label{sec:test}
{
For an observation window $W_L=[-L,L)^d$, with integer $L\geq4$,
set $J_L=\lfloor\log_2L\rfloor$, the number of degree and triangle
statistics used by the test, and let $\mathcal C_0^{(J_L)}$ denote
the projection of the null family onto these first $J_L$ coordinates.
We compare the corresponding empirical averages with this projected set.
} Two errors must remain below the separation: the error in these averages and the tail error, which is at most $2^{-J_L}$. A decreasing tolerance and a slowly increasing number of coordinates will allow a decision without knowing $\Delta_{\Q}^*$. Appendix~\ref{app:finite-references} makes the computation of the minimum explicit.

Ergodicity at a fixed radius gives no rate for a growing list of scores, and therefore does not specify a choice of $J_L$ independent of the profiles. The following uniform bound, proved by a fourth-moment argument in Appendix~\ref{app:fluctuations}, provides the required control.
\begin{lemma}[Uniform control of a per-vertex average]\label{lem:fluctuations}
Let $\xi=d_R$ or $t_R$, with normalisation $c$, and write $\mu=\E_{\mathsf M}^0(\xi/c)\in[0,1]$. If $C$ is a cube with integer vertices and volume $V\geq1$, whose radius-$R$ neighbourhood is observed, then, for $0<b\leq1$,
\begin{equation}\label{eq:fluctuation}
\sup_{\mathsf M}\mathsf M\left(\left|\pi\left(\frac{\sum_{x\in\eta\cap C}\xi(x,G)}{cN(C)}\right)-\mu\right|>b\right)\leq\frac{C_{\lambda,d}(1+R)^{15d}}{V^2b^4}.
\end{equation}
The supremum runs over $\mathbf{P}_\varphi $ and $\Q$     at fixed parameters $(\lambda,d,\gamma)$, each with its own mean $\mu$; $C_{\lambda,d}$ depends neither on the profiles, nor on the radius or the cube.
\end{lemma}
For every integer $L\geq4$, set
\begin{equation}\label{eq:test-schedule}
J_L=\lfloor\log_2L\rfloor,\quad e_L=L^{-1/8},\quad C_L=[-L+J_L,L-J_L)^d,\quad V_L=|C_L|.
\end{equation}
Our enumeration ensures $R_j\leq j\leq J_L$, so all neighbourhoods needed from $C_L$ are observed. Write $\widehat\mu_{L,j}=\widehat\mu_{j,C_L}$ and define
\begin{equation}\label{eq:test}
\mathcal D_L=\inf_{\varphi\in\F}\sum_{j=1}^{J_L}2^{-j}|\widehat\mu_{L,j}-\mu_j(\varphi)|,\qquad T_L=\Ind_{\{\mathcal D_L>3e_L\}}.
\end{equation}
The answer zero means one community, and one means two communities; for $L<4$, set $T_L=0$. Compactness ensures that the minimum exists, and the distance is Lipschitz, hence measurable. The test thus asks whether a single null profile can reproduce all observed moments; its finite computation is treated in Appendix~\ref{app:finite-references}.

To return to the canonical window of volume $n$, define
\begin{equation}\label{eq:test-volume}
L(n)=\lfloor n^{1/d}/2\rfloor,\qquad \tau_n=T_{L(n)}.
\end{equation}
Since $W_{L(n)}\subset B_n$ and $L(n)\to\infty$, this reindexing gives a decision for each observed volume, without a signal-dependent choice.
\begin{proof}[Completion of the proof of Theorem~\ref{th:disc}]
For a fixed model $\mathsf M$, set
\[\mathcal E_L(\mathsf M)=\{\max_{j\leq J_L}|\widehat\mu_{L,j}-\mu_j(\mathsf M)|\leq e_L\}.\]
Since $J_L\leq L/2$ for $L\geq4$, we have $V_L\geq L^d$. The complementary event is the union of the $J_L$ events $\{|\widehat\mu_{L,j}-\mu_j(\mathsf M)|>e_L\}$. Summing their probabilities using Lemma~\ref{lem:fluctuations}, without any independence assumption between scores, gives
\begin{equation}\label{eq:error-bound}
\mathsf M(\mathcal E_L(\mathsf M)^c)\leq p_L:=\min\{1,C_{\lambda,d}J_L(1+J_L)^{15d}L^{1/2-2d}\}.
\end{equation}
This bound is uniform over the profiles, and $\sum_{L\geq4}p_L<\infty$ for every $d\geq1$. On $\mathcal E_L(\mathsf M)$, we have $d_{*,J_L}(\widehat\mu_L,\mu(\mathsf M))\leq e_L$. Under $\mathbf{P}_\varphi $, {
using the generating null profile as a candidate in
} \eqref{eq:test} gives $\mathcal D_L\leq e_L$, hence $T_L=0$. It follows that
\begin{equation}\label{eq:type-one}
\sup_{\varphi\in\F}\mathbf{P}_\varphi (T_L=1)\leq p_L.
\end{equation}
Under a fixed alternative $\Q$, for every $\varphi\in\F$, $d_{*,J_L}(\mu(\Q),\mu(\varphi))\geq\Delta_{\Q}^*-2^{-J_L}$. On $\mathcal E_L(\Q)$, this yields $\mathcal D_L\geq\Delta_{\Q}^*-2^{-J_L}-e_L\geq\Delta_{\Q}^*-2e_L$, since $2^{-J_L}\leq2/L\leq e_L$ for $L\geq4$. Consequently,
\begin{equation}\label{eq:type-two}
\Q(T_L=0)\leq p_L\quad\text{whenever }5L^{-1/8}<\Delta_{\Q}^*.
\end{equation}
Summability and the Borel--Cantelli lemma ensure that $T_L$ is eventually correct almost surely under each fixed model, without independence between windows. Identity~\eqref{eq:test-volume} transfers this property to $\tau_n$. Together with \eqref{eq:separator}, this proves the theorem.
\end{proof}
The exponent $1/8$ is sufficient; no optimality is claimed.

{
\begin{proposition}[Finite implementation of the distinguishability test]
\label{prop:operational-main}
Assume the hypotheses of Theorem~\ref{th:disc}, with known
$(\lambda,d,\gamma)$, and the computational setting of
Appendix~\ref{app:finite-references}: $\lambda$ and $\gamma$ can be
computed with certified precision, and observed positions admit
exact geometric comparisons. For each integer $L\geq4$, there is
an explicit finite family $\mathcal S_L$ of null profiles,
independent of the generating profiles, and a randomised test
\[
\widehat T_L=\Ind_{\{\widehat D_L>4\varepsilon_L\}},\qquad
\widehat D_L=\min_{s\in\mathcal S_L}
\sum_{j=1}^{J_L}2^{-j}
|\widetilde\mu_{L,j}-\widehat m_j(s)|,
\]
where $\varepsilon_L=2^{-k_L}$ with
$k_L=\min\{k\in\N:L\leq2^{8k}\}$.
Here $\widetilde\mu_{L,j}$ approximates the observed normalised
moment to accuracy $\varepsilon_L/2$, and the reference moments
$\widehat m_j(s)$ are computed by finitely many independent
simulations, with simultaneous error at most $\varepsilon_L/2$
except on an event of probability at most $L^{-2}$.
The procedure terminates almost surely and, writing
$\overline{\mathsf M}$ for the law including these simulations,
\[
\sup_{\varphi\in\F}\overline{\mathbf P}_\varphi(\widehat T_L=1)
\leq b_L+L^{-2},\qquad
\overline{\Q}(\widehat T_L=0)\leq b_L+L^{-2}
\quad\text{if }\Delta_{\Q}^*>7\varepsilon_L,
\]
for a summable sequence $b_L$. Thus its answers are eventually
correct almost surely under each fixed admissible model.
\end{proposition}
\begin{proof}
Appendix~\ref{app:finite-references} constructs the profiles in
Proposition~\ref{prop:explicit-net}, controls the simulated moments
in Lemma~\ref{lem:reference-MC}, and proves termination and the
error bounds in Theorem~\ref{th:operational}.
\end{proof}
This replaces the ideal minimisation over all profiles by a finite
comparison, with explicit control of the additional approximation
errors; it does not assert that the number of reference profiles
is small enough for efficient computation at moderate volumes.
}

{
\begin{proof}[Proof of Corollary~\ref{cor:unknown-intensity-test}]
\phantomsection\label{rem:lambda-test}
For distinguishability, replace $\lambda$ by $\widehat\lambda_L^{\rm win}=N(W_L)/|W_L|$ in \eqref{eq:test}, using the threshold $L^{-1/16}$. The normalised moments have the form $\lambda^p A/(1+\lambda^p B)$, $p\in\{1,2\}$ and $0\leq A\leq B$; their derivative with respect to $\log\lambda$ is therefore uniformly bounded. Clipping gives the same stability for the empirical coordinates. The fourth Poisson moment yields $\Prob(|\widehat\lambda_L^{\rm win}/\lambda-1|>L^{-1/8})\leq C_\lambda L^{1/2-2d}$, a summable bound. Together with \eqref{eq:error-bound}, Borel--Cantelli shows that the computed distance eventually differs from the distance at the true intensity by at most $C_\lambda L^{-1/8}$, almost surely. This error is negligible compared with the new threshold, while the separation $\Delta_{\Q}^*>0$ remains: the test is still consistent. When $\widehat\lambda_L^{\rm win}=0$, we give it the value zero.
\end{proof}
}

 {
\paragraph{Why the guarantee is pointwise in the model.}
The preceding tests eventually succeed for every fixed admissible
alternative, but this does not give a volume sufficient for all
nonzero signals at once. The next proposition shows that this
limitation is unavoidable for any test when the alternative may
approach the null as the observation window grows; it is not an
additional assumption needed by our positive results.
}

\begin{proposition}[No uniform consistency without signal separation]\label{prop:no-uniform}
No sequence of tests has both errors tending to zero uniformly over all null profiles and all alternatives satisfying (H).
\end{proposition}
\begin{proof}
Take $\varphi(t)=e^{-t}/2$ and, for $0<b<1$, $\phin^{(b)}(t)=(1+b)e^{-t}/2$ and $\phout^{(b)}(t)=(1-b)e^{-t}/2$. Their mean is $\varphi$, and each pair is admissible. Couple the models on $B_n$ with the same positions, ages and uniform edge variables, adding labels for the alternative $Q_b$. Each edge differs with probability at most $b/2$, hence
\[\TV(\mathbf{P}_\varphi ^{(n)},Q_b^{(n)})\leq\frac b2\E\binom{N_n}{2}=\frac{b\lambda^2n^2}{4}.\]
Total variation controls every decision event. If $A$ is the event that a test answers two communities, then
\[\mathbf{P}_\varphi ^{(n)}(A)+Q_b^{(n)}(A^c)=1-\bigl(Q_b^{(n)}(A)-\mathbf{P}_\varphi ^{(n)}(A)\bigr)\geq1-\TV(\mathbf{P}_\varphi ^{(n)},Q_b^{(n)}).\]
With $b_n=(n+1)^{-3}$, total variation tends to zero, so the sum of the errors of any test remains asymptotically at least one. Uniform convergence of both errors to zero is therefore impossible.
\end{proof}
This does not contradict Theorem~\ref{th:disc}, which fixes the alternative before increasing the window, whereas the proof above weakens the signal with $n$. On a family separated by $\Delta_{\Q}^*\geq\Delta_0>0$, bound~\eqref{eq:type-two} is uniform.

\section{Recovery}\label{sec:recovery}
We first reconstruct vertices of age at least $a_0$, using slices on which the kernel varies within bounds that allow us to check the local conditions for \GBG\ (Section~\ref{sec:slice}); the intersections then align the partitions up to a global sign (Section~\ref{sec:propagation}). Voting proceeds downwards through dyadic slices towards smaller ages: its margin grows as $\lambda a_k^{-\gamma}$, while the influence of one voter on the slice being processed remains of order $\lambda$. When ages are hidden, a calibration population provides two degree ranks, one spatially truncated for initialisation and the other based on the full degree for voting (Figure~\ref{fig:recovery-route}).

\begin{figure}[htbp]
\centering\footnotesize
\setlength{\unitlength}{1cm}
\begin{picture}(14.5,2.1)
\put(0,1.2){\framebox(3.1,.75){\shortstack{Ages above $a_0$\\\GBG\ within each slice}}}
\put(3.15,1.58){\vector(1,0){.55}}
\put(3.75,1.2){\framebox(3.0,.75){\shortstack{Sign alignment\\on intersections}}}
\put(6.8,1.58){\vector(1,0){.55}}
\put(7.4,1.2){\framebox(3.2,.75){\shortstack{Dyadic voting\\towards smaller ages}}}
\put(10.65,1.58){\vector(1,0){.55}}
\put(11.25,1.2){\framebox(3.1,.75){\shortstack{Below cutoff $\zeta$\\label $+1$}}}
\put(.3,.15){\makebox(14,0.6){\shortstack{Hidden ages: independent calibration, two coloured passes\\and alignment on their shared colour.}}}
\end{picture}
\caption{From local initialisation to voting. Colours separate calibration data from the edges used for reconstruction.}
\label{fig:recovery-route}
\end{figure}

\noindent\textit{Dependencies.} We do not assume independence between vertices for the outputs of \GBG, the labels already assigned, or the estimated ages. We control error sets through inequalities that hold for every realisation, and then apply the ergodic theorem to auxiliary scores. The local independence needed for Poisson laws and calibration is specified where it is used.
{
Reconstruction repeatedly uses the dependence of degree on age
established in Lemma~\ref{lem:degree}, both to control incorrect voters
and to replace hidden ages by degree ranks.
}

\subsection{Quantitative initialisation on an age slice}\label{sec:slice}
Distinguishability detects the presence of communities but gives no individual labels; we now turn to their reconstruction. From this point on, the profiles are nonincreasing and satisfy (A0) and (H), and $d\geq2$; ages are observed and the model parameters are known.

Our first step is to apply \GBG\ on an age slice. We use its combinatorial and percolation analysis, specifying the properties needed to transfer it to the marked model. This separates the conclusions supplied by the algorithm from what remains to be checked for the preferential attachment kernel.

\noindent\textit{Review of \GBG\ and the result used.}
Consider a marked process with spatial intensity $\mu$ and independent uniform community labels. At scale $\mathcal R$, divide space into cubes of side $\ell=\mathcal R/(4\sqrt d)$; the block $L_1(z)$ consists of the $3^d$ cells at $\ell_\infty$ distance at most one from $z$, and has diameter $3\mathcal R/4$. A cell is A-Good if it contains at least $\max\{\mu\ell^d(1-\varepsilon),1\}$ vertices and no cycle of vertices in $L_1(z)$ has a product of pair classifications equal to $-1$,  {
which would make the pair classifications inconsistent around that cycle
}. This property is observable, and the algorithm then propagates relative labels through A-Good components, according to \cite[Algorithms~2--3 and Definition~14]{abbe}.

For the analysis only, a cell is T-Good if it satisfies the same  $ \max\{\mu\ell^d(1-\varepsilon),1\}$ occupancy threshold and every pair in $L_1(z)$ is classified according to the product of its true labels.  {
This implies that every cycle has product of pair classifications
$+1$, since each true label appears twice in that product, and
therefore the cell is also A-Good.
}. We use the combinatorial conclusion of \cite[Propositions~5--6 and the proof of Proposition~4]{abbe}: a T-Good cell is A-Good, and the output of \GBG\ agrees with the truth, up to a sign, on each connected T-Good component. This conclusion is taken from ABS; the local and quantitative estimates required to apply it to the marked model are established below.  

\paragraph{The three local checks.}
The analysis of \GBG\ requires the following properties, which we check before turning to connected components:
\begin{enumerate}[label=\textbf{(H\arabic*)}]
\item Under the two-point Palm law, the pair statistic is Poisson, with means $\mu\mathcal M_{\mathrm{in}}$ and $\mu\mathcal M_{\mathrm{out}}$ according to the product of the labels.
\item For admissible marks and distances at most $3\mathcal R/4$, we have $\mathcal M_{\mathrm{in}}-\mathcal M_{\mathrm{out}}\geq D>0$ and $\mathcal M_{\mathrm{in}}\leq U<\infty$.
\item The T-Good event is determined within a fixed spatial neighbourhood of the cell, independent of the intensity and volume.
\end{enumerate}
The first two properties give an exponentially small local error by Chernoff's bound, while the third preserves this control when passing to cell components. We give this step directly in a quantitative form, avoiding a later repetition of the qualitative argument in \cite[Sections~5.1--5.5]{abbe}. Truncating common neighbours ensures (H3), even for profiles with unbounded support.

\begin{lemma}[Uniform gap under monotonicity and (H)]\label{lem:signal-interval}
There exist $0<t_-<t_+<\infty$ and $b_0>0$ such that
\begin{equation}\label{eq:signal-interval}
\psi(t)\geq b_0\qquad\text{for every }t\in[t_-,t_+].
\end{equation}
\end{lemma}
\begin{proof}
Monotone profiles have only countably many discontinuities, whereas $\{\psi>0\}$ has positive measure. We may therefore choose a common continuity point $t_0>0$ in this set; in a sufficiently small neighbourhood, their difference remains above $\psi(t_0)/2$, giving the interval and $b_0$.
\end{proof}
Fix such an interval and set
\begin{equation}\label{eq:contrast}
\vartheta=\frac{t_-}{t_+}\in(0,1).
\end{equation}

We now choose slices on which the kernel varies little enough to preserve this connection gap.
Fix $t_0\in(0,1)$ and $0<w<t_0$ with $t_0+w\leq1$, and set
$I=[t_0-w,t_0+w]=[a_I,b_I],\quad |I|=2w$.
The restriction to ages in $I$ is a marked Poisson process with spatial intensity
\begin{equation}\label{eq:effective-intensity}
\mu=\lambda|I|=2w\lambda.
\end{equation}
Its vertex ages are independent and uniform on $I$, independently of the uniform community labels. Write $G_I$ for the induced graph.

\paragraph{Bounds and admissibility.}
For $I=[a_I,b_I]$, the positive exponents of $g$ sum to one, so
\begin{equation}\label{eq:kernel-sandwich}
a_I\leq g(u,v)\leq b_I\qquad(u,v\in I,\ \gamma\in(0,1)).
\end{equation}
Each nonincreasing profile then satisfies $\varphi(b_I|x-y|^d)\leq\varphi(\kappa_{XY})\leq\varphi(a_I|x-y|^d)$.

Define the spatial radii
\begin{equation}\label{eq:effective-radii}
\mathcal R_I=\left(\frac{t_+}{b_I}\right)^{1/d},\qquad
\rho_I=\left(\frac{t_-}{a_I}\right)^{1/d}.
\end{equation}
If $\rho_I\leq|x-z|\leq\mathcal R_I$, then $t_-\leq\kappa_{XZ}\leq t_+$ for every pair of ages in $I$. The interval $[\rho_I,\mathcal R_I]$ is therefore uniformly informative across the slice: this bound gives a local problem comparable with the unmarked setting.

The inequality $\rho_I<\mathcal R_I$ is equivalent to $\vartheta(t_0+w)<t_0-w$, or
\begin{equation}\label{eq:admissibility}
\frac{w}{t_0}<\frac{1-\vartheta}{1+\vartheta}.
\end{equation}

\begin{lemma}[A uniformly informative annular region]\label{lem:annulus}
Assume \eqref{eq:admissibility}. For $0<|x-y|\leq3\mathcal R_I/4$, set
\[
A_I(x,y)=\{z:\rho_I\leq|x-z|\leq\mathcal R_I,
\ \rho_I\leq|y-z|\leq\mathcal R_I\}.
\]
There is $v_I>0$, depending only on $\rho_I,\mathcal R_I,d$, such that $|A_I(x,y)|\geq v_I$ uniformly in $x,y$. For every $z\in A_I(x,y)$ and all ages of $X,Y,Z$ in $I$, both $\psi(\kappa_{XZ})$ and $\psi(\kappa_{YZ})$ are at least $b_0$.
\end{lemma}
\begin{proof}
Set $\ell_- =\max(\rho_I,3\mathcal R_I/8)<\mathcal R_I$, $\ell=(\ell_-+\mathcal R_I)/2$ and $s_I=(\mathcal R_I-\ell_-)/4>0$. Since $d\geq2$ and $|x-y|/2\leq3\mathcal R_I/8<\ell$, there is a point $z_0$ such that $|z_0-x|=|z_0-y|=\ell$. For every $z\in B(z_0,s_I)$, both distances lie between $\ell-s_I>\rho_I$ and $\ell+s_I<\mathcal R_I$. Thus $B(z_0,s_I)\subset A_I(x,y)$, and we may take $v_I=\omega_ds_I^d$. The signal bound then follows from \eqref{eq:kernel-sandwich}, \eqref{eq:effective-radii} and \eqref{eq:signal-interval}.
\end{proof}
This construction is where dimension at least two is used, while the strict inequality in \eqref{eq:admissibility} prevents the annular region from having zero volume. It is a sufficient geometric condition for this argument, not a claimed necessary condition for reconstruction on a slice.

Write $S_I(x,y)=B(x,\mathcal R_I)\cap B(y,\mathcal R_I)$. For $X=(x,u)$ and $Y=(y,v)$ with ages in $I$, define
\begin{align}
M^I_{\mathrm{in}}(X,Y)
&=\frac12\int_I\int_{S_I(x,y)}
\bigl[\phin(\kappa_{XZ})\phin(\kappa_{YZ})
+\phout(\kappa_{XZ})\phout(\kappa_{YZ})\bigr]\dd z\dd t,
\label{eq:Min}\\
M^I_{\mathrm{out}}(X,Y)
&=\frac12\int_I\int_{S_I(x,y)}
\bigl[\phin(\kappa_{XZ})\phout(\kappa_{YZ})
+\phout(\kappa_{XZ})\phin(\kappa_{YZ})\bigr]\dd z\dd t,
\label{eq:Mout}
\end{align}
where $Z=(z,t)$. In particular,
\begin{equation}\label{eq:pair-gap}
M^I_{\mathrm{in}}-M^I_{\mathrm{out}}
=\frac12\int_I\int_{S_I(x,y)}
\psi(\kappa_{XZ})\psi(\kappa_{YZ})\dd z\dd t.
\end{equation}

\begin{proposition}[Distribution and separation of pair counts]
\label{prop:pair-poisson}\label{prop:uniform-gap}
Under the two-point Palm law, conditional on the labels of $X,Y$, the number $C_I(X,Y)$ of common neighbours in $S_I(x,y)$ with ages in $I$ is Poisson, with mean $\lambda M^I_{\rm in}$ or $\lambda M^I_{\rm out}$ according to whether their signs agree or differ. For ages in $I$ and $0<|x-y|\leq3\mathcal R_I/4$,
\begin{equation}\label{eq:uniform-gap}
M^I_{\rm in}-M^I_{\rm out}\geq D_I:=wb_0^2v_I>0,
\qquad 0\leq M^I_{\rm out}\leq M^I_{\rm in}
\leq U_I:=2w\omega_d\mathcal R_I^d,
\end{equation}
with constants independent of $\gamma\in(0,1)$.
\end{proposition}
\begin{proof}
By Slivnyak's theorem, the other vertices retain their Poisson law; after averaging over their signs, independent retention as common neighbours has probability given by the integrand of \eqref{eq:Min} or \eqref{eq:Mout}. Thinning therefore gives the stated law. In \eqref{eq:pair-gap}, restricting the integral to $I\times A_I(x,y)$ gives $(|I|/2)b_0^2v_I=D_I$, while the averaged integrand in \eqref{eq:Min} is at most one and the lens has volume at most $\omega_d\mathcal R_I^d$.
\end{proof}

We now have two uniformly separated means. The pair classifier returns $+1$ if
\begin{equation}\label{eq:pair-threshold}
C_I(X,Y)\geq\frac\lambda2
\bigl(M^I_{\mathrm{in}}(X,Y)+M^I_{\mathrm{out}}(X,Y)\bigr),
\end{equation}
and $-1$ otherwise, to estimate the product $\sigma_X\sigma_Y$. The threshold is computed from the observed endpoint positions and ages, together with the known parameters.

\begin{proposition}[Uniform pair error]\label{prop:pair-error}
There is $c_I>0$, independent of $\lambda$ and $\gamma\in(0,1)$, such that the Palm probability of incorrect classification is at most $e^{-c_I\lambda}$ for every pair with $0<|x-y|\leq3\mathcal R_I/4$ and ages in $I$. One may take
\begin{equation}\label{eq:chernoff-constant}
c_I=U_I\mathfrak h\left(\frac{D_I}{2U_I}\right),\qquad
\mathfrak h(t)=(1+t)\log(1+t)-t.
\end{equation}
\end{proposition}
\begin{proof}
Write $a=M^I_{\mathrm{out}}$, $b=M^I_{\mathrm{in}}$ and $D=b-a\geq D_I$. Under opposite labels, the statistic is Poisson with mean $\lambda a$, and the threshold exceeds it by $\lambda D/2$. If $a>0$, Chernoff's bound gives $\exp[-\lambda a\mathfrak h(D/(2a))]$. If $a=0$, the statistic is almost surely zero and this error is impossible, since the threshold is strictly positive. Under equal labels, the relative lower deviation is $D/(2b)\leq1/2$. Its rate $(1-t)\log(1-t)+t$ dominates $\mathfrak h(t)$, so the error is at most $\exp[-\lambda b\mathfrak h(D/(2b))]$.

For fixed $D>0$, the function $a\mapsto a\mathfrak h(D/(2a))$ is nonincreasing: with $v=D/(2a)$, we have $\frac{\dd}{\dd v}(\mathfrak h(v)/v)
=(v-\log(1+v))/v^2\geq0$. It is also nondecreasing in $D$. Since both means are at most $U_I$, both rates are bounded below by the strictly positive constant in \eqref{eq:chernoff-constant}. Averaging the conditional error probabilities preserves this bound.
\end{proof}

For the transfer, $\mu=2w\lambda$ and $\mathcal M=M^I/(2w)$: the width factor appears only once, either in the intensity or in the age integral. The normalised bounds are $b_0^2v_I/2$ and $\omega_d\mathcal R_I^d$.

\begin{lemma}[Probability of a T-Good cell]\label{lem:tgood}
Let $\ell_I=\mathcal R_I/(4\sqrt d)$ and fix $\varepsilon\in(0,1/2)$. For sufficiently large $\lambda$,
\begin{equation}\label{eq:tgood-bound}
\Prob(z\text{ is T-Bad})
\leq e^{-2w\lambda\ell_I^d\mathfrak h(\varepsilon)}
+(2w\lambda)^2(3\ell_I)^{2d}e^{-c_I\lambda}.
\end{equation}
The bound is uniform in $z$. The T-Good field has finite dependence range depending only on $d$.
\end{lemma}
\begin{proof}
The cell count is Poisson with mean $2w\lambda\ell_I^d$; for large enough $\lambda$, the occupancy threshold exceeds one and Chernoff's bound gives the first term of \eqref{eq:tgood-bound}. For pairs, Campbell--Mecke and Proposition~\ref{prop:pair-error} give
\[
\E\!\sum_{X\neq Y\in\eta_I\cap L_1(z)}
\Ind_{\{\text{incorrectly classified pair}\}}
\leq(2w\lambda)^2|L_1(z)|^2e^{-c_I\lambda}.
\]
Since $|L_1(z)|=(3\ell_I)^d$, a union bound gives the second term, without independence between errors. All points and edges consulted lie within distance $2\mathcal R_I$ of the cell; independence of disjoint restrictions therefore gives dependence range $M=\lceil16\sqrt d\rceil+2$ in lattice coordinates, uniformly in $\lambda$.
\end{proof}

We still need to transfer the local estimates to the proportion of correctly reconstructed vertices; this is the purpose of passing through components of reliable cells.
Lemma~\ref{lem:tgood} gives $q_I(\lambda)\leq C_Ie^{-c_I\lambda}$: the factor $\lambda^2$ is absorbed by decreasing the exponent. Call a cell open when it is T-Good. Two indices are nearest neighbours if their $\ell_1$ distance is one, and $*$-neighbours if their $\ell_\infty$ distance is one. A $*$-connected set can therefore use diagonal steps; this notion arises in boundaries of open components, although propagation follows nearest-neighbour steps.

\begin{lemma}[Peierls bound under finite dependence]\label{lem:peierls-rate}
Let $(Y_z)_{z\in\mathbb Z^d}$, with $d\geq2$, be a stationary $M$-dependent field, with $M\in\N$: families indexed by two sets at $\ell_\infty$ distance greater than $M$ are independent. Assume $\Prob(Y_z=0)\leq p_{\mathrm{bad}}$ and set $K_M=(M+1)^d$. For sufficiently small $p_{\mathrm{bad}}$,
\[
\Prob(0\text{ belongs to no infinite open component})
\leq C_d p_{\mathrm{bad}}^{1/K_M}.
\]
\end{lemma}
\begin{proof}
The origin is closed with probability at most $p_{\mathrm{bad}}$. If it is open but belongs to a finite component, its outer boundary contains a $*$-connected separating set of closed sites. The number of such sets of $j$ sites is at most $A_d^j$. Indeed, a separator of size $j$ meets a fixed coordinate ray within distance $j$; once this point is chosen, a traversal of a spanning tree bounds the number of $*$-connected sets exponentially. The polynomial factor from the anchor point is absorbed into $A_d^j$. Divide the lattice into $K_M$ residue classes modulo $M+1$. At least $j/K_M$ sites of a fixed set lie in one class, where their indicators are jointly independent. The probability that the entire separator is closed is thus at most $p_{\mathrm{bad}}^{j/K_M}$. Hence,
\[
\Prob(0\text{ outside infinite open components})
\leq p_{\mathrm{bad}}+\sum_{j\geq1}(A_d p_{\mathrm{bad}}^{1/K_M})^j
\leq (1+2A_d)p_{\mathrm{bad}}^{1/K_M}
\]
provided $A_d p_{\mathrm{bad}}^{1/K_M}\leq1/2$. This is the Peierls argument for high-density percolation. ABS combine it with domination by product measures due to Liggett, Schonmann and Stacey \cite{lss}; here we retain explicit dependence on $p_{\mathrm{bad}}$ and $M$, which is needed when the scales vary with intensity.
\end{proof}

We use the boundary results of \cite[Lemma~2 and Theorem~4]{timar}: the outer boundary of a nearest-neighbour connected set is $*$-connected, and the $*$-boundary of a finite $*$-connected set visible from infinity is nearest-neighbour connected. These two properties preserve the quantitative bound in a finite window.

\begin{lemma}[A large component in a finite box]\label{lem:finite-box}
Under the assumptions of Lemma~\ref{lem:peierls-rate}, for sufficiently small $p_{\mathrm{bad}}$, let $Q_H=[-H,H]^d\cap\mathbb Z^d$ and let $\mathcal K_H$ be an open component of maximal size in the induced box. Then
\[
\frac{|Q_H\setminus\mathcal K_H|}{|Q_H|}
\leq C_d p_{\mathrm{bad}}^{1/K_M}+o_{\mathbb P}(1)
\qquad(H\to\infty).
\]
\end{lemma}
\begin{proof}
Counting $*$-connected sets and dividing into residue classes gives, with $b=A_d p_{\mathrm{bad}}^{1/K_M}<1$,
\[
\Prob(\exists\text{ a closed }*\text{-connected component of size}
\geq h\text{ meeting }Q_H)\leq C|Q_H|b^h.
\]
Indeed, such a component contains a connected set of $h$ sites anchored in the box. The same argument at a fixed site shows that closed components are almost surely finite.

On the complement of this event, join two infinite open sites of $Q_{H-2h-2}$ by a nearest-neighbour path in this cube. Each closed component encountered can be bypassed along its $*$-boundary visible from infinity, which is nearest-neighbour connected by \cite[Theorem~4]{timar}. The endpoints lie outside its filled hull because they have an infinite open path; all sites on the detour are open, otherwise they would belong to the maximal closed component. Since its diameter is less than $h$, the detour stays in its bounding box enlarged by one unit, hence in $Q_H$. Each replacement removes closed sites without adding any, so the procedure terminates. Without the box constraint, this argument also proves uniqueness of the infinite open component.

Take $h=\lceil c_{p_{\mathrm{bad}}}\log(H+2)\rceil$ with $c_{p_{\mathrm{bad}}}$ sufficiently
large. The exceptional probability tends to zero, and all
infinite sites in the inner cube then belong to one
component of $Q_H$.
{
To estimate their number, let $I_z$ be the indicator that
$z$ belongs to an infinite open component of the full
lattice configuration. The underlying field is stationary
and ergodic by finite dependence. Since $(I_z)$ is a
translation-covariant function of this field, it is also
stationary and ergodic, although it need not have finite
dependence. The lattice ergodic theorem and
Lemma~\ref{lem:peierls-rate} therefore give
\[
\frac{1}{|Q_H|}\sum_{z\in Q_H}I_z
\longrightarrow
\mathbb P(I_0=1)
\geq 1-C_d p_{\mathrm{bad}}^{1/K_M}
\qquad\text{almost surely}.
\]
Here $|Q_H|$ denotes the number of lattice sites in $Q_H$.
The boundary strip contains $O(hH^{d-1})=o(|Q_H|)$ sites,
so removing it does not change the limiting density.
On the event established above, all the remaining sites
counted by $I_z$ lie in a single open component inside
$Q_H$, which gives the claimed component-size bound.
}
This completes the proof.

\end{proof}

\begin{proposition}[Quantitative reconstruction on a slice]\label{prop:vertex-rate}
For a fixed strictly admissible slice and a fixed occupancy tolerance in $(0,1/2)$, there exist $C,c>0$ and $\lambda_0<\infty$ such that, for $\lambda\geq\lambda_0$, \GBG\ is correct up to a sign except on a proportion at most $Ce^{-c\lambda}$ of its vertices, with probability tending to one as $n\to\infty$.
\end{proposition}
\begin{proof}
Lemma~\ref{lem:finite-box}, applied to the cells of the slice, bounds their proportion outside the large component in the window by $Cq_I^{1/K_M}+o_{\mathbb P}(1)$. We must convert this to a vertex count, since cells have different occupancies. Write $N_z$ for the count in a cell and $\mathcal B_n$ for the cells outside the component. Deterministically,
\[
\sum_{z\in\mathcal B_n}N_z
\leq\sqrt{|\mathcal B_n|\sum_zN_z^2}.
\]
If $\nu=2w\lambda\ell_I^d$, then $\E N_z=\nu$ and $\E N_z^2=\nu^2+\nu$. The spatial law of large numbers for these Poisson counts therefore bounds the limiting vertex proportion by
$C\sqrt{1+\nu^{-1}}\,q_I^{1/(2K_M)}$.
For large $\lambda$, the multiplicative factor is bounded, and $q_I\leq C_Ie^{-c_I\lambda}$ makes this expression exponentially small. On the T-Good component, \GBG\ returns the true labels up to a sign, by \cite[Propositions~5--6 and the proof of Proposition~4]{abbe}, applied to the slice classifications. The boundary strip of fixed width is negligible; increasing $C$ gives the strict margin needed for the probability statement.
\end{proof}

\subsection{Propagation and reconstruction over all ages}\label{sec:propagation}
Local partitions are defined up to a sign. We align them through their intersections, then bound the contribution of ages below the cutoff. The estimators may be dependent: alignment relies on deterministic inequalities.

\begin{definition}[A cover suitable for propagation]\label{def:cover}
A family $V_1,\ldots,V_m$ covers a nonempty finite set $V$ suitably for propagation if consecutive sets intersect and $V_i\cap V_j=\varnothing$ whenever $|i-j|>1$.
\end{definition}

\begin{lemma}[Global alignment]\label{lem:global-alignment}
Let $\sigma:V\to\{-1,+1\}$ and let $(V_i)_{i=1}^m$ be such a cover. Suppose the functions $s_i:V_i\to\{-1,+1\}$ satisfy
\begin{equation}\label{eq:local-signed-overlap}
\frac1{|V_i|}\sum_{X\in V_i}s_i(X)\sigma_X\geq\alpha
\qquad (i=1,\ldots,m)
\end{equation}
for some $\alpha\in[0,1]$. Set $i_0(X)=\min\{i:X\in V_i\}$ and $\tilde\sigma_X=s_{i_0(X)}(X)$. Then
\begin{equation}\label{eq:global-alignment}
\frac1{|V|}\sum_{X\in V}\tilde\sigma_X\sigma_X\geq2\alpha-1.
\end{equation}
\end{lemma}
\begin{proof}
Set $V_0=\varnothing$ and $S=\sum_{i=2}^m|V_{i-1}\cap V_i|$. Since nonconsecutive sets are disjoint, $\sum_i|V_i|=|V|+S$ and
\[
\sum_{X\in V}\tilde\sigma_X\sigma_X
=\sum_i\sum_{X\in V_i\setminus V_{i-1}}s_i(X)\sigma_X
\geq\alpha\sum_i|V_i|-S
=\alpha|V|-(1-\alpha)S.
\]
The odd-indexed sets are pairwise disjoint, as are the even-indexed sets. Therefore $\sum_i|V_i|\leq2|V|$, so $S\leq|V|$. Dividing by $|V|$ gives \eqref{eq:global-alignment}.
\end{proof}
Alignment thus costs at most a factor of two in the error proportion, independently of the number of slices.

\begin{lemma}[Identifying the relative sign]\label{lem:sign}
Let $A,B$ be nonempty finite sets with nonempty intersection. Suppose the signed overlap of $s_A$ with $\sigma$ on $A$ is at least $\alpha$, and the absolute overlap of $s_B$ on $B$ is at least $\alpha$. Set
\[
S_{A,B}=\sum_{X\in A\cap B}s_A(X)s_B(X),\qquad
H_{A,B}=(\alpha-1)(|A|+|B|)+|A\cap B|.
\]
If $s_B$ is positively oriented, then $S_{A,B}\geq H_{A,B}$; if negatively oriented, $S_{A,B}\leq-H_{A,B}$. Thus, when $H_{A,B}>0$, the sign of $S_{A,B}$ determines the relative orientation using only estimated labels.
\end{lemma}
\begin{proof}
Removing $A\setminus B$ decreases the sum by at most its cardinality, so
\[
\sum_{X\in A\cap B}s_A(X)\sigma_X
\geq(\alpha-1)|A|+|A\cap B|.
\]
In the positive case, the same bound holds for $s_B$, with $B$ replacing $A$. For $f,g\in\{-1,+1\}$, we have $fg\geq f+g-1$. Apply this to $f=s_A\sigma$ and $g=s_B\sigma$ on the intersection and sum, obtaining $S_{A,B}\geq H_{A,B}$. The negative case reduces to the positive one by replacing the estimator with $-s_B$.
\end{proof}

We apply these two lemmas to a chain of consecutive slices, using their intersections to recover relative orientations.
For an age set $J$, write $\bar \eta _{ n}(J)=\{X=(x,u)\in\bar \eta _{ n}:u\in J\}$ and $N_n(J)=\#\bar \eta _{ n}(J)$. Independent thinning gives
\begin{equation}\label{eq:age-counts}
N_n(J)\sim\operatorname{Poisson}(\lambda n|J|).
\end{equation}
Counts associated with disjoint age sets are independent.

\begin{figure}[htbp]
\centering
\setlength{\unitlength}{1cm}
\begin{picture}(12.4,4.6)
\put(1,4){\vector(1,0){10.5}}
\multiput(1,3.9)(2,0){6}{\line(0,1){0.2}}
\put(0.8,4.3){\makebox(0.4,0){$a$}}
\put(10.8,4.3){\makebox(0.4,0){$1$}}
\put(1,3.55){\makebox(2,0){$J_0$}}
\put(3,3.55){\makebox(2,0){$J_1$}}
\put(5,3.55){\makebox(2,0){$J_2$}}
\put(7,3.55){\makebox(2,0){$J_3$}}
\put(9,3.55){\makebox(2,0){$J_4$}}
\linethickness{2pt}
\put(1,2.9){\color{black}\line(1,0){4}}
\put(3,2.15){\color{teal!75!black}\line(1,0){4}}
\put(5,1.4){\color{black}\line(1,0){4}}
\put(7,0.65){\color{teal!75!black}\line(1,0){4}}
\put(0,2.9){\makebox(0.7,0){$I_1$}}
\put(0,2.15){\makebox(0.7,0){$I_2$}}
\put(0,1.4){\makebox(0.7,0){$I_3$}}
\put(0,0.65){\makebox(0.7,0){$I_4$}}
\linethickness{0.35pt}
\put(4,2.7){\vector(0,-1){0.35}}
\put(6,1.95){\vector(0,-1){0.35}}
\put(8,1.2){\vector(0,-1){0.35}}

\end{picture}
\caption{Each slice consists of two elementary intervals. Shared vertices determine the relative signs; each vertex belongs to at most two slices.}
\label{fig:age-slices}
\end{figure}
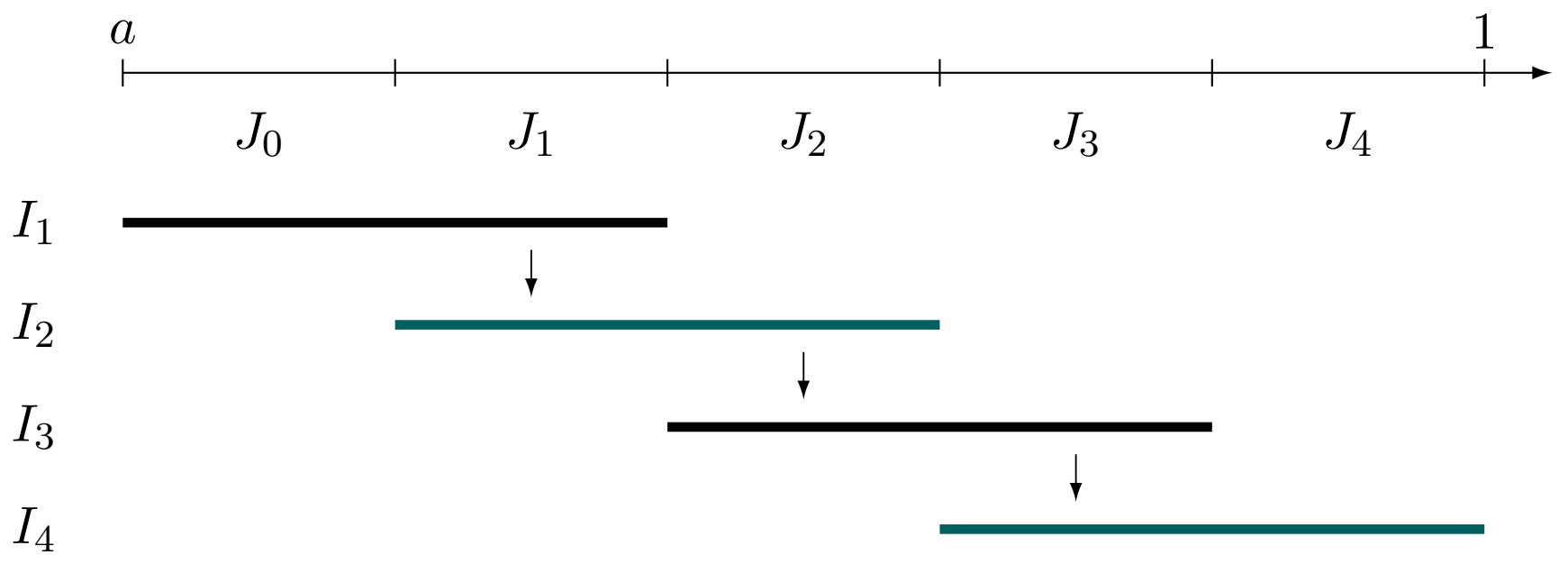

Fix a cutoff $a\in(0,1)$ and choose $w_*>0$ such that
\begin{equation}\label{eq:global-width}
w_*<\frac{a(1-\vartheta)}{2\vartheta}.
\end{equation}
Take an integer $m\geq2$ satisfying $w=(1-a)/(m+1)\leq w_*$ and define
\begin{equation}\label{eq:slices}
J_j=[a+jw,a+(j+1)w],\quad j=0,\ldots,m,\qquad
I_i=J_{i-1}\cup J_i=[t_i-w,t_i+w],\quad t_i=a+iw.
\end{equation}
The slices $I_i$ cover $[a,1]$, and only neighbours share an interval of length $w$ (Figure~\ref{fig:age-slices}). Almost surely, no vertex lies at an endpoint.

The slices in \eqref{eq:slices} are admissible: $t_m+w=1$, and the strongest constraint, at the first centre $t_1=a+w$, is equivalent to $2w\vartheta<a(1-\vartheta)$, which follows from \eqref{eq:global-width}.
The sets $V_i=\bar \eta _{ n}(I_i)$ therefore form a suitable cover whenever consecutive intersections are nonempty. This event has probability tending to one, although it need not occur for every finite $n$.

\begin{definition}[\GBG\ with propagation]\label{def:algorithm}
Given $(\bar \eta _{ n},E_n)$, the partition above and fixed occupancy tolerances, construct the estimator on $\bar \eta _{ n}([a,1])$ as follows:
\begin{enumerate}
\item For each $i=1,\ldots,m$, apply \GBG\ to the data induced by ages in $I_i$, obtaining $s_i:V_i\to\{-1,+1\}$.
\item Process $i=1,\ldots,m-1$ in this order and compute $S_i=\sum_{X\in\bar \eta _{ n}(J_i)}s_i(X)s_{i+1}(X)$. If $S_i<0$, replace $s_{i+1}$ by $-s_{i+1}$; if $S_i=0$, keep the orientation.
\item For $X\in\bar \eta _{ n}([a,1])$, return $\tilde\sigma_X=s_{i_0(X)}(X)$, where $i_0(X)=\min\{i:X\in V_i\}$.
\end{enumerate}
The fixed enumeration breaks all local ties; on an empty set, use the empty label function. The estimator is therefore a deterministic function of the observations.
\end{definition}

\begin{proposition}[Quantitative initialisation above a cutoff]
\label{prop:exp-seed}
Fix $a\in(0,1)$ and the chain~\eqref{eq:slices}. There exist $C_a,c_a>0$ and $\lambda_a<\infty$ such that, for each $\lambda\geq\lambda_a$, the propagated estimator satisfies, in a suitable global orientation,
\[
\frac{\#\{X\in\bar \eta _{ n}([a,1]):\tilde\sigma_X\neq\sigma_X\}}{N_n([a,1])}
\leq C_ae^{-c_a\lambda}+o_{\mathbb P}(1).
\]
The same bound holds with denominator $N_n$.
\end{proposition}
\begin{proof}
Proposition~\ref{prop:vertex-rate} gives an error proportion $C_ie^{-c_i\lambda}+o_{\mathbb P}(1)$ on each slice. Since the chain is finite and fixed, the maximum of the $C_i$ and the positive minimum of the $c_i$ give a common bound. For large enough $\lambda$, all local overlaps exceed a common $\alpha>3/4$ with probability tending to one; call this event $\mathcal E_3$.

It remains to show that the observations allow the correct orientations to be chosen. Set $t_n=\lambda nw$. The counts in the $J_j$ are independent Poisson variables of mean $t_n$, so all intersections are nonempty on an event $\mathcal E_1$ with probability at least $1-(m-1)e^{-t_n}$. For two neighbouring slices, write $P_i,R_i,Q_i$ for the counts in the three elementary intervals, with $R_i$ corresponding to the intersection. The quantity in Lemma~\ref{lem:sign} becomes
\[
H_i=(\alpha-1)(P_i+Q_i)+(2\alpha-1)R_i,
\qquad H_i/t_n\xrightarrow{\mathbb P}4\alpha-3>0.
\]
The event $\mathcal E_2$ on which all $H_i$ are positive therefore also has probability tending to one. On $\mathcal E_1\cap \mathcal E_2\cap \mathcal E_3$, correlations successively orient the slices according to the first one, without requiring independence of the estimators.

Lemma~\ref{lem:global-alignment} then bounds the merged error proportion by twice the largest local proportion. A finite union of failure events gives the claimed bound, and $N_n([a,1])\leq N_n$ gives the second normalisation. Increasing $C_a$, if needed, yields an inequality holding with probability tending to one rather than an $o_{\mathbb P}(1)$ remainder.
\end{proof}

\begin{corollary}[Strong recovery with observed ages]
\label{prop:strong}
Under (A0), (H), $d\geq2$ and monotonicity of the profiles, for every $\Gamma<1$ there is a finite threshold $\lambda_\Gamma$, independent of $\gamma\in(0,1)$, above which propagation between slices, extended by $+1$ below a cutoff, satisfies $\Prob(\ov_n>\Gamma)\to1$. In particular, weak recovery is possible at high intensity and strong recovery as the intensity grows, with observed ages.
\end{corollary}
\begin{proof}
Fix $\Gamma\in(0,1)$, $a=(1-\Gamma)/8$, $w_* =a(1-\vartheta)/(4\vartheta)$ and $\alpha_*=(1+\Gamma)/2$. The partition~\eqref{eq:slices}, with $w\leq w_*$, is admissible. Proposition~\ref{prop:exp-seed} guarantees overlap above $\alpha_*$ on $V_+=\bar \eta _{ n}([a,1])$, for sufficiently large intensity independently of $\gamma$. Set $V_-=\bar \eta _{ n}((0,a))$ and extend the estimator by $+1$ on $V_-$.

The Poisson counts give $|V_-|/N_n\to a$ in probability. On the event where $|V_-|/N_n<2a$ and the overlap on $V_+$ exceeds $\alpha_*$, the reverse triangle inequality gives
\[
\ov_n\geq\alpha_*\frac{|V_+|}{N_n}-\frac{|V_-|}{N_n}
=\alpha_*-(1+\alpha_*)\frac{|V_-|}{N_n}
>\alpha_*-4a=\Gamma.
\]
This event has probability tending to one. Choosing $\Gamma=1/2$ gives weak recovery, while arbitrary $\Gamma<1$ gives strong recovery. The cutoff and number of slices are fixed before $\lambda$, and then $n\to\infty$.
\end{proof}
\begin{remark}[Uniformity in $\gamma$]
The kernel bounds on a slice allow $D_I,U_I,c_I$ and the dependence range to be chosen independently of $\gamma\in(0,1)$; the same is therefore true of the constants in Proposition~\ref{prop:exp-seed} and the threshold $\lambda_\Gamma$ in the preceding corollary. This uniformity concerns initialisation and its qualitative consequence. The sharper rates below involve constants depending on $\gamma$, particularly the sum $\sum_k a_k^\gamma$.
\end{remark}

\subsection{Completing reconstruction by stratified voting}
\label{sec:vote}
Propagation initialises vertices above $a_0$; we then descend towards smaller ages through dyadic slices, each receiving votes from all previously classified neighbours. This stratification keeps the ages in a new slice comparable with its upper endpoint, allowing us to compare the oracle margin with the influence of incorrect voters on that slice alone.

Set $I_{\rm in}=\int_{\R^d}\phin(|z|^d)\dd z$ and $I_{\rm out}=\int_{\R^d}\phout(|z|^d)\dd z$, then $J=I_{\rm in}-I_{\rm out}>0$, $S=I_{\rm in}+I_{\rm out}$ and $\bar I=S/2=I_{\pbar}$. These constants are finite by (A0), so we write $m(s)=\bar I\Theta(s)$ for the mean degree without the intensity factor. For $0<s<a$, the kernel gives
\begin{equation}\label{eq:upward-PA}
\int_a^1g(s,t)^{-1}\dd t=\frac{s^{-\gamma}(1-a^\gamma)}{\gamma}.
\end{equation}
Smaller ages thus have more neighbours among previously classified vertices, and the aim is to use this growing margin without overestimating the influence of earlier errors.

Let $a_k=a_0 2^{-k}$ and $k=0,\ldots,K$, and set
$Y_k=\bar \eta _{ n}([a_k,1]),\quad \mathcal V_k=\bar \eta _{ n}([a_{k+1},a_k)),\quad \zeta=a_K$.
Initialise labels on $Y_0$ using Proposition~\ref{prop:exp-seed} and, for $X\in \mathcal V_k$, define
\begin{equation}\label{eq:stratified-vote}
V_k(X)=\sum_{Y\in Y_k}\Ind_{\{X\sim Y\}}\hat\sigma_Y,
\qquad \hat\sigma_X=\sgn V_k(X),
\end{equation}
assigning $+1$ to a zero vote, with each already labelled neighbour contributing one vote, without age weighting. Vertices of $\mathcal V_k$ are classified simultaneously, without using their internal edges, and then join the voters for the next step: for $a_0=1/4$, we process $[1/8,1/4)$, then $[1/16,1/8)$, and so on. Below the final cutoff $\zeta$, vertices receive $+1$ and are counted in the error bound. For the analysis, introduce
$\Omega_k(X)=\sum_{Y\in A\setminus\{X\}:\,u_Y\geq a_k}\Ind_{\{X\sim Y\}}\sigma_Y$, where $A$ is the reconstruction population in the infinite graph.
With hidden ages, the same definitions will use estimated sets $\widehat Y_k,\widehat{\mathcal V}_k$, while $\Omega_k$ remains defined by the true ages. Although this oracle vote is not observable, it provides a reference for measuring the effect of incorrect labels and displaced slice boundaries.

In this section all vertices participate, so $p=1$, but to reuse the same lemma with hidden ages, we allow independent selection of each vertex with probability $p$, giving a reconstruction population of intensity $p\lambda$. The colours in Section~\ref{sec:ages} implement this selection with $p=2/3$; they are not communities. The enlarged age interval in the lemma anticipates errors in placing vertices into estimated slices.

\begin{lemma}[Uniform oracle margin on enlarged slices]\label{lem:oracle-PA}
Suppose the reconstruction population is obtained by independent colouring with proportion $p\in(0,1]$, and let $a\leq a_0\leq1/4$. Set
\begin{equation}\label{eq:oracle-constants}
U=\frac{2^{-\gamma}(1-(2a_0)^\gamma)}{\gamma}>0,
\qquad v=\frac{pJU}{4}>0.
\end{equation}
For a Palm vertex of true age $s\in[a/4,2a]$, the oracle vote using all vertices of $A$ with ages at least $a$ in the infinite graph satisfies
\[
\Prob_s\bigl(\sigma_X\Omega_a(X)\leq v\lambda a^{-\gamma}\bigr)
\leq C e^{-c\lambda a^{-\gamma}}\leq C e^{-c\lambda}.
\]
\end{lemma}
\begin{proof}
The signed oracle vote is the difference of two independent Poisson counts $P_+-P_-$, with means $\Lambda I_{\rm in}$ and $\Lambda I_{\rm out}$, where $\Lambda=(p\lambda/2)\int_a^1g(s,t)^{-1}\dd t$. This statement holds under Palm, conditional only on the age and sign of the vertex; to bound its mean from below, restrict the integral to $t\geq2a\geq s$:
\[
\int_a^1g(s,t)^{-1}\dd t
\geq s^{-\gamma}\frac{1-(2a)^\gamma}{\gamma}
\geq U a^{-\gamma}.
\]
For $h>0$, exponential Markov gives
\[
\Prob(P_+-P_-\leq\Lambda J/2)
\leq\exp\{\Lambda[\,hJ/2+I_{\rm in}(e^{-h}-1)
+I_{\rm out}(e^h-1)\,]\}.
\]
The bracketed expression is zero at $h=0$ and has derivative $-J/2<0$ there. Choose a fixed $h$ small enough to make it negative. Since $v\lambda a^{-\gamma}\leq\Lambda J/2$, the claimed bound follows.
\end{proof}

A voter's total degree is too crude a bound on its influence when its label is wrong: only edges to the slice currently being processed can create new errors. Set
\[
T_k=[a_k/4,2a_k],\qquad m_k=v\lambda a_k^{-\gamma},\qquad
D_k(Y)=\#\{X\in A\setminus\{Y\}:X\sim Y,\ u_X\in T_k\}.
\]
The degrees $D_k$ are defined in the infinite graph, and $N_n^A=\#(A\cap B_n)$. The enlargement to $T_k$ will also be useful when ages are estimated.

\begin{lemma}[Influence on the current slice]\label{lem:target-degree}
Set $\Lambda_*=14e^2p\bar I$ and $\lambda_B=(7e^2p\bar I)^{-1}$. For $\lambda\geq\lambda_B$ and every $k$, a vertex of age $t\geq a_k/2$ satisfies
\[
\E_t\bigl[D_k(Y)\Ind_{\{D_k(Y)>\Lambda_*\lambda\}}\bigr]
\leq 7ep\bar I\lambda e^{-\Lambda_*\lambda}.
\]
Thus, if
$\mathcal T_n^{(k)}=\sum_{Y\in A\cap B_n:\,u_Y\geq a_k/2}
D_k(Y)\Ind_{\{D_k(Y)>\Lambda_*\lambda\}}$,
then, for fixed $k,\lambda$,
\begin{equation}\label{eq:target-tail}
\frac{\mathcal T_n^{(k)}}{N_n^A}\leq7ep\bar I\lambda e^{-\Lambda_*\lambda}+o_{\mathbb P}(1),
\qquad
\sum_{Y\in W}D_k(Y)\leq\Lambda_*\lambda|W|+\mathcal T_n^{(k)}
\end{equation}
for every $W\subset\{Y\in A\cap B_n:u_Y\geq a_k/2\}$, even if chosen from the graph.
\end{lemma}
\begin{proof}
Under Palm at age $t$, $D_k$ is Poisson with mean
$p\lambda\bar I\int_{T_k}g(t,s)^{-1}\dd s$. Now $g(t,s)\geq t\wedge s\geq a_k/4$ and $|T_k|=7a_k/4$, so this mean is at most $7p\bar I\lambda$, uniformly in $k$. For a Poisson variable $Z$ of mean $\mu$, the identity
$\E[Z\Ind_{\{Z\geq j\}}]=\mu\Prob(Z\geq j-1)$
and the bound $\Prob(Z\geq x)\leq e^{-x}$ for $x\geq e^2\mu$ give the estimate: with $j=\lfloor\Lambda_*\lambda\rfloor+1$, we have $j-1\geq\Lambda_*\lambda-1\geq7e^2p\bar I\lambda$. The Poisson bound itself follows from exponential Markov with parameter $\log(x/\mu)\geq2$. The exceptional score is translation-covariant and integrable, so the spatial ergodic theorem gives the first relation in \eqref{eq:target-tail}. The second follows for each realisation by separating degrees below the threshold from the others.
\end{proof}

\begin{lemma}[Edges lost at the window boundary]\label{lem:vote-boundary}
Write $\ell_n(X)=\#\{Y\in A\setminus B_n:Y\sim X\}$ for $X\in A\cap B_n$. If $B$ is a calibration population independent of $A$, with intensity $\rho_{\mathrm{cal}}\lambda$, set
$G_n=\{X\in A\cap B_n:X\text{ has a neighbour in }B\setminus B_n\}$.
For fixed $\lambda$,
\begin{equation}\label{eq:boundary-loss}
\frac1n\sum_{X\in A\cap B_n}\ell_n(X)\xrightarrow{\mathbb P}0,
\qquad \frac{|G_n|}{n}\xrightarrow{\mathbb P}0,
\qquad \frac1n\sum_{Y\in G_n}D(Y)\xrightarrow{\mathbb P}0,
\end{equation}
where $D(Y)$ is the degree within $A$ in the infinite graph. Outside $G_n$, the full calibration degree is the same in the window and in the infinite graph.
\end{lemma}
\begin{proof}
For $s>0$, introduce
\[
F(s,r)=\int_0^1\int_{|z|\geq r}\pbar(g(s,t)|z|^d)\dd z\dd t.
\]
We have $F(s,r)\downarrow0$ as $r\to\infty$ and $F(s,r)\leq\bar I\Theta(s)$, whose integral over $s$ is finite by (A0). If $r(x)=\dist(x,B_n^c)$, Campbell--Mecke bounds the expectation of the first term in \eqref{eq:boundary-loss} by
\[
p^2\lambda^2\int_0^1\frac1n\int_{B_n}F(s,r(x))\dd x\dd s.
\]
For fixed $R$, the strip $r(x)\leq R$ occupies a volume fraction $O(Rn^{-1/d})$, and elsewhere $F(s,r(x))\leq F(s,R)$. Let $n\to\infty$ and then $R\to\infty$; dominated convergence and Markov give the first limit. The number of vertices in $G_n$ is bounded by the number of lost calibration edges, to which the same calculation applies.

Conditional on the age and sign of $Y$, its degree within $A$ is independent of its edges to $B$. Consequently,
\[
\E_{y,t}[D(Y)\Ind_{\{Y\in G_n\}}]
\leq p\lambda\bar I\Theta(t)\min\{1,\rho_{\mathrm{cal}}\lambda F(t,r(y))\}.
\]
The same boundary decomposition, now dominated by $p\lambda\bar I\Theta(t)$, proves the third limit without requiring a second moment of the degree. Finally, a calibration degree changes only if one of its edges leaves the window.
\end{proof}

\begin{proposition}[Recursion with observed ages]\label{prop:observed-recursion}
Orient the truth once according to the initialisation and set
$W_k=\{X\in Y_k:\hat\sigma_X\neq\sigma_X\}$, $e_{k,n}=|W_k|/N_n$.
With $p=1$, $c_v=4\Lambda_*/v$ and $\lambda\geq\lambda_B$, for every finite $K$, simultaneously for $k<K$,
\begin{equation}\label{eq:observed-recursion}
e_{k+1,n}\leq(1+c_va_k^\gamma)e_{k,n}
+C_1a_ke^{-c_1\lambda}+C_2a_k^\gamma e^{-\Lambda_*\lambda}
+o_{\mathbb P}(1).
\end{equation}
The constants depend neither on $k$ nor on $K$.
\end{proposition}
\begin{proof}
The observed vote compares with the infinite oracle vote through
$|V_k(X)-\Omega_k(X)|\leq2\#\{Y\in W_k:Y\sim X\}+\ell_n(X)$.
Unless the oracle margin fails, a new error therefore requires at least $m_k/4$ incorrectly labelled neighbours, or $\ell_n(X)\geq m_k/2$. The number of vertices of the first kind is at most
\[
\frac4{m_k}\sum_{Y\in W_k}D_k(Y)
\leq c_va_k^\gamma|W_k|+\frac{4a_k^\gamma}{v\lambda}\mathcal T_n^{(k)},
\]
since voters have age at least $a_k$. Lemma~\ref{lem:target-degree} controls the second term, while Lemma~\ref{lem:vote-boundary} makes the proportion affected by the boundary negligible. Finally, oracle failures are counted on the true slice $[a_k/2,a_k)$, of length $a_k/2$; Lemma~\ref{lem:oracle-PA} and ergodicity bound their proportion by $C_1a_ke^{-c_1\lambda}+o_{\mathbb P}(1)$. Adding previous errors gives \eqref{eq:observed-recursion}. The factor $a_k^\gamma$ in front of the influence of $W_k$ is essential: the degree towards the slice stays of order $\lambda$, while the margin is $v\lambda a_k^{-\gamma}$.
\end{proof}

\subsection{Universal reconstruction: replacing parameters with ranks}
\label{sec:ages}
We now remove ages, profiles and parameters $(\gamma,\beta,\lambda)$ from the algorithm's inputs, retaining positions, edges and dimension. Mean degree remains strictly decreasing with age; instead of inverting it using known parameters, we estimate a vertex's age by the proportion of vertices whose calibration degree is at least as large as its own. Pair-classification thresholds are also replaced by observed quantiles, while voting retains its previous form, although initialisation requires scales that gradually adapt to any fixed model.

Give vertices independent uniform colours in $\{1,2,3\}$, independently of the graph. In the first pass, $A$ consists of colours 1 and 2, and $B$ of colour 3; in the second, $A$ consists of colours 2 and 3, and $B$ of colour 1. Reconstruct $A$ using only its edges to $B$ for calibration, then its internal edges for \GBG\ and voting. Thus $p=2/3$ and $\rho_{\mathrm{cal}}=1/3$ are the reconstruction and calibration proportions, while the shared colour will align the two outputs; colours are neither communities nor estimates of them.

\paragraph{A fully observable procedure.}
Set $a_0=1/4$ and $\widehat\lambda_n=N_n/n$. If $\widehat\lambda_n>0$, compute $\Lambda_n=2^{\lfloor\log_2\widehat\lambda_n\rfloor}$. For $\Lambda=\Lambda_n\geq e^e$, choose the scales from the estimated intensity as follows:
\begin{equation}\label{eq:universal-scales}
m=\left\lceil(1-a_0/2)\log\Lambda\right\rceil,\quad
w=\frac{1-a_0/2}{m+1},\quad
\ell=\frac1{\log\Lambda},\quad R=\log\log\Lambda,\quad K=\lceil\Lambda\rceil.
\end{equation}
Logarithms are natural unless stated otherwise; cells have side $\ell$, and $R_c=4\sqrt d\,\ell$. The quantity $\Lambda_n$ is the largest power of two not exceeding $\widehat\lambda_n$; it allows scales to be chosen from a discrete family. These choices contain no unknown constants: slices and cells shrink at high intensity, but their mean occupancy is of order $\lambda/(\log\lambda)^{d+1}$ and tends to infinity, while the calibration radius grows slowly.

For $X\in A\cap B_n$, let $D_B^{(n)}(X)$ be its number of neighbours in $B\cap B_n$, and $D_{B,R}^{(n)}(X)$ the number within distance $R$. With $N_n^A=\#(A\cap B_n)$, define
\begin{equation}\label{eq:rank-estimates}
\begin{aligned}
S_n(j)&=\frac1{N_n^A}\sum_{Y\in A\cap B_n}\Ind_{\{D_B^{(n)}(Y)\geq j\}},&
\widehat u_X&=S_n(D_B^{(n)}(X)),\\
S_{R,n}(j)&=\frac1{N_n^A}\sum_{Y\in A\cap B_n}\Ind_{\{D_{B,R}^{(n)}(Y)\geq j\}},&
\widetilde u_X&=S_{R,n}(D_{B,R}^{(n)}(X)).
\end{aligned}
\end{equation}
The rank $\widehat u_X$ determines voting slices, while $\widetilde u_X$ is used for initialisation. Degree ties are kept without random tie-breaking, and although the truncated degree is local, its empirical rank uses the whole window; we therefore do not assign it finite dependence range.

Partition $[a_0/2,1]$ into $m+1$ elementary intervals $J_0,\ldots,J_m$ of length $w$, half-open except for the last, which includes 1, and set $I_i=J_{i-1}\cup J_i$, $1\leq i\leq m$. For $V_i=\{X\in A\cap B_n:\widetilde u_X\in I_i\}$, consider distinct ordered pairs in $V_i$ at distance at most $3R_c/4$. Denote their set by $\mathcal P_i$, and define their common-neighbour count by
\[
C_i(X,Y)=\#\{Z\in V_i\setminus\{X,Y\}: |z-x|\leq R,\ |z-y|\leq R,\ Z\sim X,\ Z\sim Y\}.
\]
For $\alpha\in\{1/4,3/4\}$, let $\widehat q_{i,\alpha}$ be the smallest integer $c$ such that at least a fraction $\alpha$ of pairs in $\mathcal P_i$ satisfy $C_i\leq c$. We then use
\begin{equation}\label{eq:empirical-pair-threshold}
\widehat t_i=\frac{\widehat q_{i,1/4}+\widehat q_{i,3/4}}2,\qquad
\mathrm{PC}_i(X,Y)=\begin{cases}+1,&C_i(X,Y)\geq\widehat t_i,\\-1,&C_i(X,Y)<\widehat t_i.\end{cases}
\end{equation}
The first quartile identifies pairs from opposite communities, and the third pairs from the same community; Appendix~\ref{app:universal} justifies this interpretation without knowing labels or assuming independent counts.

\begin{definition}[Universal reconstruction]\label{def:hidden-algorithm}
In each pass, compute \eqref{eq:universal-scales}--\eqref{eq:empirical-pair-threshold}, then proceed as follows.
\begin{enumerate}
\item Run \GBG\ in each $V_i$, using classifications $\mathrm{PC}_i$ and occupancy threshold $\max\{p\Lambda_n w\ell^d,1\}$. Align outputs on successive intersections as in Definition~\ref{def:algorithm}.
\item Keep these labels on $\widehat Y_0=\{\widehat u_X\geq a_0\}$, using $+1$ where no label is available.
\item For $k=0,\ldots,K-1$, set $a_k=a_02^{-k}$ and simultaneously classify vertices of $\widehat{\mathcal V}_k=\{a_{k+1}\leq\widehat u_X<a_k\}$ using the vote \eqref{eq:stratified-vote} from neighbours in $\widehat Y_k=\{\widehat u_Y\geq a_k\}$, with their previously assigned labels. Ties receive $+1$.
\item Assign $+1$ to vertices with rank below $a_K$.
\end{enumerate}
Orient the second output by the sign of its correlation with the first on colour 2, then keep the first output on colours 1 and 2 and the oriented second output on colour 3. A zero correlation is assigned sign $+1$.
When $\widehat\lambda_n=0$, $\Lambda_n<e^e$, or a reconstruction population is empty, the corresponding output is $+1$. Vertices in a slice with no pairs also receive label $+1$; component choices and tie-breaking follow the orders fixed in Section~\ref{sec:model}. These conventions define the procedure on every finite observation.
\end{definition}

\paragraph{Why ranks suffice for voting.}
The algorithm computes ranks without inverting a known degree function. For the analysis, let $D_B$ and $D_{B,R}$ denote the calibration degrees of a typical vertex in the infinite graph. Their survival functions, that is, the probabilities of having at least $j$ neighbours, are
\[
S(j)=\int_0^1\Prob(\operatorname{Pois}(\rho_{\mathrm{cal}}\lambda\bar I\Theta(t))\geq j)\dd t,
\quad
S_R(j)=\int_0^1\Prob(\operatorname{Pois}(\rho_{\mathrm{cal}}\lambda m_R(t))\geq j)\dd t,
\]
where $m_R(t)=\int_0^1\int_{|z|\leq R}\pbar(g(t,u)|z|^d)\dd z\dd u$. Conditional on age $t$, degree is Poisson; integration averages this law over the uniform age, as in the degree analysis of \cite[Section~4]{gracar2019}. If degree equalled its mean, which is strictly decreasing in $t$, its descending rank would be exactly $\Prob(U\leq t)=t$. We must therefore control two approximations: the empirical survival function to its theoretical counterpart as volume grows, and rank to age as intensity grows. Ties are retained: no rank is assumed exactly uniform at fixed intensity. As before, calculations use $\beta=1$; any fixed $\beta$ is absorbed into the profiles for this analysis, without transforming the algorithm's data.

\begin{lemma}[Rank accuracy and substitution for ages]\label{lem:rank-vote}
For fixed $\lambda$ and $R$, $\sup_j|S_n(j)-S(j)|$ and $\sup_j|S_{R,n}(j)-S_R(j)|$ tend to zero in probability. If $0<\rho\leq1$, then, for sufficiently large $\lambda$,
\begin{equation}\label{eq:rank-relative}
\Prob_s\!\left(S(D_B(X))\notin[s/(1+\rho),(1+\rho)s]\right)
\leq2e^{-\kappa_\rho\lambda s^{-\gamma}},\qquad
\kappa_\rho=\frac{\rho_{\mathrm{cal}}\bar I\rho^2}{144C_\gamma},\qquad C_\gamma=\frac1{1-\gamma}+\frac1\gamma.
\end{equation}
Fix $\delta\in(0,1]$, $1\leq K<\infty$, $\zeta=a_K$, and set
\begin{equation}\label{eq:rank-exceptions}
\mathcal X_{\mathrm{exc}}=\{X\in A:S(D_B(X))\notin[u_X/(1+\delta/3),(1+\delta/3)u_X]\},\quad
\Omega_n=\left\{\sup_j|S_n(j)-S(j)|\leq\delta\zeta/8\right\}.
\end{equation}
On $\Omega_n$, outside $\mathcal X_{\mathrm{exc}}\cup G_n$ from Lemma~\ref{lem:vote-boundary}, the following inclusions hold simultaneously for $k<K$:
\begin{align*}
\widehat u_X\in[a_k/2,a_k)&\ \Longrightarrow\ u_X\in T_k=[a_k/4,2a_k],\\
\widehat u_X\geq a_k&\ \Longrightarrow\ u_X\geq a_k/2,\\
\Ind_{\{\widehat u_X\geq a_k\}}\neq\Ind_{\{u_X\geq a_k\}}
&\ \Longrightarrow\ u_X\in \mathcal J_k^\delta=[a_k/(1+\delta),(1+\delta)a_k),\\
\widehat u_X<\zeta&\ \Longrightarrow\ u_X<2\zeta,\qquad
\widehat u_X\geq a_0\ \Longrightarrow\ u_X\geq2a_0/3.
\end{align*}
Equalities at true-age boundaries are negligible. Moreover, $\Prob(\Omega_n)\to1$.
\end{lemma}
The proof is given in Appendix~\ref{app:universal-ranks}. Its key point is that the exceptional set $\mathcal X_{\mathrm{exc}}$ uses the deterministic theoretical survival function $S$ and calibration edges only. Membership of a vertex in $\mathcal X_{\mathrm{exc}}$ therefore remains independent of an incident voting edge, conditional on ages and signs. This property would not follow by directly replacing $\mathcal X_{\mathrm{exc}}$ with an event defined from all empirical ranks.

\paragraph{Initialisation without knowing the signal.}
The universal scales eventually suit any fixed signal: variation of the means within a slice becomes small compared with the gap between the two pair types, and the quartiles then place the threshold between these means. We must nevertheless preserve the locality needed for percolation, despite global empirical calibration.

\begin{proposition}[Universal initialisation]\label{prop:hidden-seed}
Under the reconstruction assumptions, for sufficiently large $\lambda$, the initial output on $\widehat Y_0$ satisfies, in a suitable orientation,
\[
e_{0,n}:=\frac{\#\{X\in\widehat Y_0:\hat\sigma_X\neq\sigma_X\}}{N_n^A}
\leq\varepsilon_{\rm init}(\lambda)+o_{\mathbb P}(1),\qquad
\varepsilon_{\rm init}(\lambda)\leq
\exp\!\left[-\frac{\lambda}{(\log\lambda)^{2d+3}}\right].
\]
The intensity threshold depends on the model, but none of its parameters enters the procedure.
\end{proposition}
\begin{proof}
Appendix~\ref{app:universal-local} constructs an auxiliary field of certified cells, with dependence range $M\leq C_d(R/\ell+1)$, whose defect probability satisfies
\begin{equation}\label{eq:universal-defect}
q^\sharp\leq e^{-c\lambda w\ell^d}
+C\lambda R^d e^{-c\lambda w^2}
+C(\lambda\ell^d)^2e^{-c\lambda w}.
\end{equation}
On a global event of probability tending to one, every certified cell is T-Good for the empirical decisions. We apply Lemma~\ref{lem:finite-box} to the certificates, not to empirical ranks. With $K_M=(M+1)^d$, the chosen scales give
\[
(q^\sharp)^{1/(2K_M)}\leq
\exp\!\left[-\frac{c\lambda}{(\log\lambda)^{2d+1}(\log\log\lambda)^d}\right]=:\epsilon_3(\lambda).
\]
The count in a selected cell is bounded by the count of all vertices of $A$, whose first two moments are Poisson moments. Cauchy--Schwarz, as in Proposition~\ref{prop:vertex-rate}, together with selected density at least $cw$, therefore gives error proportion at most $C\epsilon_3/w+o_{\mathbb P}(1)$ on each slice. Here we use the same conclusion of \GBG\ as in that proposition, from \cite[Propositions~5--6 and the proof of Proposition~4]{abbe}.

To align slices, the appendix bounds each selection between true intervals enlarged or reduced by $h=c_hw$, with $0<c_h\leq1/32$, except for a proportion $Ce^{-c\lambda h^2}$. In the notation $P_i,R_i,Q_i$ of Lemma~\ref{lem:sign}, the proportions of the two outer parts are at most $w+2h+Ce^{-c\lambda h^2}+o_{\mathbb P}(1)$, and that of the intersection is at least $w-2h-Ce^{-c\lambda h^2}-o_{\mathbb P}(1)$. For $\alpha=7/8$, the sign criterion, divided by $N_n^A$, is therefore bounded below by
\[
(4\alpha-3)w-2h-Ce^{-c\lambda h^2}-o_{\mathbb P}(1)>0.
\]
Since the number of slices is finite at fixed intensity, Lemma~\ref{lem:global-alignment} transfers errors with a bounded factor, without multiplying them by the number of slices.

Finally, outside $\mathcal X_{\mathrm{exc}}\cup G_n$, a vertex of $\widehat Y_0$ has true age at least $2a_0/3$ by Lemma~\ref{lem:rank-vote}. It is covered by the local slices unless its local rank fails, an event of density $Ce^{-c\lambda h^2}$. Restriction to $\widehat Y_0$ therefore adds only these failures, those of $\mathcal X_{\mathrm{exc}}$, and a negligible boundary contribution. We obtain $e_{0,n}\leq C\epsilon_3/w+Ce^{-c\lambda}+Ce^{-c\lambda h^2}+o_{\mathbb P}(1)$; {in fact, this deterministic bound is
$o(e^{-b(\lambda)})$, where
\[
b(\lambda):=\frac{\lambda}{(\log\lambda)^{2d+3}}.
\]
Indeed, $w,h\asymp(\log\lambda)^{-1}$ and
\[
\frac{\lambda/[(\log\lambda)^{2d+1}(\log\log\lambda)^d]}
{b(\lambda)}
=\frac{(\log\lambda)^2}{(\log\log\lambda)^d}\longrightarrow\infty,
\qquad
\frac{\lambda h^2}{b(\lambda)}\longrightarrow\infty.
\]
Also $\lambda/b(\lambda)\to\infty$ and
$\log(C/w)=O(\log\log\lambda)=o(b(\lambda))$.
Thus we may choose
$\varepsilon_{\rm init}(\lambda)=o(e^{-b(\lambda)})$,
which implies the stated bound and leaves room to absorb
constant factors in the final reconstruction.}

These arguments first apply to each deterministic value $\Lambda\in[\lambda/2,\lambda]$. Since $\widehat\lambda_n\to\lambda$, its lower dyadic power is, with probability tending to one, a single fixed value if $\lambda$ is not dyadic, or one of $\lambda/2,\lambda$ otherwise. Taking the union of these two failure events suffices, without conditioning a local field on the global intensity estimate.
\end{proof}

Voting can now use exactly the same oracle comparison. We detail this transfer because it explains why nonparametric estimation does not introduce a fresh error at every slice.
\begin{lemma}[Neighbours displaced by age estimation]\label{lem:misplaced-neighbours}
Take $\mathcal X_{\mathrm{exc}}$ from Lemma~\ref{lem:rank-vote} and $0<\delta\leq\delta_1$, where $\delta_1=\min\{1,v/(64e^2p\bar I)\}$. Set
\[
\begin{aligned}
\mathcal J_k^\delta&=[a_k/(1+\delta),(1+\delta)a_k),\\
N_k^{\rm bd}(X)&=\#\{Y\in A:Y\sim X,\ u_Y\in \mathcal J_k^\delta\},\qquad
N^{\mathrm{exc}}(X)=\#\{Y\in \mathcal X_{\mathrm{exc}}:Y\sim X\}.
\end{aligned}
\]
For sufficiently large $\lambda$, in the infinite graph under Palm at age $s\in T_k$, and for either sign,
\begin{equation}\label{eq:misplaced-neighbours}
\Prob_s(N_k^{\rm bd}\geq m_k/8)\leq e^{-v\lambda/8},\qquad
\Prob_s(N^{\mathrm{exc}}\geq m_k/16)\leq\frac{128p\bar I C_\gamma}{v}e^{-\kappa_{\delta/3}\lambda}.
\end{equation}
\end{lemma}
\begin{proof}
The count $N_k^{\rm bd}$ concerns a true-age band, so it is Poisson with mean $p\lambda\bar I\int_{\mathcal J_k^\delta}g(s,t)^{-1}\dd t$. On this band, $g(s,t)\geq a_k/4$ and $|\mathcal J_k^\delta|\leq2\delta a_k$, so the mean is at most $8p\bar I\delta\lambda\leq e^{-2}m_k/8$. The Poisson bound used in Lemma~\ref{lem:target-degree} gives the first inequality.

No Poisson law is required for the second count. Under Palm at two vertices of $A$, the edge $X\sim Y$ uses a variable independent of the calibration of $Y$, once their positions, ages and signs are fixed. Moreover, the probability $\pi_\delta(t)=\Prob_t(Y\in \mathcal X_{\mathrm{exc}})$ does not depend on the sign of $Y$. Averaging over this sign, Mecke gives
\[
\E_sN^{\mathrm{exc}}(X)=p\lambda\bar I\int_0^1\frac{\pi_\delta(t)}{g(s,t)}\dd t
\leq2p\lambda\bar I\Theta(s)e^{-\kappa_{\delta/3}\lambda}
\leq8p\lambda\bar I C_\gamma a_k^{-\gamma}e^{-\kappa_{\delta/3}\lambda}.
\]
Markov's inequality at threshold $m_k/16$ concludes the proof. Colour separation is used precisely in this two-vertex factorisation; it does not make age errors independent across vertices.
\end{proof}

\begin{proposition}[Error recursion without observed ages]\label{prop:hidden-recursion}
In the orientation given by initialisation, let
$W_k=\{X\in\widehat Y_k:\hat\sigma_X\neq\sigma_X\}$ and $e_{k,n}=|W_k|/N_n^A$.
With $p=2/3$, $\rho_{\mathrm{cal}}=1/3$, $\delta\leq\delta_1$, $\lambda\geq\lambda_B$ and $c_v=4\Lambda_*/v$, for sufficiently large $\lambda$ and every finite $K$, on $\Omega_n$ from Lemma~\ref{lem:rank-vote},
\begin{equation}\label{eq:hidden-recursion}
e_{k+1,n}\leq(1+c_va_k^\gamma)e_{k,n}
+\frac{|\widehat{\mathcal V}_k\cap \mathcal X_{\mathrm{exc}}|}{N_n^A}
+C_1a_ke^{-c_1\lambda}+C_2a_k^\gamma e^{-\Lambda_*\lambda}
+o_{\mathbb P}(1),\qquad k<K.
\end{equation}
The constants are independent of $k,K$; one may take $C_2=28ep\bar I/v$, and $c_1>0$ smaller than the oracle exponents, $v/8$ and $\kappa_{\delta/3}$.
\end{proposition}
\begin{proof}
Start with a comparison valid for every realisation. Set
$Y_k^{(n)}=\{Y\in A\cap B_n:u_Y\geq a_k\}$,
$R_k=\widehat Y_k\triangle Y_k^{(n)}$ and $W'_k=W_k\setminus(\mathcal X_{\mathrm{exc}}\cup G_n)$.
An incorrect label changes a vote by two units, while adding or removing a voter changes it by one; comparison with the infinite oracle therefore gives
$|V_k-\Omega_k|\leq2\#\{Y\in W_k:Y\sim X\}+\#\{Y\in R_k:Y\sim X\}+\ell_n(X)$.
On $\Omega_n$, the inclusions in Lemma~\ref{lem:rank-vote} give
$R_k\subset\{u\in \mathcal J_k^\delta\}\cup \mathcal X_{\mathrm{exc}}\cup G_n$.
Writing $N_{W'}(X)$ and $N^G(X)$ for neighbour counts in $W'_k$ and $G_n$, we obtain
\begin{equation}\label{eq:vote-perturbation}
|V_k(X)-\Omega_k(X)|
\leq2N_{W'}(X)+3N^{\mathrm{exc}}(X)+3N^G(X)+N_k^{\rm bd}(X)+\ell_n(X).
\end{equation}

A vertex of $\widehat{\mathcal V}_k\setminus(\mathcal X_{\mathrm{exc}}\cup G_n)$ has true age in $T_k$. On this deterministic band, introduce
\[
\mathcal E_k(X)=\{\sigma_X\Omega_k(X)\leq m_k\}
\cup\{N_k^{\rm bd}(X)\geq m_k/8\}
\cup\{N^{\mathrm{exc}}(X)\geq m_k/16\}.
\]
If $\mathcal E_k(X)$ does not occur, $N_{W'}<m_k/4$, $N^G<m_k/48$ and $\ell_n<m_k/16$, the right-hand side of \eqref{eq:vote-perturbation} is less than
$m_k(1/2+3/16+3/48+1/8+1/16)=15m_k/16$, so the vote is correct. New errors are therefore contained in
\[
(\widehat{\mathcal V}_k\cap \mathcal X_{\mathrm{exc}})\cup G_n\cup\mathcal D_k\cup\mathcal W_k\cup\mathcal G_k\cup\mathcal L_k,
\]
where, among vertices of $A\cap B_n$,
$\mathcal D_k=\{u_X\in T_k,\mathcal E_k(X)\}$,
$\mathcal W_k=\{u_X\in T_k,N_{W'}(X)\geq m_k/4\}$,
$\mathcal G_k=\{N^G(X)\geq m_k/48\}$ and
$\mathcal L_k=\{\ell_n(X)\geq m_k/16\}$.

Every $Y\in W'_k$ has true age at least $a_k/2$ by Lemma~\ref{lem:rank-vote}. Counting incidences towards $T_k$, then applying \eqref{eq:target-tail},
\[
\frac{|\mathcal W_k|}{N_n^A}
\leq\frac4{m_kN_n^A}\sum_{Y\in W'_k}D_k(Y)
\leq c_va_k^\gamma e_{k,n}+C_2a_k^\gamma e^{-\Lambda_*\lambda}+o_{\mathbb P}(1).
\]
The score defining $\mathcal D_k$ is translation-covariant and bounded. Lemma~\ref{lem:oracle-PA}, \eqref{eq:misplaced-neighbours} and $|T_k|=7a_k/4$ give, by ergodicity,
$|\mathcal D_k|/N_n^A\leq C_1a_ke^{-c_1\lambda}+o_{\mathbb P}(1)$.
Finally,
\[
|\mathcal G_k|\leq\frac{48}{v\lambda}\sum_{Y\in G_n}D(Y),\qquad
|\mathcal L_k|\leq\frac{16}{v\lambda}\sum_{X\in A\cap B_n}\ell_n(X).
\]
These two terms and $|G_n|$ are $o_{\mathbb P}(N_n^A)$ by Lemma~\ref{lem:vote-boundary}. Adding the existing errors gives \eqref{eq:hidden-recursion}. There are finitely many remainder terms, since $\lambda$ and $K$ are fixed before $n\to\infty$.
\end{proof}

Let $s^{(1)},s^{(2)}$ be the outputs of the two passes in Definition~\ref{def:hidden-algorithm}. Orient the second by the sign of their correlation on colour 2, then keep the first on colours 1 and 2 and the oriented second on colour 3.

To check this alignment, orient each output to minimise its errors. If each has error proportion at most $e$ on its population, their total error count is at most $(4e/3+o(1))N_n$. The intersection contains $(1/3+o(1))N_n$ vertices. For $e<1/8$ with a strict margin, fewer than half these vertices can therefore have disagreeing oriented labels. The observed correlation recovers the relative sign. The error count of the merged output is at most the sum of the two error counts; a union bound combines their guarantees without assuming independence.

To obtain deterministic estimators at fixed parameters, realise the colours through a measurable function $\chi(\omega,\cdot)$ of a seed $\omega$ uniform on $[0,1]$, independent of the graph, providing independent colours to lexicographically ordered vertices. If $F_n$ is the failure event and $f_n(\omega)=\Prob(F_n\mid\omega)$, then $\int_0^1f_n(\omega)\dd\omega=\Prob(F_n)$. There is therefore a seed $\omega_n$ such that $f_n(\omega_n)\leq\Prob(F_n)+1/n$; fixing it gives a deterministic estimator whose failure probability still tends to zero. This choice depends on the model and the required accuracy: it is an existence statement, without an efficient search for a universal colour sequence.

\section{Rates for the overlap deficit}\label{sec:rate}
Voting does not multiply errors by a fixed constant at every step. In both settings, the recursions give total amplification bounded by
\begin{equation}\label{eq:bounded-amplification}
\Gamma_*:=\prod_{k\geq0}(1+c_va_k^\gamma)
\leq\exp\!\left(\frac{c_va_0^\gamma}{1-2^{-\gamma}}\right)<\infty.
\end{equation}
This constant depends neither on $K$ nor on $\lambda$, but depends on the model and on the setting through $p$ and $v$. The series $\sum_ka_k=2a_0$ and $\sum_ka_k^\gamma=a_0^\gamma/(1-2^{-\gamma})$ also allow local failures to be summed. The final rate is therefore the initialisation rate, plus exponentially small errors and the fraction left below the last cutoff.

\subsection{Observed ages and known parameters}
Proposition~\ref{prop:exp-seed} gives $e_{0,n}\leq Ce^{-c_0\lambda}+o_{\mathbb P}(1)$, and iterating \eqref{eq:observed-recursion} then yields $e_{K,n}\leq Ce^{-c\lambda}+o_{\mathbb P}(1)$, with constants independent of $K$. Below the cutoff $\zeta=a_02^{-K}$, the proportion assigned $+1$ tends to $\zeta$, so, since the overlap deficit is at most twice the error proportion in a given orientation,
\begin{equation}\label{eq:rate-balance}
1-\ov_n\leq C(\zeta+e^{-c\lambda})+o_{\mathbb P}(1).
\end{equation}
Choose $K=\lceil\lambda/\log2\rceil$, so that $\zeta\leq e^{-\lambda}$. Increasing $C$ and decreasing $c$ to leave a strict margin for the remainder in probability gives the first assertion of Theorem~\ref{th:main-rates}. A sufficient intensity for deficit $\varepsilon$ satisfies
\begin{equation}\label{eq:lambda-observed}
\lambda_{\rm suff}(\varepsilon)\leq\max\{\lambda_0,c^{-1}\log(C/\varepsilon)\}.
\end{equation}

{
\paragraph{When intensity is estimated.}
\phantomsection\label{rem:lambda-observed}
Proposition~\ref{prop:intensity-estimation} gives the observable
estimator $\widehat\lambda_n=N_n/n$, with variance $\lambda/n$
and almost-sure convergence to $\lambda$. The agnostic procedure
already incorporates this estimate through its dyadic scale
$\Lambda_n$, as treated in the proof of
Proposition~\ref{prop:hidden-seed}; the question here is whether
the procedure with observed ages retains its exponential rate
when only $\lambda$ is unknown.
For a fixed accuracy $\varepsilon$, choose $K$ so that
$a_02^{-K}\leq\varepsilon/8$, independently of $\lambda$.
Substituting $\widehat\lambda_n$ into the local thresholds would
then require a local event on which all decisions are correct
simultaneously for every intensity in some fixed interval
$[(1-\delta)\lambda,(1+\delta)\lambda]$, with $\delta>0$.
Indeed, the failure probability is at most the probability that
$\widehat\lambda_n$ leaves this interval plus that of failure of
the uniform local event, without conditioning the local field
on the global estimator. An exponential bound for this local
event would preserve the form of \eqref{eq:lambda-observed},
with possibly different constants; we do not claim this
additional stability result here.

}
\subsection{Universal reconstruction without ages or profiles}
Fix $\delta=\delta_1$ in Lemma~\ref{lem:misplaced-neighbours}, and use $\mathcal X_{\mathrm{exc}}$ from \eqref{eq:rank-exceptions}. Estimated slices are disjoint, so their age failures are counted only once:
\[
\sum_{k<K}\frac{|\widehat{\mathcal V}_k\cap \mathcal X_{\mathrm{exc}}|}{N_n^A}
\leq\frac{|\mathcal X_{\mathrm{exc}}\cap B_n|}{N_n^A}
\leq2e^{-\kappa_{\delta/3}\lambda}+o_{\mathbb P}(1).
\]
Recursion~\eqref{eq:hidden-recursion} and \eqref{eq:bounded-amplification} then give
\[
e_{K,n}\leq\Gamma_*\varepsilon_{\rm init}(\lambda)+Ce^{-c\lambda}+o_{\mathbb P}(1).
\]
Outside $\mathcal X_{\mathrm{exc}}\cup G_n$, Lemma~\ref{lem:rank-vote} shows that a vertex with rank below $\zeta$ has true age below $2\zeta$. The fraction left unprocessed is thus at most $2\zeta+2e^{-\kappa_{\delta/3}\lambda}+o_{\mathbb P}(1)$. Now $K=\lceil\Lambda_n\rceil\geq\lambda/2$ with probability tending to one, so $\zeta\leq a_02^{-\lambda/2}$.

For each pass, the total error proportion is therefore at most
$\Gamma_*\varepsilon_{\rm init}(\lambda)+Ce^{-c\lambda}+o_{\mathbb P}(1)$. Alignment on colour 2 multiplies this error bound by at most $4/3$, and conversion to overlap by two. Finally,
\begin{equation}\label{eq:universal-rate}
1-\ov_n\leq\frac83\Gamma_*\varepsilon_{\rm init}(\lambda)+Ce^{-c\lambda}+o_{\mathbb P}(1).
\end{equation}
{The proof of Proposition~\ref{prop:hidden-seed}
gives $\varepsilon_{\rm init}(\lambda)=o(e^{-b(\lambda)})$,
where $b(\lambda)=\lambda/(\log\lambda)^{2d+3}$.
Since $\Gamma_*$ is independent of $\lambda$ and
$e^{-c\lambda}=o(e^{-b(\lambda)})$, the deterministic terms in
\eqref{eq:universal-rate} are at most $\tfrac12 e^{-b(\lambda)}$
for sufficiently large $\lambda$. At each such fixed intensity,
the remainder tends to zero in probability as $n\to\infty$;
the strict margin therefore gives
\[
\Prob\!\left(\ov_n\geq
1-\exp\!\left[-\frac{\lambda}{(\log\lambda)^{2d+3}}\right]\right)
\longrightarrow1.
\]
This proves the second assertion of Theorem~\ref{th:main-rates}.} The intensity threshold depends on the model; the algorithm is given neither the required accuracy of survival functions nor a minimum volume common to all models.

This bound is faster than every power of $\lambda^{-1}$, but does not give a bound $Ce^{-c\lambda}$ for the universal procedure. The logarithmic loss comes from initialisation: the selected fraction and cell volume decrease, reducing mean occupancy by the factor $w\ell^d$, while dependence range, measured in cells, grows as $R/\ell$. It does not come from accumulated voting errors. For large $L_\varepsilon=\log(1/\varepsilon)$, an intensity of order
\[
\lambda_{\rm suff}(\varepsilon)
=O_{\text{model}}\!\left(L_\varepsilon(\log L_\varepsilon)^{2d+3}\right)
\]
is sufficient. Indeed, for $\lambda=A L_\varepsilon(\log L_\varepsilon)^{2d+3}$, we have $\log\lambda\leq2\log L_\varepsilon$ for large enough $L_\varepsilon$; choosing $A>2^{2d+3}$ makes $\lambda/(\log\lambda)^{2d+3}>L_\varepsilon$.

\begin{corollary}[Hidden ages with scales adapted to the model]\label{cor:tuned-hidden}
When profiles and structural parameters are known, fixed choices of $w,\ell,R$ in the rank-based procedure give, without observing ages, constants $C,c,\lambda_0>0$ such that
\[
\Prob(\ov_n\geq1-Ce^{-c\lambda})\longrightarrow1\qquad(\lambda\geq\lambda_0).
\]
Intensity may still be estimated by $N_n/n$. Deterministic estimators with this guarantee exist for a fixed model.
\end{corollary}
\begin{proof}
First choose $R\geq R_0$ sufficiently large, then $\ell$ and the partition width $w$ sufficiently small for the conditions and $V\leq D_*/8$ in Lemma~\ref{lem:universal-margin}. More precisely, its bounds allow $T(R/2)$, $\omega(4\sqrt d\ell)$ and $w$ to be made small separately, while ensuring $R\geq2R_c$. Fix these values and $h=c_hw$. Keep ranks, quartiles, the occupancy threshold $\max\{p\Lambda_nw\ell^d,1\}$ and $K=\lceil\Lambda_n\rceil$. The defect probability \eqref{eq:universal-defect} is then $Ce^{-c\lambda}$ and the range $M$ stays fixed. Quartiles remain localised with probability tending to one as volume grows. The proof of Proposition~\ref{prop:hidden-seed} therefore gives exponentially accurate initialisation, and \eqref{eq:universal-rate} transfers this order to the final output. The averaging argument at the end of Section~\ref{sec:ages} gives the stated deterministic existence. Here the scales depend on the model; this does not claim a single fixed choice valid for all profiles.
\end{proof}

\begin{remark}[Scope of the results]
Universal reconstruction does not require separate estimates of $\gamma$ or the profiles, nor a distinction between $\beta$ and a dilation of their argument: replacing $\beta$ by $c\beta$ and $f(t)$ by ${f(ct)}$ preserves connections. It uses their monotonicity, unlike the distinguishability test, which accepts measurable profiles but still assumes known $\gamma$. The two results therefore have different scopes. Finally, all our statements take the volume limit before the high-intensity limit; they provide neither an optimal threshold, nor uniform guarantees as the signal vanishes, nor a volume convergence rate common to all models.
\end{remark}

\appendix
\section{Proof of the fluctuation estimate  (Lemma  \ref{lem:fluctuations})}\label{app:fluctuations}
Fix $\mathsf M,R$ and $\xi^{\mathrm{norm}}=\xi/c$. The cube $C$ is a union of $V$ unit cells $C_z=z+[0,1)^d$; set $Y_z=\sum_{x\in\eta\cap C_z}\xi^{\mathrm{norm}}(x,G)$   and $S_C=\sum_zY_z$. Campbell gives $\E_{\mathsf M}Y_z=\lambda\mu$ and $\E_{\mathsf M}S_C=\lambda V\mu$. We control the fourth moment, whose $V^{-2}$ decay will make error probabilities summable.

\paragraph{Local contribution and dependence.}
The variable $Y_z$ depends only on the induced graph in $H_z=z+[-R,1+R]^d$. With $K_z=N(H_z)$ of mean $\theta_R=\lambda(1+2R)^d$, we have $0\leq Y_z\leq K_z^3$, since degrees contribute at most $K_z^2$ and triangles at most $K_z^3$. The factorial moments of a Poisson variable are $\theta^j$, so its twelfth moment is a polynomial in $\theta$, giving $\E K^{12}\leq C_{12}(1+\theta)^{12}$. Jensen and $|s-t|^4\leq8(|s|^4+|t|^4)$ yield
\begin{equation}\label{eq:local-fourth}
\E_{\mathsf M}|Y_z-\E_{\mathsf M}Y_z|^4\leq16\E_{\mathsf M}Y_z^4
\leq C_{\lambda,d}(1+R)^{12d}=:A_R.
\end{equation}
Colour the indices according to their $k_R=q_R^d$ residue classes modulo $q_R=\lceil2R\rceil+2$. Within a class, the $H_z$ are disjoint and the $Y_z$ independent. For the centred sum $Z_b$ over the $m_b$ cells of colour $b$, expansion of the fourth moment retains only terms of orders four and two-two, hence
\[
\E_{\mathsf M}Z_b^4=\sum_z\E_{\mathsf M}(Y_z-\E_{\mathsf M}Y_z)^4
+6\sum_{z<w}\E_{\mathsf M}(Y_z-\E_{\mathsf M}Y_z)^2\E_{\mathsf M}(Y_w-\E_{\mathsf M}Y_w)^2
\leq3m_b^2A_R.
\]
Different colours may be dependent; the inequality $|\sum_bz_b|^4\leq k_R^3\sum_b|z_b|^4$, with $\sum_bm_b^2\leq V^2$, suffices to obtain
\begin{equation}\label{eq:numerator-fourth}
\E_{\mathsf M}\left|\frac{S_C}{\lambda V}-\mu\right|^4
\leq\frac{3k_R^3A_R}{\lambda^4V^2}
\leq\frac{C_{\lambda,d}(1+R)^{15d}}{V^2}.
\end{equation}

\paragraph{Random denominator.}
With $X=S_C/(\lambda V)$ and $Z=N(C)/(\lambda V)$, the fourth Poisson moment gives
\begin{equation}\label{eq:denominator-fourth}
\E_{\mathsf M}|Z-1|^4=(\lambda V)^{-3}+3(\lambda V)^{-2}
\leq C_\lambda V^{-2}.
\end{equation}
On $|X-\mu|\leq b/4$ and $|Z-1|\leq b/4$, where $0<b\leq1$, we have $Z\geq3/4$ and
$|S_C/N(C)-\mu|\leq(|X-\mu|+\mu|Z-1|)/Z\leq2b/3$.
Projection does not increase this error, so the event in Lemma~\ref{lem:fluctuations}, including $N(C)=0$, is contained in the union of the two complementary events. Markov applied to \eqref{eq:numerator-fourth} and \eqref{eq:denominator-fourth} gives the bound $C_{\lambda,d}(1+R)^{15d}/(V^2b^4)$, uniformly over profiles and without independence between numerator and denominator.

\section{Explicit construction of a finite-comparison test}
\label{app:finite-references}
Compactness does not compute the minimum; we therefore construct a finite list of profiles and simulate their moments with controlled error. The parameters $(\lambda,d,\gamma)$ are known, $\lambda,\gamma$ are computable with certified precision, and $\beta=1$. Positions are given exactly in a representation allowing geometric comparisons; auxiliary draws are independent of the graph. We establish almost-sure termination and consistency, without claiming efficiency or robustness to location errors.

Fix the first $J$ scores of Section~\ref{sec:unknown} and a tolerance $\varepsilon\in(0,1]$. We seek an explicitly defined finite family $\mathcal S(J,\varepsilon)$ such that
\begin{equation}\label{eq:finite-net}
\forall \varphi\in\F\quad\exists s\in\mathcal S(J,\varepsilon):
\quad \max_{j\leq J}|\mu_j(\varphi)-\mu_j(s)|\leq\varepsilon.
\end{equation}
Only moments are approximated, using the regularity of their integration kernels rather than that of the unknown measurable profile.

Write $F(t)=\varphi(e^t)$. In logarithmic coordinates, each normalised moment has the form
\begin{equation}\label{eq:log-moment}
\mu_j(\varphi)=\int_{\R^{p_j}}W_j(\mathbf t)
\prod_{i=1}^{p_j}F(t_i)\dd\mathbf t,
\qquad p_j=1\text{ or }3,
\end{equation}
Here $W_j$ is the density of the logarithms of connection arguments, multiplied by the moment's normalising factor: $W_j\geq0$ and $\|W_j\|_1\leq1$. It is not the logarithm of a density.

\subsection{Densities in logarithmic coordinates}
For an integrable function $w$ whose distributional derivatives are finite measures, write
$\mathcal V(w)=\sum_i |D_iw|(\R^p)$.
This variation accounts for both derivatives inside regularity regions and jumps at their boundaries. It concerns the integration kernel, not the unknown profile.

\begin{lemma}[Variation of moment kernels]\label{lem:log-bv}
The kernels in \eqref{eq:log-moment} can be chosen so that
\[
\mathcal V(W_j)\leq B_\gamma,
\qquad B_\gamma=\frac{36}{\gamma(1-\gamma)},
\]
independently of the radius, $\lambda$ and the null profile.
\end{lemma}
\begin{proof}
Fix three distinct positions. On one strict ordering of log-ages, the density is $h_{\log}(\mathbf v)=e^{v_1+v_2+v_3}\Ind_{\{v_1<v_2<v_3<0\}}$. Interior derivatives contribute $1/2$ to its variation, while the faces $v_3=0$, $v_1=v_2$, $v_2=v_3$ contribute $1/2$, $1/3$, $2/3$, respectively; thus $\mathcal V(h_{\log})=2$, by integration by parts and truncation.

The map to logarithms of arguments is affine, with matrix $M_\gamma$ from Lemma~\ref{lem:density}. For the image $q$ of a density under $M$, the differentiation rule, extended by regularisation, gives $\mathcal V(q)\leq\|M^{-1}\|_\infty\mathcal V(h_{\log})$, where the norm is the maximum absolute row sum. Cofactors of $M_\gamma$ are at most one in absolute value, so $\|M_\gamma^{-1}\|_\infty\leq3/[\gamma(1-\gamma)]$. Summing over the six age orderings gives $B_\gamma$; spatial constants merely translate densities.

For a degree, \eqref{eq:convolution} expresses the density of $-\log g(U,V)$ as the convolution of exponential densities with rates $2$ and $1/\gamma$. The first has total variation $4$, which convolution does not increase. Finally, the kernels $W_j$ are mixtures of these densities over positions, multiplied by $\lambda\omega_dR_j^d/(1+\lambda\omega_dR_j^d)$ or $(\lambda\omega_dR_j^d)^2/[1+(\lambda\omega_dR_j^d)^2]$, both at most one. The variational definition and Tonelli preserve the bounds.
\end{proof}

\begin{lemma}[Averages on a logarithmic grid]\label{lem:cell-average}
Partition $\R$ into intervals of length $h>0$ and let $\mathbf{P}_h $ denote cell averaging, using the product partition on $\R^p$. If $W\in L^1(\R^p)$ and $\mathcal V(W)<\infty$, then
$\|W-\mathbf{P}_h W\|_1\leq h\mathcal V(W)$.
Consequently, for every function $F:\R\to[0,1]$,
\[
\left|\int W\prod_{i=1}^pF(t_i)
-\int W\prod_{i=1}^p(\mathbf{P}_h F)(t_i)\right|
\leq2h\mathcal V(W).
\]
\end{lemma}
\begin{proof}
On an interval $I$ of length $h$, Jensen gives
$\int_I|w-w_I|\leq h^{-1}\int_{I^2}|w(x)-w(y)|\dd x\dd y\leq h|Dw|(I)$.
Summing and averaging successively in each coordinate, using contraction in $L^1$, gives the first bound. Extension to functions of bounded variation follows by regularisation. On each product cell, $\mathbf{P}_h W$ is constant, so
$\int(\mathbf{P}_h W)\prod_iF(t_i)=\int(\mathbf{P}_h W)\prod_i(\mathbf{P}_h F)(t_i)$.
Replacing $W$ with $\mathbf{P}_h W$ on both sides costs at most $2\|W-\mathbf{P}_h W\|_1$, giving the second bound.
\end{proof}

\subsection{A finite grid with certified error}
Write $R_{\min}=\min_{j\leq J}R_j>0$ and $R_{\max}=\max_{j\leq J}R_j$. Choose
\begin{equation}\label{eq:explicit-grid}
a=\varepsilon/48,\quad
t_-=R_{\min}^d a^2,\quad t_+=(2R_{\max})^d,
\quad N=\lceil12/\varepsilon\rceil.
\end{equation}
Choose rationals $0<\gamma_-<\gamma<\gamma_+<1$ and an integer $B>36/[\gamma_-(1-\gamma_+)]$. For integers $k\geq0$, $\ell\geq1$ such that $2^{-k}<t_-$ and $2^\ell>t_+$, set $A_-=-k$, $A_+=\ell$; then $e^{A_-}<t_-$ and $e^{A_+}>t_+$, since $e>2$. These choices use rational comparisons for a rational tolerance $\varepsilon$; otherwise, replace it with a smaller rational. Define
\[
h=\frac{\varepsilon}{16B},\qquad m=\left\lceil\frac{A_+-A_-}{h}\right\rceil,
\qquad I_k=[A_-+kh,A_-+(k+1)h),\quad 0\leq k<m.
\]

The family $\mathcal S(J,\varepsilon)$ contains all profiles equal to $v_k/N$ on $\exp(I_k)$, where $(v_0,\ldots,v_{m-1})\in\{0,\ldots,N\}^m$, and zero elsewhere. It contains exactly $M=(N+1)^m$ profiles, enumerated in lexicographic order of the integer vectors. It does not depend on the generating profile.

\begin{proposition}[Uniform approximation of moments]\label{prop:explicit-net}
The family above satisfies \eqref{eq:finite-net}, even with $\varepsilon/2$ in place of $\varepsilon$.
\end{proposition}
\begin{proof}
Truncate $F(t)=\varphi(e^t)$ to the union of the $I_k$. No relevant upper argument is removed, since $\kappa_{ij}\leq(2R_{\max})^d$. At the lower end, set $\rho=R_{\min}a^{1/d}$. If all three ages are at least $a$ and all three spatial distances at least $\rho$, then $g\geq a$ and $\kappa_{ij}\geq a\rho^d=t_-$. Failure of the age condition costs at most $3a$; each distance $|X|,|Y|,|X-Y|$ is below $\rho$ with probability at most $(\rho/R)^d\leq a$, conditioning on $X$ for the last. Truncation therefore costs at most $6a=\varepsilon/8$ for a triangle and $3a$ for a degree.

Extend the grid to the whole real line and replace the truncated profile by its cell averages. Since its support boundaries are cell boundaries, only the $m$ retained cells remain, and Lemmas~\ref{lem:log-bv}--\ref{lem:cell-average} bound the error by $2hB_\gamma\leq\varepsilon/8$. Finally, round averages down to multiples of $1/N$; the difference between products of $p$ factors in $[0,1]$ is at most the sum of their differences, so this costs at most $3/N\leq\varepsilon/4$. The resulting candidate belongs to $\mathcal S(J,\varepsilon)$, and the total error is at most $\varepsilon/2$ for all retained moments.
\end{proof}

\subsection{Computing references by controlled simulation}
For a candidate $s$, draw independent uniform ages and, according to the score, one or two independent uniform positions in $B(0,R)$. Set
\[
Z_s=\begin{cases}s(\kappa_{01}),&\xi_j=d_R,\\
s(\kappa_{01})s(\kappa_{02})s(\kappa_{12}),&\xi_j=t_R,
\end{cases}
\qquad
A_j=\frac{(\lambda\omega_dR^d)^{p_j'}}{1+(\lambda\omega_dR^d)^{p_j'}},
\]
with $p_j'=1$ for a degree and $p_j'=2$ for a triangle. Palm formulas give $\mu_j(s)=A_j\E Z_s$, with $0\leq Z_s\leq1$; we therefore simulate configurations of two or three vertices, rather than a large graph.  { $ A_{ j}$ involved in the algorithm? $ \lambda $ required?}

\begin{lemma}[Simultaneous accuracy of simulated moments]\label{lem:reference-MC}
For a family of $M$ candidates, $J$ scores, rational tolerance $\varepsilon\in(0,1]$ and $\delta\in(0,1)$, take
\begin{equation}\label{eq:MC-size}
S\geq\frac{8}{\varepsilon^2}\log\frac{2JM}{\delta}.
\end{equation}
Using $S$ independent draws per score, replace $\E Z_s$ by its empirical average and $A_j$ by a rational $\widetilde A_j\in[0,1]$ with error at most $\varepsilon/4$. The resulting rational values $\widehat m_j(s)$ satisfy
\[
\Prob\left(\max_{j\leq J,\,s\in\mathcal S}
|\widehat m_j(s)-\mu_j(s)|>\varepsilon/2\right)\leq\delta.
\]
\end{lemma}
\begin{proof}
Hoeffding bounds the probability of an average error exceeding $\varepsilon/4$ by $2e^{-S\varepsilon^2/8}$ for each pair $(j,s)$. A union bound over the $JM$ pairs costs at most $\delta$, and approximation of $A_j$ adds at most $\varepsilon/4$. Configurations may be shared between candidates for the same score, since the union bound does not require independence.
\end{proof}

\noindent\textit{Computational model.} An integer $S$ satisfying \eqref{eq:MC-size} is chosen using a certified bound; positions are sampled by rejection from a cube. Coordinates and arguments are evaluated by intervals until profile cells are identified: their boundaries have probability zero, so each comparison terminates almost surely. Fixed floating-point precision would require additional rounding-error control.

\subsection{Complete routine and decision}
To make the threshold rational as well, use, for $L\geq4$,
\begin{equation}\label{eq:operational-schedule}
J_L=\lfloor\log_2 L\rfloor,\qquad
\varepsilon_L=2^{-k_L},\quad
k_L=\min\{k\in\N:L\leq2^{8k}\},\qquad
\delta_L=L^{-2}.
\end{equation}
Thus $\tfrac12L^{-1/8}<\varepsilon_L\leq L^{-1/8}$. The inner cube remains $C_L$ from \eqref{eq:test-schedule}. We have $2^{-J_L}\leq\varepsilon_L$: for $4\leq L<8$, the two values are $1/4$ and $1/2$, respectively, while for $L\geq8$ this follows from $2^{-J_L}\leq2/L<\tfrac12L^{-1/8}$.

Construct $\mathcal S_L=\mathcal S(J_L,\varepsilon_L)$ and the references $\widehat m_j(s)$ from Lemma~\ref{lem:reference-MC}. Also approximate each projected observed statistic $\widehat\mu_{L,j}$ by a rational $\widetilde\mu_{L,j}\in[0,1]$ with error at most $\varepsilon_L/2$, then set
\begin{equation}\label{eq:operational-test}
\widehat D_L=\min_{s\in\mathcal S_L}
\sum_{j=1}^{J_L}2^{-j}
|\widetilde\mu_{L,j}-\widehat m_j(s)|,
\qquad
\widehat T_L=\Ind_{\{\widehat D_L>4\varepsilon_L\}}.
\end{equation}
The final comparison is rational, including equality, and integer counts and computable normalisations provide the required approximations; exact geometric comparisons terminate almost surely, and reduce to squared distances for rational coordinates.

In practice, derive $J_L,\varepsilon_L,\delta_L$ from the window, construct the finite grid and its simulated references, then compute observed degrees and triangles. Looping over candidates gives the minimum in \eqref{eq:operational-test}: answer two communities if it exceeds $4\varepsilon_L$, one otherwise, and zero by convention if $L<4$. These references can be prepared independently of the graph.

\paragraph{Guarantee of the procedure.}
Write $\overline P$ for the joint law of the graph with law $P$ and the independent auxiliary draws. Consistency now refers to this joint law.

\begin{theorem}[Explicit finite-comparison test]\label{th:operational}
Under the assumptions of Theorem~\ref{th:disc}, procedure~\eqref{eq:operational-test} terminates almost surely in the computational model above. There is a summable sequence
\[
b_L=\min\{1,C_{\lambda,d}J_L(1+J_L)^{15d}L^{1/2-2d}\}
\]
such that
\begin{align*}
\sup_{\varphi\in\F}\overline {\mathbf{ P}}_\varphi(\widehat T_L=1)&\leq b_L+\delta_L,\\
\overline Q(\widehat T_L=0)&\leq b_L+\delta_L
\quad\text{if }\Delta_{\Q}^*>7\varepsilon_L.
\end{align*}
For each fixed model, answers are eventually correct almost surely. The same holds on $B_n$ with $L=L(n)$ from \eqref{eq:test-volume}, reusing the same routine and draws for a given index $L$.
\end{theorem}
\begin{proof}
Let $\mathcal E_L(P)$ be the event that all exact observed moments are within $\varepsilon_L$ of their expectations, and $\mathcal F_L$ the event that all references are within $\varepsilon_L/2$ of theirs. Lemma~\ref{lem:fluctuations}, a union bound over $j\leq J_L$, and $\varepsilon_L>\tfrac12L^{-1/8}$ give $P(\mathcal E_L(P)^c)\leq b_L$, while Lemma~\ref{lem:reference-MC} gives $\Prob(\mathcal F_L^c)\leq\delta_L$.

Under $\mathbf{ P}_\varphi$, Proposition~\ref{prop:explicit-net} provides a candidate whose moments each differ from those of $\varphi$ by at most $\varepsilon_L$. On $\mathcal E_L(\mathbf{ P}_\varphi)\cap \mathcal F_L$, the statistical and grid errors and the two computational errors give
\[
\widehat D_L\leq\varepsilon_L+\varepsilon_L+\varepsilon_L/2+\varepsilon_L/2
=3\varepsilon_L<4\varepsilon_L,
\]
since the weights sum to less than one. Under $Q$, every candidate is a null profile, so $d_{*,J_L}(\mu(Q),\mu(s))\geq\Delta_{\Q}^*-2^{-J_L}$. Statistical and computational errors subtract at most $2\varepsilon_L$, hence
\[
\widehat D_L\geq\Delta_{\Q}^*-2^{-J_L}-2\varepsilon_L
\geq\Delta_{\Q}^*-3\varepsilon_L>4\varepsilon_L
\]
if $\Delta_{\Q}^*>7\varepsilon_L$. A union bound on failure events gives the stated risks. Since $\sum_Lb_L$ and $\sum_L\delta_L$ converge, Borel--Cantelli and $\varepsilon_L\to0$ give almost-sure consistency; reindexing by $L(n)$ transfers it to $B_n$. Finally, enumerations are finite and each comparison terminates almost surely in the computational model above.
\end{proof}

Enumerating $(N+1)^m$ profiles may be prohibitive; reducing this list requires a certificate for the minimisation error. The threshold is also conservative, since it cannot be exceeded when $4\varepsilon_L\geq1$; the result establishes a consistent finite procedure, without an efficiency guarantee at moderate sizes.

The construction assumes known structural parameters; substitution of estimated intensity into the ideal test is treated at the end of Section~\ref{sec:test}.

\section{Ranks and certification of universal initialisation}\label{app:universal}
This appendix justifies the procedure's two substitutions: degrees become age ranks, and connection means are replaced by quartiles. All constants $C,c>0$ may depend on the fixed model; they do not enter the algorithm. We work with a deterministic value $\Lambda\in[\lambda/2,\lambda]$, with scales fixed before the volume limit. The estimates below are uniform on this interval because $w\asymp(\log\lambda)^{-1}$, $\ell\asymp(\log\lambda)^{-1}$ and $R\asymp\log\log\lambda$, not because the whole procedure must be monotone.

\subsection{From degree to age rank}\label{app:universal-ranks}
\begin{proof}[Proof of Lemma~\ref{lem:rank-vote}]
First compute $S_n'$ as $S_n$, but using infinite-graph degrees. The two degrees agree outside $G_n$, so
$\sup_j|S_n(j)-S_n'(j)|\leq|G_n|/N_n^A=o_{\mathbb P}(1)$ by Lemma~\ref{lem:vote-boundary}. For each integer $j$, ergodicity applied to the bounded score $\Ind_{\{D_B(X)\geq j\}}$ gives $S_n'(j)\to S(j)$ almost surely, simultaneously for all integers. Given $\epsilon>0$, choose $j_\epsilon$ such that $S(j_\epsilon)\leq\epsilon$; for $j\geq j_\epsilon$, monotonicity gives
\[
|S_n'(j)-S(j)|\leq S_n'(j_\epsilon)+S(j_\epsilon),
\]
and only finitely many indices remain. This proves uniform convergence. The argument is identical for $S_{R,n}$; only vertices within distance $R$ of the boundary can have a changed truncated degree. Their count is random, but its expectation is $O(\lambda Rn^{1-1/d})$, so their proportion tends to zero in probability. Independence between degrees is not needed.

For relative accuracy, write $\mu(t)=\rho_{\mathrm{cal}}\lambda\bar I\Theta(t)$ and $D=D_B(X)$ under Palm at age $s$. We have $D\sim\operatorname{Pois}(\mu(s))$, $\mu(t)\leq \rho_{\mathrm{cal}}\lambda\bar I C_\gamma t^{-\gamma}$ and, for $x>1$ within the age domain,
\[
\Theta(s)-\Theta(xs)\geq s^{-\gamma}\frac{x-1}{x},\qquad
\Theta(s/x)-\Theta(s)\geq s^{-\gamma}\log x.
\]
These inequalities follow from $\Theta'(t)=-t^{-\gamma-1}$. We use the Poisson bounds
$\Prob(Z-\mu\geq t)\leq e^{-t^2/(2(\mu+t))}$ and
$\Prob(Z-\mu\leq-t)\leq e^{-t^2/(2\mu)}$.
If
\[
\lambda\geq\lambda_1(\rho):=\frac{\log(4/\rho)+1/\gamma}{\kappa_\rho},
\]
then $e^{-\kappa_\rho\lambda s^{-\gamma}}\leq\rho s/4$ for every $s\leq1$, since, with $x=s^{-\gamma}\geq1$, $\log(1/s)=\gamma^{-1}\log x\leq x/\gamma$.

For overestimation, $(1+\rho)s>1$ is impossible. Otherwise, set $s_+=(1+\rho/2)s$ and $b_+=(\mu(s)-\mu(s_+))/2\geq \rho_{\mathrm{cal}}\lambda\bar I s^{-\gamma}\rho/6$. For $j\geq\mu(s_+)+b_+$, stochastic monotonicity of Poisson laws implies
\[
S(j)\leq s_++\exp\!\left[-\frac{b_+^2}{2\mu(s)}\right]
\leq s_++\rho s/2\leq(1+\rho)s.
\]
Thus $S(D)>(1+\rho)s$ requires $D<\mu(s)-b_+$, with probability at most $e^{-2\kappa_\rho\lambda s^{-\gamma}}$.
For underestimation, set $s_-=s/(1+\rho/2)$ and $b_-=(\mu(s_-)-\mu(s))/2\geq \rho_{\mathrm{cal}}\lambda\bar I s^{-\gamma}\rho/6$. If $j\leq\mu(s_-)-b_-$, then
\[
S(j)\geq s_-\left(1-e^{-b_-^2/(2\mu(s_-))}\right)
\geq\frac{1-\rho/4}{1+\rho/2}s\geq\frac{s}{1+\rho}.
\]
We used $\mu(s_-)\leq2\rho_{\mathrm{cal}}\lambda\bar I C_\gamma s^{-\gamma}$ and $\rho\leq1$. The remaining event requires $D>\mu(s)+b_-$ and has probability at most $e^{-\kappa_\rho\lambda s^{-\gamma}}$. Their union gives \eqref{eq:rank-relative}.

It remains to check the inclusions needed for voting. On $\Omega_n$, outside $\mathcal X_{\mathrm{exc}}\cup G_n$, with $\eta=\delta\zeta/8$, we have
\begin{equation}\label{eq:additive-rank-bracket}
\frac{u_X}{1+\delta/3}-\eta\leq\widehat u_X\leq(1+\delta/3)u_X+\eta.
\end{equation}
For $k<K$, $\eta\leq\delta a_k/16$, so
\[
\frac{a_k/2-\eta}{1+\delta/3}\geq\frac{a_k}4,
\quad (1+\delta/3)(a_k+\eta)<2a_k,
\quad\frac{a_k-\eta}{1+\delta/3}\geq\frac{a_k}2.
\]
These are the first two inclusions. If the threshold $a_k$ is crossed incorrectly, \eqref{eq:additive-rank-bracket} gives either $u_X\leq(1+\delta/3)(a_k+\eta)\leq(1+\delta)a_k$, or $u_X\geq(a_k-\eta)/(1+\delta/3)\geq a_k/(1+\delta)$. Finally, $\widehat u_X<\zeta$ gives $u_X\leq(1+\delta/3)(\zeta+\eta)<2\zeta$, while $\widehat u_X\geq a_0$ gives $u_X\geq(a_0-\eta)/(1+\delta/3)\geq2a_0/3$. Since $K$ is fixed before $n\to\infty$, the accuracy $\eta$, however small, is fixed: uniform convergence suffices for $\Prob(\Omega_n)\to1$.
\end{proof}

\begin{lemma}[Local rank accuracy]\label{lem:local-rank}
There exist $R_0,c_*,c'>0$, depending on the model, such that for $R\geq R_0$,
\[
m_R(s')-m_R(s)\geq c_*(s-s'),\qquad a_0/16\leq s'\leq s\leq1.
\]
For $0<h\leq a_0/16$ and $\lambda\geq\log(2/h)/(c'h^2)$, if $\tau_X^\infty=S_R(D_{B,R}(X))$, then
\[
\sup_{s\geq a_0/8}\Prob_s(|\tau_X^\infty-s|>h)\leq2e^{-c'\lambda h^2},\qquad
\sup_{s<a_0/8}\Prob_s(\tau_X^\infty>a_0/8+h)\leq e^{-c'\lambda h^2}.
\]
\end{lemma}
\begin{proof}
Set $\Psi(y)=\int_0^y\pbar(t)\dd t$ and $L(y)=\Psi(y)-y\pbar(y)$. For a nonincreasing integrable profile, $L$ is nondecreasing and tends to $\bar I/\omega_d>0$, since $y\pbar(y)\to0$. Choose $y_1$ such that $L(y_1)\geq\bar I/(2\omega_d)$ and $R_0=(16y_1/a_0)^{1/d}$. The function
\[
\Phi_R(b)=\int_{|z|\leq R}\pbar(b|z|^d)\dd z=\frac{\omega_d}{b}\Psi(bR^d)
\]
is absolutely continuous on compact subsets of $(0,\infty)$, and $\Phi_R'(b)=-\omega_d L(bR^d)/b^2$ almost everywhere. For $s\geq a_0/16$ and $t\in[1/2,1]$, $a_0/16\leq g(s,t)\leq1$, so $\Phi_R'(g)\leq-\bar I/2$, while $\partial_s g\geq\min(\gamma,1-\gamma)a_0/16$. For other $t$, the contribution is nonincreasing. Integration gives the stated slope with $c_*=\bar I a_0\min(\gamma,1-\gamma)/64$. Passing under the integral is justified by bounding difference quotients by the integrable function $Ct^{-\gamma}$; no differentiability of the profiles is assumed.

Write $\mu_R(t)=\rho_{\mathrm{cal}}\lambda m_R(t)$ and $M_0=\bar I C_\gamma(a_0/16)^{-\gamma}$. On $[a_0/16,1]$, $\mu_R\leq \rho_{\mathrm{cal}}\lambda M_0$. For $s\geq a_0/8$, comparison with $s+h/2$ and $s-h/2$ follows the preceding proof: half-gaps between means are at least $\rho_{\mathrm{cal}}\lambda c_*h/4$, and each tail probability is at most $e^{-c'\lambda h^2}$ for $c'=\rho_{\mathrm{cal}}c_*^2/(32M_0)$. The condition on $\lambda$ makes this term smaller than $h/2$, allowing the passage from survival probability to rank. If $s+h>1$, overestimation is impossible. For $s<a_0/8$, degree is stochastically larger than at age $a_0/8$, and rank stochastically smaller; the same comparison with $a_0/8+h/2$ proves the last bound without bounding mean degree at very small ages.
\end{proof}

Call a vertex \emph{locally $h$-reliable} if $u_X\geq a_0/8$ and $|\tau_X^\infty-u_X|\leq h/2$, or if $u_X<a_0/8$ and $\tau_X^\infty\leq a_0/8+h/2$. This event depends only on calibration in $B(x,R)$, and its failure probability is at most $2e^{-c'\lambda h^2/4}$, uniformly in age, for sufficiently large $\lambda$. On
\[
\Omega_n'=\{\sup_j|S_{R,n}(j)-S_R(j)|\leq h/2\},
\]
a reliable vertex at least distance $R$ from the boundary satisfies $|\widetilde u_X-u_X|\leq h$ if $u_X\geq a_0/8$; otherwise its rank is below $a_0/2$ and it is selected in no slice. Lemma~\ref{lem:rank-vote} gives $\Prob(\Omega_n')\to1$.

\subsection{Separating pairs and certifying cells}\label{app:universal-local}
Fix a slice $I$ of length $2w$ and $h\leq w/32$. Write $I^{+h}$ for its enlargement by $h$, intersected with $(0,1]$, and $I^{-h}$ for the interval obtained by removing $h$ at each endpoint. Thus $|I^{+h}|\leq2w+2h$ and $|I^{-h}|=2w-2h$; this convention also handles the last slice, whose upper endpoint is 1. Endpoint choices almost surely do not change true-age counts.

Under Palm at two distinct vertices $X,Y\in A$, define the ideal count $C^{\rm id}(X,Y)$ as $C_i(X,Y)$, but selecting other vertices by their true age in $I$. It is Poisson with mean $p\lambda M_\nu(\rho,u,v)$, where $\rho=|x-y|$ and $\nu=0,1$ denotes opposite or equal signs, respectively. Explicitly, set $S_{xy}=B(x,R)\cap B(y,R)$ and, in the integral over $(z,t)$, $f_X=\phin(g(u,t)|z-x|^d)$, $h_X=\phout(g(u,t)|z-x|^d)$, with analogous notation $f_Y,h_Y$. Then
\[
M_1=\frac12\int_I\int_{S_{xy}}(f_Xf_Y+h_Xh_Y)\dd z\dd t,
\qquad M_0=\frac12\int_I\int_{S_{xy}}(f_Xh_Y+h_Xf_Y)\dd z\dd t.
\]
Both Palm vertices are excluded from all auxiliary counts below, preserving the Poisson laws of the other points. On $\mathcal D=\{0<\rho\leq3R_c/4,\ u,v\in I^{+h}\}$, set
\[
D_* =\inf_{\mathcal D}(M_1-M_0),\qquad
V=\max_{\nu=0,1}(\sup_{\mathcal D}M_\nu-\inf_{\mathcal D}M_\nu),\qquad
U=\sup_{\mathcal D}M_1.
\]

\begin{lemma}[Uniform margin for local thresholds]\label{lem:universal-margin}
There exist $c_0,C>0$ such that, for sufficiently small $w,R_c$ and sufficiently large $R$, uniformly along the chain,
\[
D_*\geq c_0w,\qquad U\leq Cw,\qquad
V\leq Cw\{w+\omega(R_c)+T(R/2)\},
\]
where $\omega(r)\to0$ as $r\downarrow0$ and $T(r)\to0$ as $r\to\infty$. In particular, the scales in \eqref{eq:universal-scales} ensure $V\leq D_*/8$ for sufficiently large $\lambda$. With
\begin{equation}\label{eq:quartile-brackets}
b=\frac{p\lambda D_*}{32},\quad
L_\nu=p\lambda\inf_{\mathcal D}M_\nu-2b,\quad U_\nu=p\lambda\sup_{\mathcal D}M_\nu+2b,\quad
T_- =\frac{L_0+L_1}2,\quad T_+=\frac{U_0+U_1}2,
\end{equation}
we have $U_0<L_1$ and
\[
p\lambda\inf M_1-T_+\geq5p\lambda D_*/16,\qquad
T_--p\lambda\sup M_0\geq5p\lambda D_*/16.
\]
\end{lemma}
\begin{proof}
All ages involved are at least $a_0/4$. For $f\in\{\phin,\phout\}$ and $b\in[a_0/4,1]$, write $f_b(z)=f(b|z|^d)$. These functions are integrable, and their translations are uniformly continuous in $L^1$: the change of variables $z\mapsto b^{1/d}z$ gives
\[
\sup_b\sup_{|y|\leq r}\|f_b(\cdot+y)-f_b\|_1\leq\frac4{a_0}\sup_{|y|\leq r}\|f_1(\cdot+y)-f_1\|_1.
\]
Take $\omega$ to be the maximum of these moduli for the two profiles. Similarly, $T(r)=\sup_{b,f}\int_{|z|\geq r}f_b(z)\dd z\to0$. The bound $U\leq4(I_{\rm in}+I_{\rm out})w/a_0$ follows by bounding each product by its first factor.

For variation with age, monotonicity gives exactly
$\|f_b-f_{b'}\|_1=I_f|b^{-1}-{b'}^{-1}|$, while $|\partial_u g(u,t)^{-1}|\leq16/a_0^2$ on the compact set considered. Both endpoint ages vary over an interval of length at most $4w$, so integration over $t$ costs $Cw^2$. Variation in $\rho$ costs $Cw\omega(R_c)$ for the profiles and $CwT(R/2)$ for the change in domain $S_{xy}$, provided $R\geq2R_c$. This gives the bound on $V$.

Finally, Lemma~\ref{lem:signal-interval} gives $[t_-,t_+]$ and $b_0>0$ on which $\psi\geq b_0$. As $w\to0$, ratios between values of $g$ on each enlarged slice tend uniformly to one, since ages stay above $a_0/4$. There is therefore an annular region of strictly positive width, uniform along the chain, such that both arguments of $\psi$ remain in $[t_-,t_+]$ for small enough $\rho$. Lemma~\ref{lem:annulus} places within it a ball of volume at least $v_0>0$. Its radii are uniformly bounded, so it lies in $S_{xy}$ for large enough $R$. Hence
$M_1-M_0=\frac12\int_I\int_{S_{xy}}\psi\psi\geq wb_0^2v_0$.

These choices give $V/D_*\leq1/8$. We have $\inf M_1-\sup M_0\geq D_*-V$ and
$\inf M_1-(\sup M_0+\sup M_1)/2\geq D_*/2-V$.
After subtracting $2b/(p\lambda)=D_*/16$, the margin is at least $5D_*/16$. The opposite calculation is identical and also gives $L_1-U_0>0$.
\end{proof}

\begin{lemma}[Location of empirical quartiles]\label{lem:universal-quartiles}
One may choose a model-dependent constant $c_h\in(0,1/32]$ such that, with $h=c_hw$, for sufficiently large $\lambda$,
\[
\Prob\bigl(\widehat t_i\in[T_-^{(i)},T_+^{(i)}]\text{ for all slices }i\bigr)\longrightarrow1.
\]
The constant $c_h$ is used only in the proof; these choices ensure $R\geq R_0$, $h\leq a_0/16$ and $\lambda h^2\gg\log(1/h)$ for sufficiently large $\lambda$, as required by Lemma~\ref{lem:local-rank}.
\end{lemma}
\begin{proof}
For a fixed slice, let $\mathcal Q^\pm$ be the ordered pairs within distance $3R_c/4$ whose true ages lie in $I^{\pm h}$. Remove the set $\mathcal B$ of pairs with an unreliable endpoint or an endpoint within distance $2R$ of the boundary. On $\Omega_n'$, the rank bounds give
\begin{equation}\label{eq:pair-sandwich}
\mathcal Q^-\setminus\mathcal B\subset\mathcal P_i\subset\mathcal Q^+\cup\mathcal B.
\end{equation}
By Mecke and ergodicity, $|\mathcal Q^\pm|/n\to Q^\pm=(p\lambda)^2|I^{\pm h}|^2\omega_d(3R_c/4)^d$, and pairs of each sign type have half this density. The density of $\mathcal B$ is at most
$C(p\lambda)^2R_c^d e^{-c'\lambda h^2/4}$, up to an $o_{\mathbb P}(1)$ remainder; the boundary part is negligible. This two-vertex Palm estimate in $A$ is valid because inserting these vertices does not change the calibration population $B$. Since $|I^{-h}|/|I^{+h}|\geq31/33$ and $\lambda h^2\gg\log(1/w)$, \eqref{eq:pair-sandwich} implies that each sign type occupies between $3/8$ and $5/8$ of selected pairs, with probability tending to one.

For a selected pair outside $\mathcal B$, let $N_X^{\rm bd}$ count neighbours $Z\in A\cap B(x,R)\setminus\{X,Y\}$ with age in $I^{+h}\setminus I^{-h}$, and let $N_X^{\rm nf}$ count unreliable neighbours in the same ball. Then, on $\Omega_n'$,
\begin{equation}\label{eq:local-count-perturbation}
|C_i(X,Y)-C^{\rm id}(X,Y)|\leq N_X^{\rm bd}+N_X^{\rm nf}.
\end{equation}
Restriction to $B(x,R)$ matters: it will also make certificates local. Under Palm, the first count is Poisson with mean at most $16p\lambda\bar I h/a_0$. Choose
$c_h\leq\min\{1/32,e^{-2}c_0a_0/(1024\bar I)\}$; this mean is then at most $e^{-2}b/2$. The Poisson bounds and $D_*\geq c_0w$, $U\leq Cw$ give
\[
\Prob(|C^{\rm id}-p\lambda M_\nu|>b)\leq2e^{-c\lambda w},\qquad
\Prob(N_X^{\rm bd}\geq b/2)\leq e^{-c\lambda w}.
\]
Reliability of $Z$ uses only its calibration edges. Conditional on ages and signs, it is independent of the edge $X\sim Z$; its probability is independent of the sign of $Z$. Mecke therefore gives $\E N_X^{\rm nf}\leq C\lambda e^{-c'\lambda h^2/4}$ for $u_X\geq a_0/4$, then
\[
\Prob(N_X^{\rm nf}\geq b/2)\leq\frac Cw e^{-c'\lambda h^2/4}.
\]
The probability of the union of the three failures is thus at most
$\pi_\lambda=Ce^{-c\lambda w}+Cw^{-1}e^{-c'\lambda h^2/4}\to0$, uniformly on $\mathcal D$. Outside these failures, $C_i\in[L_\nu,U_\nu]$ by \eqref{eq:local-count-perturbation}. Ergodicity of pair scores and Mecke therefore show that the fraction of selected pairs of a given type falling outside its interval is at most
$(33/31)^2\pi_\lambda+o_\lambda(1)+o_{\mathbb P}(1)<1/16$.
Here $o_\lambda(1)$ comes from unreliable endpoints, not from independence between counts.

To conclude, consider the counting measure $\mu=\mu_0+\mu_1$ of the values of $C_i$, separated by true signs. Each $\mu_\nu$ carries between $3/8$ and $5/8$ of the total mass, with at most $1/16$ of that mass outside $[L_\nu,U_\nu]$. Since $U_0<L_1$, the mass below $L_0$ is at most $1/8$, and the mass at or below $U_0$ is at least $3/8-1/16>1/4$. Similarly, the mass strictly below $L_1$ is at most $5/8+1/16<3/4$, and the mass at or below $U_1$ is at least $7/8$. Thus the first quartile lies in $[L_0,U_0]$ and the third in $[L_1,U_1]$, locating their midpoint. A finite union bound gives the conclusion for the whole chain.
\end{proof}

\begin{lemma}[Uniform local certification]\label{lem:universal-cells}
There is a stationary field of certified cells, with dependence range $M\leq C_d(R/\ell+1)$, satisfying \eqref{eq:universal-defect}. On $\Omega_n'$ and the event in Lemma~\ref{lem:universal-quartiles}, every certified cell at least distance $3R$ from the boundary is T-Good for the observable decisions.
\end{lemma}
\begin{proof}
A cell $z$ is certified if the following three conditions hold in the infinite graph:
\begin{enumerate}
\item It contains at least $p\lambda w\ell^d$ vertices of $A$ with true ages in $I^{-h}$.
\item Every vertex of $A$ within distance $R$ of the block $L_1(z)$ is locally $h$-reliable.
\item For every distinct pair in $L_1(z)$ with ages in $I^{+h}$ and type $\nu$, we have $|C^{\rm id}-p\lambda M_\nu|\leq b$ and $N_X^{\rm bd}<b/2$.
\end{enumerate}
The first condition ensures observed occupancy, since $\Lambda\leq\lambda$ and, for large $\lambda$, $p\lambda w\ell^d\geq1$. The second ensures that all selected vertices in the block have true ages in $I^{+h}$ and makes $N_X^{\rm nf}$ zero in \eqref{eq:local-count-perturbation}. The count error is therefore at most $3b/2$, strictly below the margin $5p\lambda D_*/16$ in Lemma~\ref{lem:universal-margin}. Every selected pair classification in the block is correct, for any threshold in $[T_-,T_+]$. Products around cycles are then positive, also giving A-Good.

The certificate uses only the induced graph within distance $2R+3\sqrt d\ell$ of the cell: calibration of neighbouring vertices needs an additional radius $R$, and $N_X^{\rm bd}$ was truncated at $R$. Restrictions to disjoint regions are independent, giving the stated range. The first count is Poisson with mean $p\lambda(2w-2h)\ell^d$, and Chernoff bounds its failure by $e^{-c\lambda w\ell^d}$. Mecke and Lemma~\ref{lem:local-rank} bound failure of the second condition by $C\lambda R^d e^{-c\lambda h^2}$. Finally, two-point Mecke applied to the two Poisson tails in the preceding proof bounds the third failure by $C(\lambda\ell^d)^2e^{-c\lambda w}$. Their union gives \eqref{eq:universal-defect}, since $h=c_hw$.

No conditioning on empirical ranks or thresholds was used to establish finite dependence. We constructed a local field that certifies all decisions simultaneously when their global parameters lie in the permitted intervals; this distinction is precisely what allows the percolation step.
\end{proof}

\subsection*{AI statement}
AI has only been used for redactional purposes and checking  correctness of proofs   in Section  \ref{sec:propagation}. The article topic and proof ideas have been provided by authors exclusively, and we take full responsability for the content of the article and the proofs.


\begin{thebibliography}{99}
\bibitem{abbe-sbm}
E.~Abbe,
\emph{Community detection and stochastic block models: recent developments},
Journal of Machine Learning Research \textbf{18}(177) (2018), 1--86.
\href{https://www.jmlr.org/papers/v18/16-480.html}{jmlr.org/papers/v18/16-480.html}.

\bibitem{abbe}
E.~Abbe, F.~Baccelli and A.~Sankararaman,
\emph{Community detection on Euclidean random graphs},
Information and Inference: A Journal of the IMA \textbf{10}(1) (2021),
109--160. \href{https://doi.org/10.1093/imaiai/iaaa009}
{doi:10.1093/imaiai/iaaa009}.
\href{https://arxiv.org/abs/1706.09942}{arXiv:1706.09942}.

\bibitem{alaoglu}
L.~Alaoglu,
\emph{Weak topologies of normed linear spaces},
Annals of Mathematics, Second Series, \textbf{41}(1) (1940), 252--267.

\bibitem{barabasi}
{A.-L. Barabási and R. Albert.
Emergence of scaling in random networks.
\emph{Science}, 286(5439):509--512, 1999.}

\bibitem{bickelchen}
P.~J.~Bickel and A.~Chen,
\emph{A nonparametric view of network models and Newman--Girvan and other
modularities},
Proceedings of the National Academy of Sciences \textbf{106}(50) (2009),
21068--21073.
\href{https://doi.org/10.1073/pnas.0907096106}{doi:10.1073/pnas.0907096106}.

\bibitem{gracar2019}
P.~Gracar, A.~Grauer, L.~L\"uchtrath and P.~M\"orters,
\emph{The age-dependent random connection model},
Queueing Systems \textbf{93} (2019), 309--331.
\href{https://arxiv.org/abs/1810.03429}{arXiv:1810.03429}.

\bibitem{gracar2021}
P.~Gracar, L.~L\"uchtrath and P.~M\"orters,
\emph{Percolation phase transition in weight-dependent random connection models},
Advances in Applied Probability \textbf{53}(4) (2021), 1090--1114.
\href{https://arxiv.org/abs/2003.04040}{arXiv:2003.04040}.

\bibitem{ingstersuslina}
Y.~I.~Ingster and I.~A.~Suslina,
\emph{Nonparametric Goodness-of-Fit Testing Under Gaussian Models},
Lecture Notes in Statistics, vol.~169, Springer, New York, 2003.
\href{https://doi.org/10.1007/978-0-387-21580-8}
{doi:10.1007/978-0-387-21580-8}.

\bibitem{lastpenrose}
G.~Last and M.~Penrose.
\emph{Lectures on the Poisson Process}.
Institute of Mathematical Statistics Textbooks, vol.~7.
Cambridge University Press, 2017.

\bibitem{lss}
T.~M.~Liggett, R.~H.~Schonmann and A.~M.~Stacey,
\emph{Domination by product measures},
The Annals of Probability \textbf{25}(1) (1997), 71--95.
\href{https://doi.org/10.1214/aop/1024404279}{doi:10.1214/aop/1024404279}.

\bibitem{lovasz}
L.~Lov\'asz,
\emph{Large Networks and Graph Limits},
American Mathematical Society Colloquium Publications, vol.~60,
American Mathematical Society, Providence, RI, 2012.

\bibitem{lovaszszegedy}
L.~Lov\'asz and B.~Szegedy,
\emph{Limits of dense graph sequences},
Journal of Combinatorial Theory, Series B \textbf{96}(6) (2006), 933--957.
\href{https://arxiv.org/abs/math/0408173}{arXiv:math/0408173}.

\bibitem{newman}
M.~E.~J.~Newman,
\emph{The structure and function of complex networks},
SIAM Review \textbf{45}(2) (2003), 167--256.
\href{https://doi.org/10.1137/S003614450342480}{doi:10.1137/S003614450342480}.

\bibitem{nguyen1979ergodic}
X.~X. Nguyen and H. Zessin.
\newblock Ergodic theorems for spatial processes.
\newblock \emph{Zeitschrift f{\"u}r Wahrscheinlichkeitstheorie
und verwandte Gebiete}, 48:133--158, 1979.

\bibitem{timar}
Á.~Timár,
\emph{Boundary-connectivity via graph theory},
2011 version, \href{https://arxiv.org/abs/0711.1713}{arXiv:0711.1713}.

\end{thebibliography}
\end{document}